\documentclass{amsart}
\usepackage{amsfonts}
\usepackage{amsmath}
\newtheorem{theorem}{Theorem}
\newtheorem{lemma}[theorem]{Lemma}
\newtheorem{corollary}[theorem]{Corollary} 

\begin{document}
\title[Your comment to wonkpark@hotmail.com]{Biholomorphic mapping on the boundary I}
\maketitle
\author{Won K. Park}
\address{Department of Mathematics, University of Seoul}
\email{wonkpark@hotmail.com}

\begin{abstract}
We present a new proof of Chern-Ji's mapping theorem on a strongly
pseudoconvex domain with differentiable spherical boundary. We show that a
proper holomorphic self mapping of a strongly pseudoconvex domain with the
real analytic boundary is biholomorphic.
\end{abstract}

\addtocounter{section}{-1}

\section{Introduction and Preliminaries}

We shall show that a bounded domain $D$ is biholomorphic to an open ball $%
B^{n+1}$ whenever the boundary $bD$ is locally biholomorphic to the boundary
of an open ball $B^{n+1}.$

\begin{theorem}
\label{first}Let $D$ be a simply connected bounded domain in $\Bbb{C}^{n+1}$
with differentiable spherical boundary $bD.$ Suppose that there is a
biholomorphic mapping 
\begin{equation*}
\phi \in H\left( U\cap D\right) \cap C^{1}\left( U\cap \overline{D}\right) 
\end{equation*}
for a connected open neighborhood $U$ of a point $p\in bD$ satisfying 
\begin{equation*}
\phi \left( U\cap bD\right) \subset bB^{n+1}.
\end{equation*}
Then the mapping $\phi $ is analytically continued to a biholomorphic
mapping from $D$ onto $B^{n+1}.$
\end{theorem}

Our result is a new proof of a weaker version of Chern-Ji's mapping theorem 
\cite{CJ}. The main steps of our proof come as follows: We show that the
inverse mapping $\phi ^{-1}$ is analytically continued on the unit ball $%
B^{n+1}$ to be a locally biholomorphic mapping 
\begin{equation*}
\varphi :B^{n+1}\rightarrow D.
\end{equation*}
We show that the mapping $\varphi $ is a proper holomorphic mapping onto a
universal covering Riemann domain over $D.$ Thus the mapping $\varphi $ is a
biholomorphic mapping whenever $D$ is simply connected.

We shall study on a proper holomorphic mapping $\phi $ between strongly
pseudoconvex bounded domains $D,D^{\prime }$ with real analytic boundaries $%
bD,bD^{\prime }.$

\begin{theorem}
\label{second}Let $D,D^{\prime }$ be strongly pseudoconvex bounded domains
in $\Bbb{C}^{n+1}$ with real analytic boundaries $bD,bD^{\prime }$ and $\phi
:D\rightarrow D^{\prime }$ be a proper holomorphic mapping. Then the mapping 
$\phi $ is locally biholomorphic. If $D=D^{\prime }$, then the mapping $\phi 
$ is a biholomorphic self mapping.
\end{theorem}

Our result is a new proof of a weaker version of Pinchuk's mapping theorem 
\cite{Pi}. The main steps of our proof come as follows: We show that the
mapping $\phi $ is analytically continued along any path on $bD$ as a
locally biholomorphic mapping when $bD$ is nonspherical so that the mapping $%
\phi :D\rightarrow D^{\prime }$ are locally biholomorphic. From the study of
Theorem \ref{first}, we show that the same is true when $bD$ is spherical so
that the mapping $\phi :D\rightarrow D^{\prime }$ are locally biholomorphic.
For the case of $D=D^{\prime },$ we show that the boundary $bD$ is
necessarily spherical whenever the claim is not true. Then we show that
there is a sequence of automorphisms $\phi _{j}\in Aut\left( D\right) $ and
a sequence of points $p_{j}$ on a compact subset $K\subset \subset D$ such
that 
\begin{equation*}
\phi _{j}\left( p_{j}\right) \rightarrow bD
\end{equation*}
whenever the boundary $bD$ is spherical and the claim is not true. We apply
Wong-Rosay Theorem so that the domain $D$ is biholomorphic to an open ball $%
B^{n+1}.$ Then we obtain a contradiction that the mapping $\phi $ induces a
nonautomorphic proper self mapping of an open ball $B^{n+1}$ whenever the
boundary $bD$ is spherical and the claim is not true.

We remark that Theorem \ref{first} is a weaker version of the following
theorem:

\begin{theorem}[cf. Chern-Ji [CJ]]
Let $D$ be a simply connected bounded domain in $\Bbb{C}^{n+1}$ with
continuous spherical boundary $bD.$ Suppose that there is a biholomorphic
mapping 
\begin{equation*}
\phi \in H\left( U\cap D\right) \cap C\left( U\cap \overline{D}\right)
\end{equation*}
for a connected open neighborhood $U$ of a point $p\in bD$ satisfying 
\begin{equation*}
\phi \left( U\cap bD\right) \subset bB^{n+1}.
\end{equation*}
Then the mapping $\phi $ is analytically continued to a biholomorphic
mapping from $D$ onto $B^{n+1}.$
\end{theorem}

We remark that Theorem \ref{second} is a weaker version of the following
theorem:

\begin{theorem}[cf. Pinchuk [Pi]]
Let $D,D^{\prime }$ be strongly pseudoconvex bounded domains in $\Bbb{C}%
^{n+1}$ with the boundaries $bD,bD^{\prime }$ of class $C^{2}$ and $\phi
:D\rightarrow D^{\prime }$ be a proper holomorphic mapping. Then the mapping 
$\phi $ is locally biholomorphic. If $D=D^{\prime }$, then the mapping $\phi 
$ is a biholomorphic self mapping.
\end{theorem}

We presented parts of this article in the spring meeting of Korean
Mathematical Society in 2000. This article is a preliminary version. Any
comment shall be greatly appreciated. Send your comments to
wonkpark@hotmail.com.

\subsection{Canonical normalizing mapping}

Let $M$ be a nondegenerate analytic real hypersurface in $\Bbb{C}^{n+1}$.
For each point $p\in M,$ there is a complex tangent hyperplane $H_{p}\subset
T_{p}M$ so that there is a unit tangent vector $v_{p}\in T_{p}M$
perpendicular to the complex tangent hyperplane $H_{p}$ with respect to the
usual riemannian metric in $\Bbb{C}^{n+1}=\Bbb{R}^{2n+2}.$ Then we can take
a unique distinguished chain $\gamma _{p}$ tangential to the direction $%
v_{p} $ and passing through the complex tangent hyperplane $H_{p}$ at the
point $p$ on $M$(cf. \cite{Pa3}). Further, there is a distinguished normal
parametrization on the chain $\gamma _{p}$ having the same values up to
order $2$ of the straight real line to the direction $v_{p}$ with the usual
euclidean parametrization.. Therefore, we can take a distinguished
normalizing mapping $\mu _{p}$ to Moser normal form 
\begin{equation*}
\mu _{p}:M\rightarrow \mu _{p}\left( M\right)
\end{equation*}
sending the germ $M$ at the point $p$ to a normal form such that $\mu
_{p}(p) $ is the origin and $\mu _{p}\left( \gamma _{p}\right) $ is on the
straightened chain of the normal form.

We take a local orientation near the point $p\in M$ so that the tangent
vector $v_{p}$ extends to a smooth unit vector field $v$ on $M\cap U$ for an
open neighborhood $U$ of the point $p$ such that $v_{q}\in T_{q}M$ is a unit
tangent vector perpendicular to the complex tangent hyperplane $H_{q}$ for
each $q\in M\cap U.$ Then we obtain a family of the distinguished
normalizing mapping $\mu _{q}$ for $q\in M\cap U$ : 
\begin{equation*}
\mu _{q}:M\rightarrow \mu _{q}\left( M\right)
\end{equation*}
associated the unit tangent vector $v_{q}\in T_{q}M\backslash H_{q}.$ The
distinguished normalizing mapping $\mu _{q}$ shall be called the canonical
normalizing mapping associated Moser normal form.

\begin{lemma}
\label{canonicalmapping}Let $p_{j}\in M$ be a sequence of points converging
to a point $p\in M.$ Then there is a positive real number $\delta >0$ such
that

\begin{enumerate}
\item  the mapping $\mu _{p_{j}}$ and its inverse $\mu _{p_{j}}^{-1}$ are
analytically continued respectively on 
\begin{equation*}
B\left( p_{j};\delta \right) \quad \text{and}\quad B\left( 0;\delta \right)
\end{equation*}
as a biholomorphic mapping,

\item  the real hypersurface $\mu _{p_{j}}\left( M\right) $ is analytically
continued on $B\left( 0;\delta \right) $ by its defining equation,

\item  the sequence $\mu _{p_{j}}^{-1}$ uniformly converges to $\mu
_{p}^{-1} $ on $B\left( 0;\delta \right) $ as a biholomorphic mapping,

\item  the sequence $\mu _{p_{j}}\left( M\right) $ uniformly converges on $%
B\left( 0;\delta \right) $ to the real hypersurface $\mu _{p}\left( M\right)
.$
\end{enumerate}
\end{lemma}

Let $H$ be the local automorphism group at the origin of the real
hyperquadric 
\begin{equation*}
v=\langle z,z\rangle .
\end{equation*}
The isotropy subgroup $Aut_{p}\left( M\right) $ is naturally identified to
the isotropy subgroup $Aut_{0}\left( \mu _{p}\left( M\right) \right) $ by
the following relation: 
\begin{equation*}
\mu _{p}\circ \phi \circ \mu _{p}^{-1}\in Aut_{0}\left( \mu _{p}\left(
M\right) \right) \quad \text{for}\quad \phi \in Aut_{p}\left( M\right) .
\end{equation*}
Note that every biholomorphic mapping $\varphi $ between real hypersurfaces
in normal form is faithfully represented by a natural group action of the
isotropy subgroup $H$(cf. \cite{Pa3}) such that 
\begin{equation*}
\varphi =N_{e}\quad \text{for}\quad e\in H.
\end{equation*}
Because $\mu _{p}\left( M\right) $ is in normal form, there is a natural
identification of a local automorphism $\phi \in Aut_{p}\left( M\right) $ to
an element 
\begin{equation*}
\left( U_{\phi },a_{\phi },\rho _{\phi },r_{\phi }\right) \in H,
\end{equation*}
where $U_{\phi },a_{\phi },\rho _{\phi },r_{\phi }$ are the normalizing
parameters(cf. \cite{Pa3}) of the mapping 
\begin{equation*}
\mu _{p}\circ \phi \circ \mu _{p}^{-1}\in Aut_{0}\left( \mu _{p}\left(
M\right) \right) .
\end{equation*}

\begin{lemma}
If the isotropy subgroup $Aut_{p}\left( M\right) $ is compact, then $%
Aut_{p}\left( M\right) $ is isomorphic to the subgroup 
\begin{equation*}
\left\{ \left( U_{\phi },a_{\phi },\rho _{\phi },r_{\phi }\right) \in H:\phi
\in Aut_{p}\left( M\right) \right\}
\end{equation*}
as a Lie group.
\end{lemma}

\subsection{Preliminary Lemmas}

We have lemmas on the automorphism of $B^{n+1}.$

\begin{lemma}
Let $p,q$ be two distinct points on $bB^{n+1}$ and $\phi \in Aut\left(
bB^{n+1}\right) $ be a local automorphism of $bB^{n+1}$ such that 
\begin{equation*}
\phi \left( p\right) \neq q.
\end{equation*}
Then there is a unique decomposition 
\begin{equation*}
\phi =\psi \circ \varphi
\end{equation*}
where 
\begin{equation*}
\varphi \in Aut_{p}\left( bB^{n+1}\right) ,\quad \psi \in Aut_{q}\left(
bB^{n+1}\right)
\end{equation*}
and the local automorphism $\psi $ acts trivially on the complex tangent
hyperplane of $bB^{n+1}$ at the fixed point $q.$
\end{lemma}

\proof%
Note that the isotropy subgroup $Aut_{q}\left( bB^{n+1}\right) $ acts on $%
bB^{n+1}\backslash q$ transitively. Further, there is a unique element $\psi
\in Aut_{q}\left( bB^{n+1}\right) $ for each point displacement on $%
bB^{n+1}\backslash q$ by requiring the element $\psi $ acts trivially on the
complex tangent hyperplane of $bB^{n+1}$ at the point $q$(cf. \cite{Pa1})$.$

Let's put $p^{\prime }=\phi \left( p\right) .$ Since $p^{\prime }\neq q,$ we
take a unique automorphism $\psi \in Aut_{q}\left( bB^{n+1}\right) $ such
that 
\begin{equation*}
\psi \left( p^{\prime }\right) =p.
\end{equation*}
Then $\varphi \equiv \psi \circ \phi \in Aut_{p}\left( bB^{n+1}\right) $ so
that 
\begin{equation*}
\phi =\psi ^{-1}\circ \varphi .
\end{equation*}
This completes the proof.%
\endproof%

\begin{lemma}
\label{decomposition}Let $p,q$ be two distinct points on $bB^{n+1}$ and $%
\phi \in Aut\left( bB^{n+1}\right) $ be a local automorphism of $bB^{n+1}$
such that 
\begin{equation*}
p^{\prime }\equiv \phi \left( p\right) \neq q.
\end{equation*}
Then there is a unique decomposition 
\begin{equation*}
\phi =\varphi \circ \psi
\end{equation*}
where 
\begin{equation*}
\varphi \in Aut_{p^{\prime }}\left( bB^{n+1}\right) ,\quad \psi \in
Aut_{q}\left( bB^{n+1}\right)
\end{equation*}
and the local automorphism $\psi $ acts trivially on the complex tangent
hyperplane of $bB^{n+1}$ at the fixed point $q.$
\end{lemma}

\proof%
By Lemma \ref{decomposition}, there is a decomposition 
\begin{equation*}
\phi ^{-1}=\psi \circ \varphi
\end{equation*}
where 
\begin{equation*}
\varphi \in Aut_{p^{\prime }}\left( bB^{n+1}\right) ,\quad \psi \in
Aut_{q}\left( bB^{n+1}\right)
\end{equation*}
and the local automorphism $\psi $ acts trivially on the complex tangent
hyperplane of $bB^{n+1}$ at the fixed point $q.$ Hence we obtain 
\begin{equation*}
\phi =\varphi ^{-1}\circ \psi ^{-1}
\end{equation*}
where 
\begin{equation*}
\varphi ^{-1}\in Aut_{p^{\prime }}\left( bB^{n+1}\right) ,\quad \psi
^{-1}\in Aut_{q}\left( bB^{n+1}\right) .
\end{equation*}
This completes the proof.%
\endproof%

\begin{lemma}
\label{either}Let $\phi _{j}$ be a sequence of automorphisms of $B^{n+1}$.
Suppose that the sequence $\phi _{j}$ converges to a holomorphic mapping $%
\lambda $ uniformly on every compact subset of $B^{n+1}.$ Then the mapping $%
\lambda $ is either a constant mapping or an automorphism of $B^{n+1}.$
\end{lemma}

\proof%
Note that the mapping $\lambda $ satisfies 
\begin{equation}
\lambda \left( B^{n+1}\right) \subset \overline{B^{n+1}}
\label{boundarypoint}
\end{equation}
and a complex line is mapped to a complex line under the biholomorphic
automorphism of the unit ball $B^{n+1}.$

Suppose that there is a complex line $\pi $ such that 
\begin{equation*}
\pi \cap B^{n+1}\neq \emptyset \quad \text{and}\quad \lambda \left( \pi \cap
B^{n+1}\right) =\emptyset .
\end{equation*}
Note that $bB^{n+1}$ is strongly pseudoconvex so that, by the condition \ref
{boundarypoint}, there is a point $q\in bB^{n+1}$ satisfying 
\begin{equation*}
q=\lambda \left( \pi \right) \cap \overline{B^{n+1}}.
\end{equation*}
Then we obtain 
\begin{equation*}
\left. \lambda \right| _{\pi \cap B^{n+1}}=q.
\end{equation*}
Let $x$ be an interior point of $\pi \cap B^{n+1}$ and $p^{\prime }$ be an
arbitrary point of $B^{n+1}$ so that we take a complex line $\pi ^{\prime }$
passes through $x$ and $p^{\prime }.$ Since $B^{n+1}$ is strongly
pseudoconvex, the maximum modulus theorem of one complex variable yields 
\begin{equation*}
\left. \lambda \right| _{\pi ^{\prime }\cap B^{n+1}}=q
\end{equation*}
so that the mapping $\lambda $ is a constant mapping.

Suppose that the mapping $\lambda $ is not a constant mapping. Then, for a
complex line $\pi $ satisfying 
\begin{equation*}
\pi \cap B^{n+1}\neq \emptyset ,
\end{equation*}
there is a real number $\varepsilon >0$ such that 
\begin{equation}
\left| \phi _{j}\left( \pi \cap B^{n+1}\right) \right| \geq \varepsilon
\label{finitesize}
\end{equation}
where $\left| \phi _{j}\left( \pi \cap B^{n+1}\right) \right| $ is the area
of the analytic disk 
\begin{equation*}
\phi _{j}\left( \pi \cap B^{n+1}\right) .
\end{equation*}
We take a point $p\in \pi \cap bB^{n+1}$ so that 
\begin{equation}
\phi _{j}\left( p\right) \rightarrow p^{\prime }\in bB^{n+1}.
\label{approaching}
\end{equation}
Then we take a point $p^{\prime \prime }\in bB^{n+1}$ such that 
\begin{equation*}
p^{\prime \prime }\notin \overline{\left\{ \phi _{j}\left( p\right) :j\in 
\Bbb{N}^{+}\right\} },
\end{equation*}
if necessary, passing to a subsequence. We have the following decomposition 
\begin{equation*}
\phi _{j}=\varphi _{j}\circ \psi _{j}
\end{equation*}
where 
\begin{equation*}
\varphi _{j}\in Aut_{p^{\prime }}\left( bB^{n+1}\right) \quad \text{and}%
\quad \psi _{j}\in Aut_{p^{\prime \prime }}\left( bB^{n+1}\right)
\end{equation*}
where the automorphisms $\psi _{j}$ act trivially on the complex tangent
hyperplane at the fixed point $p^{\prime \prime }.$ Then we obtain 
\begin{equation*}
U_{\psi _{j}}=id_{n\times n},\quad \rho _{\psi _{j}}=1.
\end{equation*}
By the condition \ref{approaching}, there is a real number $e>0$ such that 
\begin{equation*}
\left| a_{\psi _{j}}\right| \leq e,\quad \left| r_{\psi _{j}}\right| \leq e.
\end{equation*}
By the condition \ref{finitesize}, there is a real number $e>0,$ if
necessary, increasing $e,$ such that 
\begin{equation*}
\left| a_{\varphi _{j}}\right| \leq e.
\end{equation*}
Since the mapping $\lambda $ is not a constant mapping, there is a real
number $e>0,$ if necessary, increasing $e,$ such that 
\begin{equation*}
e^{-1}\leq \left| \rho _{\varphi _{j}}\right| \leq e,\quad \left| r_{\varphi
_{j}}\right| \leq e.
\end{equation*}
Since $bB^{n+1}$ is strongly pseudoconvex, we have 
\begin{equation*}
\left| U_{\varphi _{j}}\right| =1.
\end{equation*}
Then the Jacobian determinant $\det \varphi _{j}^{\prime }$ is uniformly
bounded from the zero on an open neighborhood of the point $p\in bB^{n+1}.$
By Hurwitz theorem, the mapping $\lambda $ is locally biholomorphic and,
further, the mapping $\lambda $ is one-to-one. Hence the mapping $\lambda $
is an automorphism of the unit ball $B^{n+1}.$ This completes the proof.%
\endproof%

We have lemmas on the chain of $bB^{n+1}.$

\begin{lemma}
Let $\gamma :[0,1]\rightarrow bB^{n+1}$ be a chain-segment on $bB^{n+1}$.
Then there is a complex line $\pi $ such that 
\begin{equation*}
\gamma [0,1]\subset \pi \cap bB^{n+1}.
\end{equation*}
\end{lemma}

\proof%
We take a point $p\in \gamma [0,1].$ Then the chain-segment $\mu _{p}\circ
\gamma [0,1]$ on $\mu _{p}\left( bB^{n+1}\right) $ is on a complex line $\pi
^{\prime }$(cf. \cite{Pa3}). Since the canonical normalizing mapping $\mu
_{p}$ of the sphere $bB^{n+1}$ is a fractional linear mapping(cf. \cite{Pa3}%
), we take 
\begin{equation*}
\pi =\mu _{p}^{-1}\left( \pi ^{\prime }\right)
\end{equation*}
so that 
\begin{equation*}
\gamma [0,1]\subset \pi \cap bB^{n+1}.
\end{equation*}
This completes the proof.%
\endproof%

\begin{lemma}
Let $\pi $ be a complex line such that 
\begin{equation*}
\pi \cap B^{n+1}\neq \emptyset .
\end{equation*}
Then the circle $\pi \cap bB^{n+1}$ is a chain on $bB^{n+1}.$
\end{lemma}

\proof%
We take a point $p\in \pi \cap bB^{n+1}.$ Then we obtain 
\begin{equation*}
\mu _{p}\left( \pi \right) \cap \mu _{p}\left( bB^{n+1}\right) =\mu
_{p}\left( \pi \cap bB^{n+1}\right) \neq \emptyset .
\end{equation*}
Since the canonical normalizing mapping $\mu _{p}$ of the sphere $bB^{n+1}$
is a fractional linear mapping(cf. \cite{Pa3}), $\mu _{p}\left( \pi \right) $
is a complex line so that $\mu _{p}\left( \pi \right) \cap \mu _{p}\left(
bB^{n+1}\right) $ is a chain(cf. \cite{Pa3}). Thus the circle 
\begin{equation*}
\pi \cap bB^{n+1}=\mu _{p}^{-1}\left( \mu _{p}\left( \pi \right) \cap \mu
_{p}\left( bB^{n+1}\right) \right)
\end{equation*}
is a chain as well. This completes the proof.%
\endproof%

\begin{lemma}
Let $\gamma $ be a chain passing through a point $p\in bB^{n+1}$ and $\delta
_{\gamma }$ be an analytic disk such that 
\begin{equation*}
\gamma =\pi \cap bB^{n+1}\quad \text{and}\quad \delta _{\gamma }=\pi \cap
B^{n+1}
\end{equation*}
where $\pi $ is a complex line. Let $\theta _{\gamma }$ be the angle between
the tangent vector of $\gamma $ at the point $p$ and a unit vector $v_{p}$
perpendicular to the complex tangent hyperplane at the point $p\in bB^{n+1}$
and $\left| \delta _{\gamma }\right| $ be the area of the analytic disk $%
\delta _{\gamma }$. Then 
\begin{equation*}
\left| a_{\gamma }\right| \equiv \left| \tan \theta _{\gamma }\right|
\rightarrow \infty
\end{equation*}
if and only if 
\begin{equation*}
\left| \delta _{\gamma }\right| \rightarrow 0.
\end{equation*}
\end{lemma}

\proof%
We easily see that 
\begin{equation*}
\theta _{\gamma }\rightarrow \frac{\pi }{2}
\end{equation*}
if and only if 
\begin{equation*}
\left| a_{\gamma }\right| \equiv \left| \tan \theta _{\gamma }\right|
\rightarrow \infty .
\end{equation*}
Note that the complex line $\pi $ would be on the complex tangent hyperplane
of $bB^{n+1}$ at the point $p$ if 
\begin{equation*}
\theta _{\gamma }=\frac{\pi }{2}.
\end{equation*}
Since $bB^{n+1}$ is strongly pseudoconvex and $B^{n+1}$ is strongly convex,
the complex tangent hyperplane of $bB^{n+1}$ at the point $p$ has no
intersection to $B^{n+1}.$ Thus we easily see 
\begin{equation*}
\theta _{\gamma }\rightarrow \frac{\pi }{2}
\end{equation*}
if and only if 
\begin{equation*}
\left| \delta _{\gamma }\right| \rightarrow 0.
\end{equation*}
This completes the proof.%
\endproof%

We may require the following well-known results in this article(cf. \cite{Kr}%
, \cite{Ra}, \cite{Bo}).

\begin{lemma}[Lewy, Pinchuk]
\label{Lewy-Pinchuk}Let $D,D^{\prime }$ be domains with strongly
pseudoconvex real analytic boundaries $bD,bD^{\prime }$ and $U$ be a
connected open neighborhood of a point $p\in bD.$ Suppose that there is a
holomorphic mapping $\phi $ on $U\cap D$ such that 
\begin{equation*}
\phi \in H\left( U\cap D\right) \cap C^{1}\left( U\cap \overline{D}\right)
,\quad \phi \left( U\cap bD\right) \subset bD^{\prime }
\end{equation*}
and the induced mapping $\phi :U\cap bD\rightarrow bD^{\prime }$ is CR
diffeomorphic. Then the mapping $\phi $ is analytically continued on $U,$ if
necessary, shrinking $U.$
\end{lemma}

\begin{lemma}[Lewy]
\label{localhull}Let $D$ be a domain with a strongly pseudoconvex boundary $%
bD$ and $U$ be an open connected neighborhood of a point $p\in bD.$ Suppose
that there is a holomorphic mapping $\phi $ on $U\cap bD.$ Then there is an
open neighborhood $V$ of the point $p$ such that the mapping $\phi $ is
analytically continued onto 
\begin{equation*}
V\cap \overline{D}.
\end{equation*}
\end{lemma}

\begin{lemma}[Wong, Rosay]
\label{Wong-Rosay}Let $D$ be a strongly pseudoconvex bounded domain. Suppose
that there is a compact set $K\subset \subset D$ and a sequence $p_{j}\in K$
and automorphisms $\phi _{j}\in Aut\left( D\right) $ such that 
\begin{equation*}
\phi _{j}\left( p_{j}\right) \rightarrow bD.
\end{equation*}
Then the domain $D$ is biholomorphic to an open unit ball $B^{n+1}.$
\end{lemma}

\begin{lemma}[Bell-Catlin, Diederich-Fornaess]
\label{b-regularity}Let $D,D^{\prime }$ be strongly pseudoconvex bounded
domains with the boundaries $bD,bD^{\prime }$ of class $C^{\infty }$ and $%
\phi :D\rightarrow D^{\prime }$ be a proper holomorphic mapping. Then $\phi
:D\rightarrow D^{\prime }$ is a locally biholomorphic mapping and the
induced mapping $\phi :bD\rightarrow bD^{\prime }$ is a locally CR
diffeomorphism.
\end{lemma}

\section{Analytic Continuation on a Sphere}

\subsection{Analytic continuation with finiteness}

Let $D$ be a domain in $\Bbb{C}^{n+1},n\geq 1,$ with real analytic boundary $%
bD.$ The boundary $bD$ shall be called spherical if, for each point $p\in
bD, $ there is a connected open neighborhood $U$ of the point $p$ and a
biholomorphic mapping $\phi $ on $U$ such that 
\begin{equation*}
\phi \left( U\cap bD\right) \subset bB^{n+1}.
\end{equation*}
Note that a domain $D$ with spherical real analytic boundary is necessarily
strongly pseudoconvex.

\begin{lemma}
\label{sphere}Let $p$ be a point of $bB^{n+1}$ and $U$ be a connected open
neighborhood of the point $p$. Suppose that there is a biholomorphic mapping 
$\phi $ on $U$ such that 
\begin{equation}
\phi \left( U\cap bB^{n+1}\right) \subset bB^{n+1}.  \label{sph}
\end{equation}
Then the mapping $\phi $ is analytically continued on an open neighborhood
of the closed ball $\overline{B^{n+1}}.$
\end{lemma}

\proof%
Note that each local automorphism $\varphi \in Aut_{p}\left( bB^{n+1}\right) 
$ for any point $p\in bB^{n+1}$ is necessarily birational such that $\varphi 
$ is analytically continued on an open neighborhood of $\overline{B^{n+1}}$
as a biholomorphic mapping(cf.\cite{Pa1})$.$

Let's put $q=\phi \left( p\right) .$ We take a point $r\in bB^{n+1}$ such
that $r\neq p,$ $r\neq q.$ Then we take an automorphism 
\begin{equation*}
\psi \in Aut_{r}\left( bB^{n+1}\right)
\end{equation*}
satisfying 
\begin{equation*}
\psi \left( q\right) =p
\end{equation*}
so that 
\begin{equation*}
\varphi \equiv \psi \circ \phi \in Aut_{p}\left( bB^{n+1}\right) .
\end{equation*}
Then the mapping $\psi ^{-1}\circ \varphi $ is an automorphism of $B^{n+1}$
and an analytic continuation of the mapping $\phi $ such that 
\begin{equation*}
\phi =\psi ^{-1}\circ \varphi \quad \text{on}\quad U.
\end{equation*}
This completes the proof.%
\endproof%

\begin{theorem}
\label{anypath}Let $D$ be a domain in $\Bbb{C}^{n+1}$ with spherical real
analytic boundary $bD.$ Suppose that there is a connected open neighborhood $%
U$ of a point $p\in bD$ and a biholomorphic mapping $\phi $ on $U$ such that 
$\phi \left( U\cap bD\right) \subset bB^{n+1}.$ Then the mapping $\phi $ is
analytically continued along any path on $bD$ as a local biholomorphic
mapping.
\end{theorem}

\proof%
Suppose that the assertion is not true. Then there would be a path $\gamma
:[0,1]\rightarrow bD$ such that $\gamma \left( 0\right) \in U\cap bD$ and
the germ of a biholomorphic mapping $\phi $ at the point $\gamma \left(
0\right) $ is analytically continued along the subpath $\gamma [0,\tau ]$
with all $\tau <1$ as a local biholomorphic mapping, but not the whole path $%
\gamma [0,1].$

Since $bD$ is spherical, by definition, there exist a connected open
neighborhood $V$ of the point $\gamma (1)$ and a biholomorphic mapping $%
\varphi $ on $V$ such that 
\begin{equation*}
\varphi (V\cap bD)\subset bB^{n+1}.
\end{equation*}
We take $\lambda \in [0,1)$ such that 
\begin{equation*}
\gamma (\tau )\in V\cap bD\quad \text{for all }\tau \in [\lambda ,1]
\end{equation*}
and we take a sufficiently small connected open neighborhood $W$ of the
point $\gamma \left( \lambda \right) $ such that $\phi $ is analytically
continued on $W\subset V$ along the path $\gamma [0,\lambda ]$ and 
\begin{equation*}
\varphi \left( W\right) \cap bB^{n+1}\neq \emptyset .
\end{equation*}
Then we have 
\begin{equation*}
\psi \left( \varphi \left( W\right) \cap bB^{n+1}\right) \subset bB^{n+1}
\end{equation*}
where 
\begin{equation*}
\psi =\phi \circ \varphi ^{-1}.
\end{equation*}
By Lemma \ref{sphere}, $\psi $ is analytically continued on an open
neighborhood of $bB^{n+1}$ as a local biholomorphic mapping. By abuse of
notation, the mapping $\psi \circ \varphi $ is biholomorphic on the open set 
$V$ such that 
\begin{equation*}
\psi \circ \varphi =\phi \quad \text{on }W.
\end{equation*}
Thus the germ $\psi \circ \varphi $ at the point $\gamma (1)$ is an analytic
continuation of the germ $\phi $ at the point $\gamma \left( 0\right) $
along the path $\gamma [0,1]$ as a local biholomorphic mapping. This
contradiction completes the proof.%
\endproof%

\begin{lemma}
\label{boundary}Let $D$ be a bounded domain in $\Bbb{C}^{n+1}$ with
spherical real analytic boundary $bD$ such that the fundamental group $\pi
_{1}\left( bD\right) $ is finite. Suppose that there is a connected open
neighborhood $U$ of a point $p\in bD$ and a biholomorphic mapping $\phi $ on 
$U$ such that $\phi \left( U\cap bD\right) \subset bB^{n+1}.$ Then $\phi $
is analytically continued to a biholomorphic mapping from $D$ onto $B^{n+1}.$
\end{lemma}

\proof%
By Lemma \ref{anypath}, the mapping $\phi $ is analytically continued along
any path on $bD$ as a local biholomorphic mapping. Let $E$ be the path space
of $bD$ pointed at the point $p$ mod homotopy so that $E$ is a universal
covering of $bD$ with a natural CR structure and a CR projection $\varphi
:E\rightarrow bD.$ Then there is a unique CR lift $\psi :E\rightarrow
bB^{n+1}$ as the analytic continuation of the biholomorphic mapping $\phi $.
Note that $\psi :E\rightarrow bB^{n+1}$ is an open mapping because $\phi $
and $\varphi $ are both locally CR diffeomorphisms.

Since $bD$ is finitely connected, $E$ is necessarily compact so that the
mapping $\psi :E\rightarrow bB^{n+1}$ is surjective. Further, the mapping $%
\psi :E\rightarrow bB^{n+1}$ is a simple covering map because $\psi $ is
locally a CR diffeomorphism and the sphere $bB^{n+1}$ is simply connected.
Hence there exists a locally biholomorphic mapping $\lambda
:bB^{n+1}\rightarrow bD$ defined by 
\begin{equation*}
\lambda =\varphi \circ \psi ^{-1}\quad \text{on}\quad bB^{n+1}.
\end{equation*}

By Hartogs extension theorem, the mapping $\lambda $ uniquely extends to the
open ball $B^{n+1}$ as a local biholomorphic mapping and, further, the
extension is smooth up to the boundary. Hence $\lambda $ is well defined as
a locally biholomorphic mapping on an open neighborhood of $\overline{B^{n+1}%
}$ by Lemma \ref{Lewy-Pinchuk}.

We obtain a proper mapping $\lambda :B^{n+1}\rightarrow D$ so that $\lambda $
is a globally branched covering and the branched locus of $\lambda $ cannot
be bounded by $bD.$ Since $\lambda ^{-1}=\phi $ is locally biholomorphic on $%
bD$, $\lambda :B^{n+1}\rightarrow D$ and $\lambda :bB^{n+1}\rightarrow bD$
are finite coverings respectively of $D$ and $bD.$ Since the closed ball $%
\overline{B^{n+1}}$ has the fixed point property, the mapping $\lambda
:B^{n+1}\rightarrow D$ is globally one-to-one. Otherwise, there would be a
nontrivial deck transform of $\overline{B^{n+1}}$ which is continuous on $%
\overline{B^{n+1}}$ without a fixed point. Hence the mapping $\lambda
^{-1}:D\rightarrow B^{n+1}$ is biholomorphic with $\lambda ^{-1}=\phi $ on $%
bD.$ This completes the proof.%
\endproof%

\begin{theorem}
\label{germs}Let $D$ be a bounded domain in $\Bbb{C}^{n+1}$ with a connected
spherical real analytic boundary $bD.$ Suppose that there is a biholomorphic
mapping $\phi $ on a connected open neighborhood $U$ of a point $p\in bD$
satisfying 
\begin{equation*}
\phi \left( U\cap bD\right) \subset bB^{n+1}
\end{equation*}
such that the analytic continuation of $\phi $ on the boundary $bD$ yields
finitely many germs at each point on $bD.$ Then $D$ is necessarily simply
connected and the mapping $\phi $ is analytically continued to a
biholomorphic mapping from $D$ onto $B^{n+1}.$
\end{theorem}

\proof%
We claim that there is a finite covering space $E_{1}$ of $bD$ with a
natural CR structure and a CR projection $\varphi _{1}:E_{1}\rightarrow bD$
and a local CR diffeomorphism $\psi _{1}:E_{1}\rightarrow bB^{n+1}$
satisfying the relation $\psi _{1}=\phi \circ \varphi _{1}.$ Then the
desired result follows from this claim by the same argument in the proof of
Lemma \ref{boundary}.

Let $E$ be the path space of $bD$ pointed at the point $p\in bD$ mod
homotopy so that $E$ is a universal covering of $bD$ with a natural CR
structure and a CR projection $\varphi :E\rightarrow bD.$ Then there is a
unique CR lift $\psi :E\rightarrow bB^{n+1}:$ 
\begin{equation*}
\begin{array}{lll}
E &  &  \\ 
\downarrow \varphi & \overset{\psi }{\searrow } &  \\ 
bD & \overset{\phi }{\longrightarrow } & bB^{n+1}
\end{array}
\end{equation*}
satisfying the relation $\psi =\phi \circ \varphi $. Note that $\psi
:E\rightarrow bB^{n+1}$ is an open mapping because $\phi $ and $\varphi $
are both local CR diffeomorphisms. Let $F$ be the image of the mapping $\psi 
$ such that $F=\psi \left( E\right) .$ Then $F$ is an open subset of $bD.$

Suppose that $bF\neq \emptyset .$ Then we take a point $p\in bF$ and a
sequence $p_{j}\in F$ such that $p_{j}\rightarrow p.$ Thus there exists a
point $q_{j}\in bD$ and germs of biholomorphic mappings $\phi _{j}$ such
that 
\begin{equation*}
\phi _{j}\left( q_{j}\right) \rightarrow p
\end{equation*}
where $\phi _{j}$ are analytic continuations of the mapping $\phi $ on $bD.$
Since $bD$ is compact, there exist a point $q\in bD$ and a subsequence $%
q_{m_{j}}$ of $q_{j}$ such that $q_{m_{j}}\rightarrow q.$ Further, by
passing to a subsequence, if necessary, we may assume that 
\begin{equation*}
\phi ^{*}=\phi _{m_{j}}\quad \text{for all }m_{j}
\end{equation*}
because the analytic continuation of the mapping $\phi $ yields only
finitely many germs at the point $q\in bD.$ Hence we obtain 
\begin{equation*}
\phi ^{*}\left( q\right) =\lim_{j\rightarrow \infty }\phi _{m_{j}}\left(
q_{m_{j}}\right) =p\in F.
\end{equation*}
This contradiction implies that $bF=\emptyset ,$ i.e., $\psi \left( E\right)
=bB^{n+1}.$

For each point $q\in bB^{n+1},$ there is a subset $X_{q}\subset E$ such that 
\begin{equation*}
X_{q}=\left\{ p\in E:\psi \left( p\right) =q\right\} .
\end{equation*}
Since $\psi $ is a CR diffeomorphism, $X_{q}$ is necessarily a discrete set
on $E.$ Then we define a subset $Y_{q}\subset bD$ such that 
\begin{equation*}
Y_{q}=\left\{ \varphi \left( p\right) \in bD:p\in X_{q}\right\} .
\end{equation*}
Suppose that $Y_{q}$ has an accumulation point $y\in bD.$ Then there is a
sequence of point $p_{j}\in bD$ satisfying 
\begin{equation*}
p_{j}\rightarrow y\quad \text{and}\quad p_{j}\neq y,
\end{equation*}
and biholomorphic mappings $\phi _{j}$ such that 
\begin{equation*}
\phi _{j}\left( p_{j}\right) =q
\end{equation*}
where the mapping $\phi $ at the point $p_{j}$ is the analytic continuation
of the mapping $\phi .$ Because the analytic continuation of $\phi $ yields
only finitely many germs at the point $y,$ we can take a subsequence of $%
\phi $ and a biholomorphic mapping $\phi ^{*}$ such that $\phi ^{*}=\phi
_{m_{j}}.$ Then we have 
\begin{equation*}
\phi ^{*}\left( y\right) =\lim_{j\rightarrow \infty }\phi _{m_{j}}\left(
p_{m_{j}}\right) =q.
\end{equation*}
Since $\phi ^{*}$ is locally biholomorphic, it is impossible that $\phi
^{*}\left( p_{m_{j}}\right) =q=\phi ^{*}\left( y\right) $ and $%
p_{m_{j}}\rightarrow y,p_{m_{j}}\neq y$ at the same time. Thus we find that
the set $Y_{q}$ is finite.

Therefore, the analytic continuation of the mapping $\phi $ on $bD$ is
mapped to a point of $bB^{n+1}$ only at finitely many points of $bD.$ Since
the analytic continuation of the mapping $\phi $ yields finitely many germs
at each point on $bD,$ only finitely many germs of the analytic continuation
of the mapping $\phi $ on $bD$ are mapped to each point of $bB^{n+1}.$ Then,
by the compactness of $bB^{n+1},$ we obtain a finite covering space $E_{1}$
of $bD$ satisfying all conditions in the claim. This completes the proof.%
\endproof%

\subsection{First Dogginal Lemma}

\begin{lemma}[First Scaling Lemma]
\label{scaling}Let $p$ be a point of the boundary $bB^{n+1}$ and $p_{j},$ $%
j\in \Bbb{N}^{+},$ be a sequence of points of $bB^{n+1}$ such that $%
p_{j}\neq p$ for all $j$ and $p_{j}\rightarrow p$ as $j\rightarrow \infty $
to a direction transversal to the complex tangent hyperplane at the point $%
p\in bB^{n+1}.$ Let $\varepsilon _{j}$ be the euclidean distance between the
two points $p_{j}$ and $p,$ and $\delta _{j}$ be the analytic disk 
\begin{equation*}
\delta _{j}=\pi _{j}\cap B^{n+1}
\end{equation*}
where $\pi _{j}$ is the complex line passing through the two points $p_{j}$
and $p.$ Suppose that there is a sequence $p_{j}^{\prime }$ of points of $%
bB^{n+1}$ satisfying 
\begin{equation*}
p_{j}^{\prime }\rightarrow p^{\prime }\in bB^{n+1},
\end{equation*}
and a sequence of biholomorphic automorphisms $\phi _{j}\in Aut\left(
B^{n+1}\right) $ satisfying 
\begin{equation*}
\phi _{j}\left( p_{j}^{\prime }\right) =p_{j}
\end{equation*}
such that the sequence $\phi _{j}$ converges to a constant mapping and the
area $\left| \phi _{j}^{-1}\left( \delta _{j}\right) \right| $ of the
analytic disks 
\begin{equation*}
\phi _{j}^{-1}\left( \delta _{j}\right)
\end{equation*}
is bounded from the below, i.e., there is a real number $c>0$ satisfying 
\begin{equation*}
\left| \phi _{j}^{-1}\left( \delta _{j}\right) \right| \geq c.
\end{equation*}
Then there is a subsequence $\phi _{m_{j}}$ and a sequence of local
automorphisms $\sigma _{j}\in Aut_{p_{m_{j}}}\left( B^{n+1}\right) $ such
that 
\begin{equation*}
U_{\sigma _{j}}=id_{n\times n},\quad \rho _{\sigma _{j}}=\varepsilon
_{m_{j}},\quad a_{\sigma _{j}}=0,\quad r_{\sigma _{j}}=0
\end{equation*}
and the composition 
\begin{equation*}
\sigma _{j}^{-1}\circ \phi _{m_{j}}:B^{n+1}\rightarrow B^{n+1}
\end{equation*}
uniformly converges to an automorphism of the unit ball $B^{n+1}.$
\end{lemma}

\proof%
Note that there is a subsequence $\phi _{m_{j}}$ which converges to the
point $p\in bB^{n+1}$ uniformly on every compact subset of the unit ball $%
B^{n+1}.$ We take a point $p^{\prime }\in bB^{n+1}$ such that $p^{\prime
}\neq p$ and 
\begin{equation*}
p^{\prime }\notin \overline{\left\{ p_{m_{j}}:j\in \Bbb{N}^{+}\right\} },
\end{equation*}
if necessary, passing to a subsequence.. Then, by Lemma \ref{decomposition},
there is a unique decomposition of the automorphism $\phi _{m_{j}}$ such
that 
\begin{equation*}
\phi _{m_{j}}=\varphi _{j}\circ \psi _{j}
\end{equation*}
where 
\begin{equation*}
\varphi _{j}\in Aut_{p_{m_{j}}}\left( bB^{n+1}\right) ,\quad \psi _{j}\in
Aut_{p^{\prime }}\left( bB^{n+1}\right)
\end{equation*}
and the local automorphism $\psi _{j}$ acts trivially on the complex tangent
hyperplane of $bB^{n+1}$ at the fixed point $p^{\prime }.$

Let $U_{\psi _{j}},\rho _{\psi _{j}},a_{\psi _{j}},r_{\psi _{j}}$ be the
normalizing parameters of the local automorphism $\psi _{j}.$ Since the
local automorphism $\psi _{j}$ acts trivially on the complex tangent
hyperplane of $bB^{n+1}$ at the fixed point $p^{\prime },$ we obtain 
\begin{equation*}
U_{\psi _{j}}=id_{n\times n},\quad \rho _{\psi _{j}}=1.
\end{equation*}
Since the sequence $\phi _{m_{j}}$ uniformly converges to the point $p,$
there is a real number $e>0$ such that 
\begin{equation}
\left| a_{\psi _{j}}\right| \leq e,\quad \left| r_{\psi _{j}}\right| \leq e.
\label{estimate1}
\end{equation}

Let $\pi _{j}$ be the complex line passing through the two points $p$ and $%
p_{m_{j}}.$ Since the points $p_{j}$ converges to the point $p$ to a
direction transversal to the complex tangent hyperplane at the point $p\in
bB^{n+1},$ the analytic disks 
\begin{equation*}
\delta _{j}=\pi _{j}\cap B^{n+1}
\end{equation*}
is uniformly bounded in their area from the below by a positive number.

Let $U_{\varphi _{j}},\rho _{\varphi _{j}},a_{\varphi _{j}},r_{\varphi _{j}}$
be the normalizing parameters of the local automorphism $\varphi _{j}.$ Note
that the analytic disk 
\begin{equation*}
\phi _{m_{j}}^{-1}\left( \delta _{j}\right)
\end{equation*}
is mapped by $\phi _{m_{j}}$ onto the analytic disk $\delta _{j},$ where the
areas of the analytic disks in both classes 
\begin{equation*}
\left| \phi _{m_{j}}^{-1}\left( \delta _{j}\right) \right| \quad \text{and}%
\quad \left| \delta _{j}\right|
\end{equation*}
are bounded from the below. Since the chain $\phi _{m_{j}}^{-1}\left(
b\delta _{j}\right) $ is mapped by $\phi _{m_{j}}$ to the chain $b\delta
_{j},$ there is a real number $e>0,$ if necessary, increasing $e,$ such that 
\begin{equation}
\left| a_{\varphi _{j}}\right| \leq e.  \label{estimate2}
\end{equation}
Since the area of the analytic disks $\delta _{j}$ is bounded by a positive
number, the point $p$ is attracted to the center of the local automorphism $%
Aut_{p_{m_{j}}}\left( bB^{n+1}\right) .$ Then, by passing to a subsequence,
if necessary and increasing $e$, there is a real number $e>0$ such that 
\begin{equation}
\left| r_{\varphi _{j}}\right| \leq e  \label{estimate3}
\end{equation}
and 
\begin{equation*}
\rho _{\varphi _{j}}\rightarrow 0\quad \text{as }j\rightarrow \infty .
\end{equation*}
Since $bB^{n+1}$ is strongly pseudoconvex, we obtain 
\begin{equation*}
\left| U_{\varphi _{j}}\right| =1.
\end{equation*}

Let $\eta _{j}$ be a local automorphism in $Aut_{p_{m_{j}}}\left(
bB^{n+1}\right) $ defined by the normalizing parameters 
\begin{equation*}
U_{\eta _{j}}=id_{n\times n},\quad \rho _{\eta _{j}}=\rho _{\varphi
_{j}},\quad a_{\eta _{j}}=0,\quad r_{\eta _{j}}=0.
\end{equation*}
Then the composition 
\begin{equation*}
\varphi _{j}^{\prime }\equiv \eta _{j}^{-1}\circ \varphi _{j}\in
Aut_{p_{m_{j}}}\left( bB^{n+1}\right)
\end{equation*}
has the same normalizing parameters of the mapping $\varphi _{j}$ except for 
$\rho _{\varphi _{j}^{\prime }}=1,$ i.e., 
\begin{equation*}
U_{\varphi _{j}^{\prime }}=U_{\varphi _{j}},\quad \rho _{\varphi
_{j}^{\prime }}=1,\quad a_{\varphi _{j}^{\prime }}=a_{\varphi _{j}},\quad
r_{\varphi _{j}^{\prime }}=r_{\varphi _{j}}.
\end{equation*}
Therefore, by passing to a subsequence, if necessary, the sequence 
\begin{equation*}
\tau _{j}\equiv \eta _{j}^{-1}\circ \phi _{m_{j}}:B^{n+1}\rightarrow B^{n+1}
\end{equation*}
converges by Hurwitz theorem to a locally biholomorphic mapping. Since $\tau
_{j}$ are automorphism of $B^{n+1},$ the sequence $\tau _{j}$ converges to
an automorphism of $B^{n+1}.$

Let $q\in \overline{B^{n+1}}$ be the limit point of the sequence 
\begin{equation*}
\tau _{j}\left( p_{m_{j}}^{\prime }\right) \rightarrow q\in \overline{B^{n+1}%
}.
\end{equation*}
We set 
\begin{equation*}
\lambda _{j}=\frac{\varepsilon _{m_{j}}}{\rho _{\varphi _{j}}}.
\end{equation*}
The sequence $\lambda _{j}$ converges to the distance between the two
distinct points $q$ and $p\in bB^{n+1}$ so that there is a real number $e>0$
such that 
\begin{equation}
e^{-1}\leq \left| \lambda _{j}\right| \leq e,  \label{estimate4}
\end{equation}
if necessary, increasing $e.$

Let $\sigma _{j}$ be a local automorphism in $Aut_{p_{m_{j}}}\left(
bB^{n+1}\right) $ defined by the normalizing parameters 
\begin{equation*}
U_{\sigma _{j}}=id_{n\times n},\quad \rho _{\sigma _{j}}=\varepsilon
_{m_{j}},\quad a_{\sigma _{j}}=0,\quad r_{\sigma _{j}}=0.
\end{equation*}
Then the composition 
\begin{equation*}
\varphi _{j}^{\prime \prime }\equiv \sigma _{j}^{-1}\circ \varphi _{j}\in
Aut_{p_{m_{j}}}\left( bB^{n+1}\right)
\end{equation*}
has the same normalizing parameters of the mapping $\varphi _{j}$ except for 
$\rho _{\varphi _{j}^{\prime \prime }}=\lambda _{j}^{-1},$ i.e., 
\begin{equation*}
U_{\varphi _{j}^{\prime \prime }}=U_{\varphi _{j}},\quad \rho _{\varphi
_{j}^{\prime \prime }}=\lambda _{j}^{-1},\quad a_{\varphi _{j}^{\prime
\prime }}=a_{\varphi _{j}},\quad r_{\varphi _{j}^{\prime \prime
}}=r_{\varphi _{j}}.
\end{equation*}
Therefore, by the estimates \ref{estimate1}, \ref{estimate2}, \ref{estimate3}%
, \ref{estimate4}, the sequence 
\begin{equation*}
\tau _{j}^{\prime }\equiv \sigma _{j}^{-1}\circ \phi
_{m_{j}}:B^{n+1}\rightarrow B^{n+1}
\end{equation*}
uniformly converges to an automorphism of $B^{n+1}.$ This completes the
proof.%
\endproof%

\begin{lemma}
\label{alongchain}Let $D$ be a bounded domain in $\Bbb{C}^{n+1}$ with
spherical real analytic boundary $bD.$ Suppose that there is a connected
open neighborhood $U$ of a point $p\in bD$ and a biholomorphic mapping $\phi 
$ on $U$ such that $\phi \left( U\cap bD\right) \subset bB^{n+1}$ and the
inverse mapping $\phi ^{-1}$ on $bB^{n+1}$ is analytically continued along
every chain of $bB^{n+1}.$ Then the mapping $\phi $ is analytically
continued to a biholomorphic mapping from $D$ onto $B^{n+1}.$
\end{lemma}

\proof%
The chain on $bB^{n+1}$ is characterized to be the intersection of a complex
line on $bB^{n+1}$. Thus the chains on $bB^{n+1}$ form a continuous family
so that the analytic continuity of the inverse mapping $\phi ^{-1}$ along
every chain on $bB^{n+1}$ is equivalent to the analytic continuity along any
path on $bB^{n+1}$.

Note that the inverse mapping $\phi ^{-1}$ is analytically continued along
any path on $bB^{n+1}$ and, by Hartogs extension theorem, the branching
locus of a proper mapping cannot be bounded by $bD.$ Since $bB^{n+1}$ is
simply connected, by the monodromy theorem, the mapping $\phi ^{-1}$ is
analytically continued to, by abuse of notation, a locally biholomorphic
proper mapping $\phi ^{-1}:B^{n+1}\rightarrow D$ such that $\phi ^{-1}\left(
bB^{n+1}\right) =bD.$ By the fixed point property of the closed ball $%
\overline{B^{n+1}},$ the proper mapping $\phi ^{-1}:B^{n+1}\rightarrow D$ is
globally one-to-one so that the mapping $\phi ^{-1}:B^{n+1}\rightarrow D$ is
biholomorphic. This completes the proof.%
\endproof%

\begin{lemma}
\label{circle}Let $D$ be a domain in $\Bbb{C}^{n+1}$ with spherical real
analytic boundary $bD$ and $U$ be a connected open neighborhood of a point $%
p\in bD.$ For a chain $\gamma $ on $bD$ passing through the point $p$ and
tangential to the direction with an angle $\theta _{\gamma }$ with respect
to a unit vector $v_{p}$ at the point $p$ perpendicular to the complex
tangent hyperplane, we denote $\left| a_{\gamma }\right| \equiv \left| \tan
\theta _{\gamma }\right| $. Then there is a real number $e>0$ such that
every chain $\gamma $ passing through the point $p$ with $\left| a_{\gamma
}\right| \geq e$ is the boundary of a nonsingular analytic disk $\delta
_{\gamma }\subset U\cap D$ such that 
\begin{equation*}
\gamma =\overline{\delta _{\gamma }}\cap bD.
\end{equation*}
\end{lemma}

\proof%
We take a biholomorphic mapping $\phi $ on the open set $U,$ if necessary,
shrinking $U$ such that 
\begin{equation*}
\phi \left( U\cap bD\right) \subset bB^{n+1}.
\end{equation*}
Every chain $\lambda $ on $bB^{n+1}$ is an intersection with a complex line $%
\pi _{\lambda }$ such that 
\begin{equation*}
\lambda =\pi _{\lambda }\cap bB^{n+1}.
\end{equation*}
For a sufficiently large real number $e>0,$ each chain $\gamma $ with $%
\left| a_{\gamma }\right| \geq e$ is obtain by the relation 
\begin{equation*}
\gamma =\phi ^{-1}\left( \lambda \right)
\end{equation*}
where the chain $\lambda $ on $bB^{n+1}$ satisfies the condition 
\begin{equation*}
\pi _{\lambda }\cap B^{n+1}\subset \phi \left( U\cap D\right) .
\end{equation*}
Then we take 
\begin{equation*}
\delta _{\gamma }=\phi ^{-1}\left( \pi _{\lambda }\cap B^{n+1}\right) .
\end{equation*}
This completes the proof.%
\endproof%

\begin{lemma}[First Dogginal Lemma]
\label{dogginal}Let $D$ be a bounded domain in $\Bbb{C}^{n+1}$ with
spherical real analytic boundary $bD$ and $\phi $ be a biholomorphic mapping
on a connected open neighborhood $U$ of a point $p\in bB^{n+1}$ satisfying 
\begin{equation*}
\phi \left( U\cap bB^{n+1}\right) \subset bD.
\end{equation*}
Suppose that there is a chain-segment $\gamma :[0,1]\rightarrow bB^{n+1}$
such that $\gamma \left( 0\right) \in U\cap bB^{n+1}$ and the mapping $\phi $
is analytically continued along the subpath $\gamma [0,\tau ]$ for all $\tau
<1,$ but not the whole path $\gamma [0,1]$ as a local biholomorphic mapping.
Let $\pi $ be the complex line containing the chain-segment $\gamma [0,1].$
Then there is an open neighborhood $V$ along the path $\gamma [0,1]$ such
that

\begin{enumerate}
\item  $\gamma [0,\tau ]\subset V\quad $for all $\tau <1,$

\item  $bV\cap \pi \cap B\left( \gamma \left( 1\right) ;\delta \right) $ is
an angle for a sufficient small $\delta >0,$ which contains the
chain-segment $\gamma [0,1],$

\item  $bV\cap bB^{n+1}\cap B\left( \gamma \left( 1\right) ;\delta \right) $
is slanted paraboloid for a sufficiently small $\delta >0,$ which smoothly
touches the complex tangent hyperplane at the point $\gamma \left( 1\right)
, $

\item  the mapping $\phi $ is analytically continued on $V$ as a local
biholomorphic mapping.
\end{enumerate}
\end{lemma}

\proof%
By the analytic continuation of the mapping $\phi $ along the subpath $%
\gamma [0,\tau ]$ for all $\tau <1$, there is a path $\phi \circ \gamma
:[0,1)\rightarrow bD.$ Then we consider the following sequences 
\begin{eqnarray*}
p_{j} &=&\gamma \left( 1-\frac{1}{j}\right) ,\quad \text{for }j\in \Bbb{N}%
^{+}, \\
p_{j}^{\prime } &=&\phi \circ \gamma \left( 1-\frac{1}{j}\right) ,\quad 
\text{for }j\in \Bbb{N}^{+}.
\end{eqnarray*}
Since $bD$ is compact, there is a subsequence $p_{m_{j}}^{\prime }$ and a
point $p^{\prime }\in bD$ such that 
\begin{equation*}
p_{m_{j}}^{\prime }\rightarrow p^{\prime }.
\end{equation*}
By Theorem \ref{anypath}, the mapping $\phi ^{-1}$ is analytically continued
along the path $\phi \circ \gamma [0,1)\subset bD.$

Let $\varphi _{j}$ be the analytic continuation of the mapping $\phi ^{-1}$
at the point $p_{m_{j}}^{\prime }$ along the path $\phi \circ \gamma
[0,1)\subset bD.$ By Theorem \ref{anypath}, there is an open neighborhood $W$
of the point $p^{\prime }$ such that $\varphi _{j}$ is locally biholomorphic
on an open neighborhood of $W\cap bD.$ By Lemma \ref{localhull}, we may
assume that $\varphi _{j}$ is holomorphic on $W\cap D,$ if necessary,
shrinking $W.$

Since $bD$ is spherical, there is an open neighborhood $W$ of the point $%
p^{\prime },$ if necessary, shrinking $W,$ and a biholomorphic mapping $%
\varphi $ on $W$ such that 
\begin{equation*}
\varphi \left( W\cap bD\right) \subset bB^{n+1}.
\end{equation*}
Then, by Lemma \ref{sphere}, the compositions 
\begin{equation*}
\phi _{j}\equiv \varphi _{j}\circ \varphi ^{-1}:\varphi \left( W\right) \cap
bB^{n+1}\subset bB^{n+1}
\end{equation*}
are analytically continued, by abuse of notation, to automorphisms $\phi
_{j} $ of the unit ball $B^{n+1}.$ Without loss of generality, we may assume
that the sequence $\phi _{j}$ converges to a holomorphic mapping uniformly
on every compact subset of $B^{n+1}.$

We set 
\begin{equation*}
p_{m_{j}}^{\prime \prime }\equiv \varphi \left( p_{m_{j}}^{\prime }\right)
\in bB^{n+1}\quad \text{and\quad }p^{\prime \prime }\equiv \varphi \left(
p^{\prime }\right)
\end{equation*}
so that 
\begin{equation*}
p_{m_{j}}^{\prime \prime }\rightarrow p^{\prime \prime }.
\end{equation*}
Hence the relation 
\begin{equation*}
\phi _{j}\left( p_{m_{j}}^{\prime \prime }\right) =p_{m_{j}}
\end{equation*}
yields 
\begin{equation*}
\phi _{j}\left( p_{m_{j}}^{\prime \prime }\right) \rightarrow p.
\end{equation*}

By Lemma \ref{either}, the mapping $\phi _{j}$ converges to either a
constant mapping or an automorphism of $B^{n+1}.$ We claim that the sequence 
$\phi _{j}$ converges to the point $p$ uniformly on every compact subset of $%
B^{n+1}.$ Otherwise, the sequence $\phi _{j}$ converges to an automorphism
of $B^{n+1}$ so that the sequence 
\begin{equation*}
\varphi _{j}=\phi _{j}\circ \varphi
\end{equation*}
converges to a biholomorphic mapping on an open neighborhood $W$ of the
point $p^{\prime }.$ Then there is a real number $\delta $ such that the
mapping $\varphi _{j}$ and its inverse $\varphi _{j}^{-1}$ are analytically
continued respectively on 
\begin{equation*}
B\left( p_{m_{j}}^{\prime };\delta \right) \quad \text{and}\quad B\left(
p_{m_{j}};\delta \right)
\end{equation*}
as a locally biholomorphic mapping. Thus the mapping $\phi =\varphi
_{j}^{-1} $ is biholomorphic on $B\left( p_{m_{j}};\delta \right) $ for
every point $p_{m_{j}}\rightarrow p.$ This is impossible by the hypothesis
on the mapping $\phi .$ Thus, by Lemma \ref{either}, the sequence $\phi _{j}$
converges to the point $p$ on every compact subset of $B^{n+1}.$

Let $\pi _{j}$ be the complex line passing through the two points $p$ and $%
p_{m_{j}},$ and $\delta _{j}$ be the analytic disk 
\begin{equation*}
\delta _{j}=\pi _{j}\cap B^{n+1}.
\end{equation*}
Let $\varepsilon _{j}$ be the euclidean length between the two points $p$
and $p_{m_{j}}.$ Note that the area $\left| \phi _{j}^{-1}\left( \delta
_{j}\right) \right| $ of the analytic disk 
\begin{equation*}
\phi _{j}^{-1}\left( \delta _{j}\right)
\end{equation*}
is bounded from the below. Otherwise, by Lemma \ref{circle}, the mapping $%
\phi =\varphi _{j}^{-1}$ is analytically continued over the point $p$ as a
locally biholomorphic mapping. Therefore, by First Scaling Lemma, there is a
sequence $\sigma _{j}\in Aut_{p_{m_{j}}}\left( bB^{n+1}\right) $ such that 
\begin{equation*}
U_{\sigma _{j}}=id_{n\times n},\quad \rho _{\sigma _{j}}=\varepsilon
_{j},\quad a_{\sigma _{j}}=0,\quad r_{\sigma _{j}}=0
\end{equation*}
and the sequence 
\begin{equation*}
\psi _{j}\equiv \sigma _{j}^{-1}\circ \phi _{j}
\end{equation*}
converges to an automorphism of $B^{n+1}.$

Note that there is a real number $\delta >0$ such that the mapping $\varphi $
and its inverse $\varphi ^{-1}$ are biholomorphically continued respectively
on 
\begin{equation*}
B\left( p_{m_{j}}^{\prime };\delta \right) \quad \text{and}\quad B\left(
p_{m_{j}}^{\prime \prime };\delta \right) ,
\end{equation*}
if necessary, passing to a subsequence. Further, there is a real number $%
\delta >0$ such that the mapping $\psi _{j}$ and its inverse $\psi _{j}^{-1}$
are biholomorphically continued respectively to 
\begin{equation*}
B\left( p_{m_{j}}^{\prime };\delta \right) \quad \text{and}\quad B\left(
p_{m_{j}};\delta \right) ,
\end{equation*}
if necessary, shrinking $\delta .$ Then the mapping 
\begin{eqnarray*}
\phi &=&\varphi _{j}^{-1} \\
&=&\varphi ^{-1}\circ \psi _{j}^{-1}\circ \sigma _{j}^{-1}
\end{eqnarray*}
is biholomorphically continued on the open neighborhood 
\begin{equation*}
\sigma _{j}\left( B\left( p_{m_{j}};\delta \right) \right) .
\end{equation*}

For a canonical normalizing mapping $\mu _{p_{m_{j}}},$ we obtain 
\begin{equation*}
\sigma _{j}^{\prime }\equiv \mu _{p_{m_{j}}}\circ \sigma _{j}\circ \mu
_{p_{m_{j}}}^{-1}:\left\{ 
\begin{array}{l}
z^{*}=\sqrt{\varepsilon _{j}}z \\ 
w^{*}=\varepsilon _{j}w
\end{array}
\right. .
\end{equation*}
Since $p_{m_{j}}\rightarrow p,$ by Lemma \ref{canonicalmapping}, there is a
real number $\delta >0$ such that the mapping $\mu _{p_{m_{j}}}$ and its
inverse $\mu _{p_{m_{j}}}^{-1}$ are biholomorphically continued respectively
to 
\begin{equation*}
B\left( p_{m_{j}};\delta \right) \quad \text{and}\quad B\left( 0;\delta
\right) .
\end{equation*}
Hence the mapping 
\begin{equation*}
\phi \circ \mu _{p_{m_{j}}}^{-1}=\varphi ^{-1}\circ \psi _{j}^{-1}\circ \mu
_{p_{m_{j}}}^{-1}\circ \sigma _{j}^{\prime -1}
\end{equation*}
is biholomorphically continued on 
\begin{equation*}
\sigma _{j}^{\prime }\left( B\left( 0;\delta \right) \right) .
\end{equation*}
Since $p_{m_{j}}\rightarrow p,$ by Lemma \ref{canonicalmapping}, the
canonical normalizing mapping $\mu _{p_{m_{j}}}$ uniformly converges to the
mapping $\mu _{p}$ so that the mapping $\phi $ is biholomorphically
continued near the point $p_{m_{j}}$ on 
\begin{equation*}
\mu _{p_{m_{j}}}^{-1}\circ \sigma _{j}^{\prime }\left( B\left( 0;\delta
\right) \right) .
\end{equation*}
Therefore, the analytically continued region of the mapping $\phi $ along
the chain $\gamma [0,1)\subset bB^{n+1}$ contains an open set along the
chain $\gamma [0,1]$ which touches to the point $\gamma \left( 1\right) $ by
an edge shape transversal to $bB^{n+1}$ and by a slanted paraboloid shape on 
$bB^{n+1}.$ This completes the proof.%
\endproof%

\subsection{Doggaebi variety on a sphere}

Let $\phi $ be a biholomorphic mapping on an open neighborhood $U$ of a
point $p\in bD$ satisfying $\phi \left( U\cap bD\right) \subset bB^{n+1}.$
Let $L\subset bB^{n+1}$ be the singular locus of the analytic continuation
of the mapping $\phi ^{-1}$, which shall be called the Doggaebi variety
associated to the mapping $\phi .$

\begin{lemma}[Second Scaling Lemma]
Let $p$ be a point of the boundary $bB^{n+1}$ and $p_{j},$ $j\in \Bbb{N}%
^{+}, $ be a sequence of points of $bB^{n+1}$ such that $p_{j}\neq p$ for
all $j$ and $p_{j}\rightarrow p$ as $j\rightarrow \infty $ to a direction
tangential to the complex tangent hyperplane at the point $p\in bB^{n+1}.$
Let $\varepsilon _{j}$ be the euclidean distance between the two points $%
p_{j}$ and $p,$ and $\delta _{j}$ be the analytic disk 
\begin{equation*}
\delta _{j}=\pi _{j}\cap B^{n+1}
\end{equation*}
where $\pi _{j}$ is the complex line passing through the two points $p_{j}$
and $p.$ Suppose that there is a sequence $p_{j}^{\prime }$ of points of $%
bB^{n+1}$ satisfying 
\begin{equation*}
p_{j}^{\prime }\rightarrow p^{\prime }\in bB^{n+1},
\end{equation*}
and a sequence of biholomorphic automorphisms $\phi _{j}\in Aut\left(
B^{n+1}\right) $ satisfying 
\begin{equation*}
\phi _{j}\left( p_{j}^{\prime }\right) =p_{j}
\end{equation*}
such that the area $\left| \phi _{j}^{-1}\left( \delta _{j}\right) \right| $
of the analytic disks 
\begin{equation*}
\phi _{j}^{-1}\left( \delta _{j}\right)
\end{equation*}
is bounded from the below, i.e., there is a real number $c>0$ satisfying 
\begin{equation*}
\left| \phi _{j}^{-1}\left( \delta _{j}\right) \right| \geq c.
\end{equation*}
Then there is a subsequence $\phi _{m_{j}}$ and a sequence of local
automorphisms $\sigma _{j}\in Aut_{p_{m_{j}}}\left( bB^{n+1}\right) $ such
that 
\begin{equation*}
U_{\sigma _{j}}=id_{n\times n},\quad \rho _{\sigma _{j}}=\varepsilon
_{m_{j}}^{2},\quad a_{\sigma _{j}}=0,\quad r_{\sigma _{j}}=0
\end{equation*}
and the composition 
\begin{equation*}
\sigma _{j}^{-1}\circ \phi _{m_{j}}:B^{n+1}\rightarrow B^{n+1}
\end{equation*}
uniformly converges to an automorphism of the unit ball $B^{n+1}.$
\end{lemma}

\proof%
Note that there is a subsequence $\phi _{m_{j}}$ which converges to the
point $p\in bB^{n+1}$ uniformly on an open neighborhood of the closed ball $%
\overline{B^{n+1}}.$ We take a point $p^{\prime \prime }\in bB^{n+1}$ such
that $p^{\prime \prime }\neq p$ and 
\begin{equation*}
p^{\prime \prime }\notin \overline{\left\{ p_{m_{j}}:j\in \Bbb{N}%
^{+}\right\} },
\end{equation*}
if necessary, passing to a subsequence.. Then, by Corollary \ref
{decomposition}, there is a unique decomposition of the automorphism $\phi
_{m_{j}}$ such that 
\begin{equation*}
\phi _{m_{j}}=\varphi _{j}\circ \psi _{j}
\end{equation*}
where 
\begin{equation*}
\varphi _{j}\in Aut_{p_{m_{j}}}\left( bB^{n+1}\right) ,\quad \psi _{j}\in
Aut_{p^{\prime \prime }}\left( bB^{n+1}\right)
\end{equation*}
and the local automorphism $\psi _{j}$ acts trivially on the tangent
hyperplane of $bB^{n+1}$ at the fixed point $p^{\prime }.$

Let $U_{\psi _{j}},\rho _{\psi _{j}},a_{\psi _{j}},r_{\psi _{j}}$ be the
normalizing parameters of the local automorphism $\psi _{j}.$ Since the
local automorphism $\psi _{j}$ acts trivially on the tangent hyperplane of $%
bB^{n+1}$ at the fixed point $p^{\prime \prime },$ we obtain 
\begin{equation*}
U_{\psi _{j}}=id_{n\times n},\quad \rho _{\psi _{j}}=1.
\end{equation*}
Since the sequence $\phi _{m_{j}}$ uniformly converges to the point $p,$
there is a real number $e>0$ such that 
\begin{equation}
\left| a_{\psi _{j}}\right| \leq e,\quad \left| r_{\psi _{j}}\right| \leq e.
\label{estimate5}
\end{equation}

Let $\pi _{j}$ be the complex line passing through the two points $p$ and $%
p_{m_{j}},$ and $\delta _{j}$ be the analytic disks 
\begin{equation*}
\delta _{j}=\pi _{j}\cap B^{n+1}.
\end{equation*}
Let $\sigma _{j}$ be a local automorphism in $Aut_{p_{m_{j}}}\left(
bB^{n+1}\right) $ defined by the normalizing parameters 
\begin{equation*}
U_{\sigma _{j}}=id_{n\times n},\quad \rho _{\sigma _{j}}=\varepsilon
_{m_{j}}^{2},\quad a_{\sigma _{j}}=0,\quad r_{\sigma _{j}}=0.
\end{equation*}
Since the points $p_{j}$ converges to the point $p$ to a direction
tangential to the complex tangent hyperplane at the point $p\in bB^{n+1},$
the area $\left| \sigma _{j}^{-1}\left( \delta _{j}\right) \right| $ of the
analytic disk 
\begin{equation*}
\sigma _{j}^{-1}\left( \delta _{j}\right)
\end{equation*}
is bounded from below.

Let $U_{\varphi _{j}},\rho _{\varphi _{j}},a_{\varphi _{j}},r_{\varphi _{j}}$
be the normalizing parameters of the local automorphism $\varphi _{j}.$ Note
that the analytic disk 
\begin{equation*}
\phi _{m_{j}}^{-1}\left( \delta _{j}\right)
\end{equation*}
is mapped by $\sigma _{j}^{-1}\circ \phi _{m_{j}}$ onto the analytic disk 
\begin{equation*}
\sigma _{j}^{-1}\left( \delta _{j}\right)
\end{equation*}
where the areas of the analytic disks in both classes are bounded from
below. Further, the normalizing parameters $a_{\varphi _{j}^{\prime
}},r_{\varphi _{j}^{\prime }}$ of the composition $\varphi _{j}^{\prime
}=\sigma _{j}^{-1}\circ \varphi _{_{j}}$ is the same value of $a_{\varphi
_{j}},r_{\varphi _{j}},$ i.e., 
\begin{equation*}
a_{\varphi _{j}^{\prime }}=a_{\varphi _{j}},\quad r_{\varphi _{j}^{\prime
}}=r_{\varphi _{j}}
\end{equation*}
so that there is a real number $e>0,$ if necessary, increasing $e,$ such
that 
\begin{equation}
\left| a_{\varphi _{j}}\right| \leq e.  \label{estimate6}
\end{equation}
Since the area of the analytic disks $\sigma _{j}^{-1}\left( \delta
_{j}\right) $ is bounded by a positive number, the point $p$ should be
attracted to the center of the local automorphism $Aut_{p_{m_{j}}}\left(
bB^{n+1}\right) .$ Hence, by passing to a subsequence, if necessary and
increasing $e$, there is a real number $e>0$ such that 
\begin{equation}
\left| r_{\varphi _{j}}\right| \leq e  \label{estimate7}
\end{equation}
and 
\begin{equation*}
\rho _{\varphi _{j}}\rightarrow 0\quad \text{as }j\rightarrow \infty .
\end{equation*}
Since $bB^{n+1}$ is strongly pseudoconvex, we obtain 
\begin{equation*}
\left| U_{\varphi _{j}}\right| =1.
\end{equation*}

We set 
\begin{equation*}
\lambda _{j}=\frac{\varepsilon _{m_{j}}^{2}}{\rho _{\varphi _{j}}}
\end{equation*}
so that there is a real number $e>0$ such that 
\begin{equation}
e^{-1}\leq \left| \lambda _{j}\right| \leq e,  \label{estimate8}
\end{equation}
if necessary, increasing $e.$ Then the composition 
\begin{equation*}
\varphi _{j}^{\prime }\equiv \sigma _{j}^{-1}\circ \varphi _{j}\in
Aut_{p_{m_{j}}}\left( bB^{n+1}\right)
\end{equation*}
has the same normalizing parameters of the mapping $\varphi _{j}$ except for 
$\rho _{\varphi _{j}^{\prime }}=\lambda _{j}^{-1},$ i.e., 
\begin{equation*}
U_{\varphi _{j}^{\prime }}=U_{\varphi _{j}},\quad \rho _{\varphi
_{j}^{\prime }}=\lambda _{j}^{-1},\quad a_{\varphi _{j}^{\prime
}}=a_{\varphi _{j}},\quad r_{\varphi _{j}^{\prime }}=r_{\varphi _{j}}.
\end{equation*}
Therefore, by the estimates \ref{estimate5}, \ref{estimate6}, \ref{estimate7}%
, \ref{estimate8}, the sequence 
\begin{equation*}
\tau _{j}\equiv \sigma _{j}^{-1}\circ \phi _{m_{j}}:B^{n+1}\rightarrow
B^{n+1}
\end{equation*}
uniformly converges to an automorphism of $B^{n+1}.$ This completes the
proof.%
\endproof%

\begin{theorem}
\label{finite}Let $D$ be a bounded domain in $\Bbb{C}^{n+1}$ with spherical
real analytic boundary $bD.$ Suppose that there is a biholomorphic mapping $%
\phi $ on a connected open neighborhood $U$ of a point $p\in bD$ satisfying 
\begin{equation*}
\phi \left( U\cap bD\right) \subset bB^{n+1}.
\end{equation*}
Then the Doggaebi variety $L$ associated to the mapping $\phi $ is a finite
subset of $bB^{n+1}$ such that the inverse mapping $\phi ^{-1}$ is
analytically continued along any piecewise chain curve on $%
bB^{n+1}\backslash L$ as a locally biholomorphic mapping.
\end{theorem}

\proof%
By First Dogginal Lemma, the singular locus of the analytic continuation of
the inverse mapping $\phi ^{-1}$ on $bB^{n+1}$ is an integral manifold of
the subdistribution of the complex tangent hyperplanes on $bB^{n+1}$ in its
smooth part. Thus the mapping $\phi ^{-1}$ cannot be branched on $bB^{n+1}$
by a branching locus passing through the boundary $bB^{n+1}.$ Otherwise, the
intersection of the nontrivial branch of the mapping $\phi ^{-1}$ to the
boundary $bB^{n+1}$ would be a nontrivial complex submanifold on $bB^{n+1}.$
Hence the singular locus of the analytic continuation of the mapping $\phi
^{-1}$ is well defined so that the inverse mapping $\phi ^{-1}$ is
analytically continued on $bB^{n+1}\backslash L$ as a locally biholomorphic
mapping.

We take a chain-segment $\gamma :[0,1]\rightarrow bB^{n+1}$ with $\gamma
\left( 1\right) =p\in L$ such that a germ of the mapping $\phi ^{-1}$ at the
point $\gamma \left( 0\right) $ is analytically continued along all subpath $%
\gamma [0,\tau ]$ with $\tau <1,$ but not the whole path $\gamma [0,1].$
Note that the distribution of the complex tangent hyperplanes on $bB^{n+1}$
is maximally nonintegrable. Thus, by First Dogginal Lemma, the singular
locus $L$ of the analytic continuation of the mapping $\phi ^{-1}$ cannot
bound the open region of the analytic continuation of the mapping $\phi
^{-1}.$ Then the mapping $\phi ^{-1}$ is analytically continued on the
opposite side of the complex tangent hyperplane $H_{p}$ at the point $p\in
bB^{n+1}.$ Thus there is a complex line $\pi $ passing through the point $p$
and an open neighborhood $U$ of the point $p$ such that 
\begin{equation*}
\pi \cap bB^{n+1}
\end{equation*}
is a chain on $bB^{n+1}$ satisfying 
\begin{equation*}
L\cap U\cap \pi \cap bB^{n+1}=\left\{ p\right\} .
\end{equation*}

We claim that, if necessary, shrinking $U,$%
\begin{equation*}
L\cap U\cap bB^{n+1}=\left\{ p\right\}
\end{equation*}
so that the singular locus $L$ is a finite set on $bB^{n+1}.$ By First
Dogginal Lemma, there is a sequence $p_{j}\in bB^{n+1}\backslash L$ such
that $p_{j}\rightarrow p$ and the sequence $p_{j}$ converges to the point $p$
to a direction tangential to the complex tangent hyperplane at the point $%
p\in L\subset bB^{n+1}.$ Since $bD$ is compact, there is a subsequence $%
p_{m_{j}}$ and a point $p^{\prime }\in bD$ such that 
\begin{equation*}
p_{m_{j}}^{\prime }\equiv \phi ^{-1}\left( p_{m_{j}}\right) \rightarrow
p^{\prime }.
\end{equation*}
Let $\phi _{j}$ be the germ of the mapping $\phi $ at the point $%
p_{m_{j}}^{\prime }\in bD$ such that 
\begin{equation*}
\phi _{j}\left( p_{m_{j}}^{\prime }\right) =p_{m_{j}}.
\end{equation*}
Since $bD$ is spherical, there is an open neighborhood $W$ of the point $%
p^{\prime }$ and a biholomorphic mapping $\psi $ on $W$ such that 
\begin{equation*}
\psi \left( W\cap bD\right) \subset bB^{n+1}.
\end{equation*}
Then the compositions $\varphi _{j}\equiv \phi _{j}\circ \psi ^{-1}$ satisfy 
\begin{equation*}
\varphi _{j}\left( \psi \left( W\right) \cap bB^{n+1}\right) \subset
bB^{n+1},
\end{equation*}
if necessary, shrinking $W$ so that, by abuse of notation, the mapping $%
\varphi _{j}$ is an automorphism of the unit ball $B^{n+1}.$

Let $\pi _{j}$ be the complex line passing through the points $p_{m_{j}}$
and $p,$ and $\delta _{j}$ be the analytic disk 
\begin{equation*}
\delta _{j}=\pi _{j}\cap B^{n+1}.
\end{equation*}
Note that the area $\left| \varphi _{j}^{-1}\left( \delta _{j}\right)
\right| $ of the analytic disk $\varphi _{j}^{-1}\left( \delta _{j}\right) $
is bounded from the below. Otherwise, $p\in bB^{n+1}\backslash L.$ Further,
the sequence $\varphi _{j}$ converges to the point $p$ uniformly on every
compact subset of $B^{n+1},$ if necessary, passing to a subsequence.
Otherwise, the sequence $\varphi _{j}$ converges to an automorphism of the
unit ball $B^{n+1}$ so that $p\in bB^{n+1}\backslash L.$ Let $\varepsilon
_{j}$ be the euclidean length between the two points $p$ and $p_{m_{j}}.$
Then, by Second Scaling Lemma, there is a sequence of local automorphisms $%
\sigma _{j}\in Aut_{p_{m_{j}}}\left( bB^{n+1}\right) $ such that 
\begin{equation*}
U_{\sigma _{j}}=id_{n\times n},\quad \rho _{\sigma _{j}}=\varepsilon
_{j}^{2},\quad a_{\sigma _{j}}=0,\quad r_{\sigma _{j}}=0
\end{equation*}
and the composition 
\begin{equation*}
\tau _{j}\equiv \sigma _{j}^{-1}\circ \varphi _{j}:B^{n+1}\rightarrow B^{n+1}
\end{equation*}
uniformly converges to an automorphism of the unit ball $B^{n+1}.$ Thus
there is a positive real number $\delta >0$ such that the mapping $\tau _{j}$
and its inverse $\tau _{j}^{-1}$ are analytically continued respectively on 
\begin{equation*}
B\left( p_{m_{j}}^{\prime \prime };\delta \right) \quad \text{and}\quad
B\left( p_{m_{j}};\delta \right)
\end{equation*}
where 
\begin{equation*}
p_{m_{j}}^{\prime \prime }=\psi \left( p_{m_{j}}^{\prime }\right) .
\end{equation*}
Therefore the mapping $\phi ^{-1}$ at the point $p_{m_{j}}:$%
\begin{eqnarray*}
\phi ^{-1} &=&\phi _{j}^{-1} \\
&=&\psi ^{-1}\circ \tau _{j}^{-1}\circ \sigma _{j}^{-1}
\end{eqnarray*}
is analytically continued to the open neighborhood 
\begin{equation*}
\sigma _{j}\left( B\left( p_{m_{j}};\delta \right) \right) .
\end{equation*}
By the canonical normalizing mapping $\mu _{p_{m_{j}}},$ we obtain 
\begin{equation*}
\sigma _{j}^{\prime }\equiv \mu _{p_{m_{j}}}\circ \sigma _{j}\circ \mu
_{p_{m_{j}}}^{-1}:\left\{ 
\begin{array}{l}
z^{*}=\varepsilon _{j}z \\ 
w^{*}=\varepsilon _{j}^{2}w
\end{array}
.\right.
\end{equation*}
Hence the mapping $\phi ^{-1}$ at the point $p_{m_{j}}$ is analytically
continued to the open neighborhood 
\begin{equation*}
\mu _{p_{m_{j}}}^{-1}\circ \sigma _{j}^{\prime }\left( B\left( 0;\delta
\right) \right) ,
\end{equation*}
if necessary, shrinking $\delta .$ Since $\varepsilon _{j}$ is the euclidean
length between the two points $p$ and $p_{m_{j}},$ the mapping $\phi ^{-1}$
is analytically continued on an open region which touches to the point $p$
to a converging direction to the sequence $p_{j}$ by an edge shape
tangential to the complex tangent hyperplane $H_{p}$ at the point $p$ and a $%
\sqrt{\left| x\right| }$ curve shape normal to $H_{p}$ on $bB^{n+1}.$

Therefore, the singular locus $L$ of the analytic continuation of the
mapping $\phi ^{-1}$ is isolated to the direction of the complex tangent
hyperplane at each point of $L.$ Since the boundary $bB^{n+1}$ is compact,
the singular locus $L$ is a finite subset of $bB^{n+1}.$ This completes the
proof.%
\endproof%

\begin{lemma}
\label{length}Let $D$ be a bounded domain in $\Bbb{C}^{n+1}$ with spherical
real analytic boundary $bD.$ Then every chain on $bD$ is extended each
direction infinitely in its euclidean length.
\end{lemma}

\proof%
Suppose that the assertion is not true. Then there would be a path $\gamma
:[0,1]\rightarrow bD$ such that the subpath $\gamma [0,\tau ]$ is a
chain-segment for all $\tau <1$ but the whole path $\gamma [0,1]$ is not a
chain segment.

Since $bD$ is spherical, by definition, there exist a connected open
neighborhood $U$ of the point $\gamma \left( 1\right) $ and a biholomorphic
mapping $\phi $ on $U$ such that $\phi \left( U\cap bD\right) \subset
bB^{n+1}$. There is a unique closed circle $\chi $ on $bB^{n+1}$ such that $%
\chi $ is a chain on $bB^{n+1}$ and $\phi \circ \gamma [0,\tau ]\subset \chi 
$ for all $\tau <1$(cf. \cite{Pa3}).

Then the inverse image $\phi ^{-1}\left( \chi \cap \phi \left( U\right)
\right) $ under the biholomorphic mapping $\phi $ is a chain on $bB^{n+1}$%
(cf. \cite{Pa3}) such that 
\begin{equation*}
\gamma [0,1]\cap U\subset \phi ^{-1}\left( \chi \cap \phi \left( U\right)
\right) .
\end{equation*}
This is a contradiction. This completes the proof.%
\endproof%

\begin{theorem}
\label{creq}Let $D$ be a bounded domain in $\Bbb{C}^{n+1}$ with spherical
real analytic boundary $bD.$ Suppose that there is a biholomorphic mapping $%
\phi $ on a connected open neighborhood $U$ of a point $p\in bD$ satisfying 
\begin{equation*}
\phi \left( U\cap bD\right) \subset bB^{n+1}.
\end{equation*}
Let $E$ be the path space of $bD$ pointed at the point $p\in bD$ mod
homotopic relation so that $E$ be a universal covering of $bD$ with a
natural CR structure and a natural CR covering map $\psi :E\rightarrow bD$.
Then there is a unique CR equivalence $\varphi :E\rightarrow
bB^{n+1}\backslash L$ commuting the diagram 
\begin{equation*}
\begin{array}{lll}
E &  &  \\ 
\downarrow \psi & \overset{\varphi }{\searrow } &  \\ 
bD & \overset{\phi }{\rightarrow } & bB^{n+1}\backslash L
\end{array}
\end{equation*}
where $L\subset bB^{n+1}$ is the Doggaebi variety associated to the mapping $%
\phi .$
\end{theorem}

\proof%
We obtain the set $bB^{n+1}\backslash L$ when the mapping $\phi $ is
analytically continued along chains on $bD.$ Since the mapping $\varphi
:E\rightarrow bB^{n+1}$ is naturally induced by the analytic continuation of
the mapping $\phi $ on $bD$(cf. Theorem \ref{anypath}), we obtain 
\begin{equation*}
bB^{n+1}\backslash L\subset \varphi \left( E\right)
\end{equation*}

Let $\gamma :[0,1]\rightarrow bD$ be a path such that $\gamma \left(
0\right) =p.$ Then, for a sufficiently small $\varepsilon >0,$ there is an $%
\varepsilon $ open neighborhood $V$ of the path $\gamma [0,1]$ such that the
mapping $\phi $ is analytically continued on $V$ as a local biholomorphic
mapping. Then we take a piecewise chain curve $\eta :[0,1]\rightarrow bD$%
(cf. \cite{Pa3}) such that 
\begin{equation*}
\eta \left( 0\right) =\gamma \left( 0\right) ,\quad \eta \left( 1\right)
=\gamma \left( 1\right)
\end{equation*}
and 
\begin{equation*}
\eta [0,1]\subset V.
\end{equation*}
Further, we may take a continuous function $\Gamma :[0,1]\times
[0,1]\rightarrow V$ such that 
\begin{equation*}
\Gamma \left( 0,\tau \right) =\gamma \left( \tau \right) ,\quad \Gamma
\left( 1,\tau \right) =\eta \left( \tau \right) \quad \text{for }\tau \in
[0,1]
\end{equation*}
and $\Gamma \left( \cdot ,\tau \right) :[0,1]\rightarrow V$ is a
chain-segment for all $\tau \in [0,1].$ Note that every chain-segment in
this construction is finite in its euclidean length. By Lemma \ref{length},
the mapping $\phi $ is analytically continued along the whole path $\gamma
[0,1].$

Then the image of $E$ under the mapping $\varphi $ satisfies 
\begin{equation*}
\varphi \left( E\right) \subset bB^{n+1}\backslash L
\end{equation*}
which yields 
\begin{equation*}
\varphi \left( E\right) =bB^{n+1}\backslash L.
\end{equation*}
Since the mapping $\phi ^{-1}$ is analytically continued on $%
bB^{n+1}\backslash L,$ the mapping $\varphi :E\rightarrow bB^{n+1}$ is a CR
equivalence. This completes the proof.%
\endproof%

\begin{corollary}
Let $D$ be a bounded domain with spherical real analytic boundary $bD.$ Then
there is a natural embedding 
\begin{equation*}
Aut\left( D\right) \subset Aut\left( B^{n+1}\right) .
\end{equation*}
\end{corollary}

\proof%
Let $E$ be the path space of $bD$ pointed at the point $p\in bD$ mod
homotopy so that $E$ be a universal covering of $bD$ with a natural CR
structure and a natural CR covering map $\psi :E\rightarrow bD$. Then we
take a biholomorphic mapping $\phi $ on an open neighborhood $U$ of a point $%
p\in bD$ such that 
\begin{equation*}
\phi \left( U\cap bD\right) \subset bB^{n+1}.
\end{equation*}
Let $L$ be the Doggaebi variety associated to the mapping $\phi .$ Then, by
Theorem \ref{creq}, the analytic continuation of the mapping $\phi $ yields
a natural CR equivalence 
\begin{equation*}
Aut\left( E\right) \simeq Aut\left( bB^{n+1}\backslash L\right) .
\end{equation*}
Note that 
\begin{equation*}
Aut\left( D\right) =Aut(bD)\subset Aut\left( E\right) \simeq Aut\left(
bB^{n+1}\backslash L\right) \subset Aut\left( bB^{n+1}\right) =Aut\left(
B^{n+1}\right) .
\end{equation*}
This completes the proof.%
\endproof%

\begin{lemma}
\label{sph-proper}Let $D,D^{\prime }$ be bounded domains in $\Bbb{C}^{n+1}$
with spherical real analytic boundaries $bD,bD^{\prime }$, and $\phi
:D\rightarrow D^{\prime }$ be a proper holomorphic mapping. Suppose that
there is an open neighborhood $U$ of a point $p\in bD$ such that the mapping 
$\phi $ is analytically continued on $U.$ Then the mapping $\phi $ is
analytically continued along any path on $bD$ as a locally biholomorphic
mapping so that $\phi :D\rightarrow D^{\prime }$ is locally biholomorphic.
\end{lemma}

\proof%
We may obtain the result by the boundary regularity of Lemma \ref
{b-regularity} and the transformation formula for a complex Monge-Ampere
equation. Here we may give an independent proof. Note that 
\begin{equation*}
\phi \left( U\cap bD\right) \subset bD^{\prime }
\end{equation*}
so that, by shrinking $U$, if necessary, the mapping $\phi $ is
biholomorphic on $U.$ By shrinking $U$, if necessary, we take biholomorphic
mappings $\varphi ,\varphi ^{\prime }$ respectively on $U,\phi \left(
U\right) $ such that 
\begin{eqnarray*}
\varphi \left( U\cap bD\right) &\subset &bB^{n+1} \\
\varphi ^{\prime }\left( \phi \left( U\right) \cap bD^{\prime }\right)
&\subset &bB^{n+1}.
\end{eqnarray*}
Then there are Doggaebi varieties $L,L^{\prime }$ on $bB^{n+1}$ such that
the inverses $\varphi ^{-1},\varphi ^{\prime -1}$ are analytically continued
respectively on $bB^{n+1}\backslash L$ and $bB^{n+1}\backslash L^{\prime }$
as a local biholomorphic mapping. Note that the composition $\psi =\varphi
^{\prime }\circ \phi \circ \varphi ^{-1}$ satisfies the relation 
\begin{equation*}
\psi \left( \varphi \left( U\right) \cap bB^{n+1}\right) \subset bB^{n+1}.
\end{equation*}
By abuse of notation, the mapping $\psi $ is an automorphism of $bB^{n+1}$
which comments the following diagram: 
\begin{equation*}
\begin{array}{lll}
bB^{n+1}\backslash L & \overset{\psi }{\longrightarrow } & 
bB^{n+1}\backslash L^{\prime } \\ 
\downarrow \kappa &  & \downarrow \kappa ^{\prime } \\ 
bD & \overset{\phi }{\longrightarrow } & bD^{\prime }
\end{array}
.
\end{equation*}
where $\kappa :bB^{n+1}\backslash L\rightarrow bD$ is the analytic
continuation of the mapping $\varphi ^{-1}$ and $\kappa ^{\prime
}:bB^{n+1}\backslash L^{\prime }\rightarrow bD^{\prime }$ is the analytic
continuation of the mapping $\varphi ^{\prime -1}$ on the boundary $%
bB^{n+1}. $

We claim 
\begin{equation*}
L^{\prime }\subset \psi \left( L\right)
\end{equation*}
so that the mapping $\phi $ is analytically continued along any path on $bD$
as a locally biholomorphic mapping. Otherwise, there is a point $q^{\prime
}\in L^{\prime }$ such that 
\begin{equation*}
q^{\prime }\in L^{\prime }\backslash \psi \left( L\right)
\end{equation*}
From the proof of First Scaling Lemma, there is a sequence $p_{j}^{\prime
}\in bD^{\prime }$ with $p_{j}^{\prime }\rightarrow p^{\prime }\in
bD^{\prime }$ and a sequence $\varphi _{j}^{\prime }$ of the analytic
continuation of the mapping $\varphi ^{\prime }$ at the point $p_{j}^{\prime
}$ with $q_{j}^{\prime }\equiv \varphi _{j}^{\prime }\left( p_{j}^{\prime
}\right) \rightarrow q^{\prime }\in L^{\prime }\backslash \psi \left(
L\right) $ such that there is an open neighborhood $V$ of the point $%
p^{\prime }$ such that the mapping $\varphi _{j}^{\prime }$ is biholomorphic
on $V\cap D^{\prime }$ and the sequence $\varphi _{j}^{\prime }$ converges
to the point $r$ uniformly on every compact subset of $V\cap D^{\prime }.$
Since $\psi $ is an automorphism of the unit ball $B^{n+1},$ we obtain 
\begin{equation*}
\lambda _{j}\equiv \kappa \circ \psi ^{-1}\circ \varphi _{j}^{\prime }:V\cap
D^{\prime }\rightarrow D
\end{equation*}
and 
\begin{equation*}
\phi \circ \lambda _{j}=id\quad \text{on }V\cap D^{\prime }.
\end{equation*}
We set 
\begin{equation*}
p_{j}\equiv \lambda _{j}\left( p_{j}^{\prime }\right) \in bD.
\end{equation*}
Since $bD$ is compact, there is a subsequence $p_{m_{j}}\in bD$ and a point $%
p\in bD$ such that 
\begin{equation*}
p_{m_{j}}\rightarrow p.
\end{equation*}
Since $\phi $ is a proper mapping, the mapping $\phi $ is a globally finite
covering so that there is an open neighborhood $W$ of the point $p$ such
that 
\begin{equation*}
\lambda _{m_{j}}=\phi ^{-1}\quad \text{on }\phi \left( W\right) \cap D.
\end{equation*}
Hence we obtain 
\begin{equation*}
\psi =\varphi _{m_{j}}^{\prime }\circ \phi \circ \kappa
\end{equation*}

Since the mapping $\psi $ is an automorphism of the unit ball $B^{n+1},$
there is a sequence $q_{j}\in bB^{n+1}$ such that 
\begin{equation*}
\psi \left( q_{j}\right) =q_{j}^{\prime }.
\end{equation*}
Hence we obtain 
\begin{equation*}
\kappa \left( q_{m_{j}}\right) =p_{m_{j}}.
\end{equation*}
Note that there is a real number $\delta >0$ such that 
\begin{equation*}
B\left( r;2\delta \right) \cap bB^{n+1}\cap \psi \left( L\right) =\emptyset .
\end{equation*}
Passing to a subsequence, if necessary, we may assume 
\begin{equation*}
q_{m_{j}}\in \psi ^{-1}\left( B\left( r;\delta \right) \right) \cap bB^{n+1}
\end{equation*}
so that there is a point $q\in bB^{n+1}\backslash L$ satisfying 
\begin{equation*}
q_{m_{j}}\rightarrow q.
\end{equation*}
Hence we set 
\begin{equation*}
\chi _{m_{j}}\equiv \psi ^{-1}\circ \varphi _{m_{j}}^{\prime }\circ \lambda
_{m_{j}}^{-1}
\end{equation*}
and 
\begin{equation*}
\kappa \circ \chi _{m_{j}}=id.
\end{equation*}
Finally, we obtain 
\begin{equation*}
\varphi _{m_{j}}^{\prime }=\psi \circ \chi _{m_{j}}\circ \lambda _{m_{j}}.
\end{equation*}

Note that there is a real number $\delta >0$ that the mappings $\lambda
_{m_{j}}$ and the inverse mappings $\lambda _{m_{j}}^{-1}$ are analytically
continued respectively on 
\begin{equation*}
B\left( p_{m_{j}}^{\prime };\delta \right) \quad \text{and}\quad B\left(
p_{m_{j}};\delta \right) ,
\end{equation*}
and the mappings $\chi _{m_{j}}$ and the inverse mappings $\chi
_{m_{j}}^{-1} $ are analytically continued respectively on 
\begin{equation*}
B\left( p_{m_{j}};\delta \right) \quad \text{and}\quad B\left(
q_{m_{j}};\delta \right) ,
\end{equation*}
and the mapping $\psi $ and the inverse mapping $\psi ^{-1}$ are
analytically continued respectively on 
\begin{equation*}
B\left( q_{m_{j}};\delta \right) \quad \text{and}\quad B\left(
q_{m_{j}}^{\prime };\delta \right) ,
\end{equation*}
Hence the mapping $\varphi _{m_{j}}^{\prime }$ and the inverse mappings $%
\varphi _{m_{j}}^{\prime -1}=\kappa ^{\prime }$ are analytically continued
respectively on 
\begin{equation*}
B\left( p_{m_{j}}^{\prime };\delta \right) \quad \text{and}\quad B\left(
q_{m_{j}}^{\prime };\delta \right) ,
\end{equation*}
as a locally biholomorphic mapping. Since $q_{m_{j}}^{\prime }\rightarrow
q^{\prime },$ the mapping $\kappa ^{\prime }$ is analytically continued to
the point $q^{\prime }$ as a locally biholomorphic mapping so that 
\begin{equation*}
q^{\prime }\in bB^{n+1}\backslash L^{\prime }.
\end{equation*}
This is a contradiction so that 
\begin{equation*}
L^{\prime }\subset \psi \left( L\right) .
\end{equation*}
This completes the proof.%
\endproof%

\section{Bounded Domains with Spherical Boundaries}

\subsection{Differentiable spherical boundary}

Let $D$ be a domain in $\Bbb{C}^{n+1},n\geq 1,$ with a differentiable
boundary $bD.$ The boundary $bD$ shall be called spherical if, for each
point $p\in bD,$ there is a connected open neighborhood $U$ of the point $p$
and a biholomorphic mapping $\phi $ on $U\cap D$ such that 
\begin{equation}
\phi \in H\left( U\cap D\right) \cap C^{1}\left( U\cap \overline{D}\right)
,\quad \phi \left( U\cap bD\right) \subset bB^{n+1}  \label{maptosphere}
\end{equation}
and the induced mapping $\phi :U\cap bD\rightarrow bB^{n+1}$ is CR
diffeomorphic.

\begin{lemma}
\label{onsphere}Let $U$ be a connected open neighborhood of a point $p\in
bB^{n+1}$ and $\phi $ be a biholomorphic mapping on $U\cap B^{n+1}$ such
that 
\begin{equation*}
\phi \in H\left( U\cap B^{n+1}\right) \cap C^{1}\left( U\cap \overline{%
B^{n+1}}\right) ,\quad \phi \left( U\cap bB^{n+1}\right) \subset bB^{n+1}
\end{equation*}
and the induced mapping $\phi :U\cap bB^{n+1}\rightarrow bB^{n+1}$ is CR
diffeomorphic. Then the mapping $\phi $ is analytically continued to an
automorphism of the unit ball $B^{n+1}.$
\end{lemma}

\proof%
By Lemma \ref{Lewy-Pinchuk}, the mapping $\phi $ is biholomorphic on $U,$ if
necessary, shrinking $U.$ Then, by Lemma \ref{sphere}, we obtain the desired
result. This completes the proof.%
\endproof%

The chain on $bB^{n+1}$ is defined to be the intersection on $bB^{n+1}$ by a
complex line. The family of chains on $bB^{n+1}$ leaves invariant under the
action of biholomorphic automorphisms of $B^{n+1}.$ We define the chain on
the spherical differentiable boundary $bD$ of a domain $D$ to be the inverse
image of the chain on $bB^{n+1}$ under the mapping \ref{maptosphere}. By
Lemma \ref{onsphere}, the chain on the spherical differentiable boundary $bD$
is well defined so that a chain of a spherical differentiable boundary is
mapped to a chain of another spherical differentiable boundary under any
induced CR diffeomorphism.

\begin{lemma}
Let $D$ be a domain in $\Bbb{C}^{n+1}$ with spherical differentiable
boundary $bD.$ Suppose that there is a connected open neighborhood $U$ of a
point $p\in bD$ and a biholomorphic mapping $\phi $ on $U\cap D$ such that 
\begin{equation*}
\phi \in H\left( U\cap D\right) \cap C^{1}\left( U\cap \overline{D}\right)
,\quad \phi \left( U\cap bD\right) \subset bB^{n+1}.
\end{equation*}
Then the mapping $\phi $ is analytically continued along any path on $bD$ as
a local biholomorphic mapping.
\end{lemma}

\proof%
Suppose that the assertion is not true. Then there would be a path $\gamma
:[0,1]\rightarrow bD$ such that $\gamma \left( 0\right) \in U\cap bD$ and
the germ of a biholomorphic mapping $\phi $ at the point $\gamma \left(
0\right) $ is analytically continued along the subpath $\gamma [0,\tau ]$
with all $\tau <1$ as a local biholomorphic mapping, but not the whole path $%
\gamma [0,1].$

Thus we reduce the proof to a local problem near the point $\gamma \left(
1\right) \in bD.$ Then, by the definition, there is a connected open
neighborhood $V$ of the point $\gamma \left( 1\right) $ and a biholomorphic
mapping $\varphi $ on $V\cap D$ such that 
\begin{equation*}
\varphi \in H\left( V\cap D\right) \cap C^{1}\left( V\cap \overline{D}%
\right) ,\quad \varphi \left( V\cap bD\right) \subset bB^{n+1}
\end{equation*}
and the induced mapping $\varphi :U\cap bD\rightarrow bB^{n+1}$ is CR
diffeomorphic.

Then we consider the mapping 
\begin{equation*}
\phi \circ \varphi ^{-1}\in H\left( \varphi \left( V\cap D\right) \right)
\cap C^{1}\left( \varphi \left( V\cap \overline{D}\right) \right) ,\quad
\phi \circ \varphi ^{-1}\left( \varphi \left( V\cap bD\right) \right)
\subset bB^{n+1}
\end{equation*}
and the curve 
\begin{equation*}
\varphi \circ \gamma [0,1]\cap \varphi \left( V\cap bD\right) \subset
bB^{n+1}.
\end{equation*}
By Lemma \ref{Lewy-Pinchuk}, the remaining part of the proof repeats the
proof of Lemma \ref{anypath}. This completes the proof.%
\endproof%

\begin{lemma}[First Dogginal Lemma]
\label{dogginal2}Let $D$ be a bounded domain in $\Bbb{C}^{n+1}$ with
spherical differentiable boundary $bD$ and $\phi $ be a biholomorphic
mapping on $U\cap D$ for a connected open neighborhood $U$ of a point $p\in
bB^{n+1}$ satisfying 
\begin{equation*}
\phi \in H\left( U\cap B^{n+1}\right) \cap C^{1}\left( U\cap \overline{%
B^{n+1}}\right) ,\quad \phi \left( U\cap bB^{n+1}\right) \subset bD.
\end{equation*}
Suppose that there is an injective path $\gamma :[0,1]\rightarrow bB^{n+1}$
such that $\gamma [0,1]\subset bB^{n+1}$ is a chain-segment satisfying 
\begin{equation*}
\gamma \left( 0\right) \in U\cap bB^{n+1}
\end{equation*}
and the mapping $\phi $ is analytically continued along the subpath $\gamma
[0,\tau ]$ for all $\tau <1,$ but not the whole path $\gamma [0,1]$ as a
local biholomorphic mapping. Let $\pi $ be the complex line containing the
chain-segment $\gamma [0,1].$ Then there is an open neighborhood $V$ along
the path $\gamma [0,1]$ such that

\begin{enumerate}
\item  $\gamma [0,\tau ]\subset V\quad $for all $\tau <1,$

\item  $bV\cap \pi \cap B\left( \gamma \left( 1\right) ;\delta \right) $ is
an angle for a sufficient small $\delta >0,$ which contains the
chain-segment $\gamma [0,1],$

\item  $bV\cap bB^{n+1}\cap B\left( \gamma \left( 1\right) ;\delta \right) $
is paraboloid for a sufficiently small $\delta >0,$ which smoothly touches
the complex tangent hyperplane at the point $\gamma \left( 1\right) ,$

\item  the mapping $\phi $ is analytically continued on $V\cap B^{n+1}$ as a
local biholomorphic mapping 
\begin{equation*}
\phi \in H\left( V\cap B^{n+1}\right) \cap C^{1}\left( V\cap \overline{%
B^{n+1}}\right) .
\end{equation*}
\end{enumerate}
\end{lemma}

\proof%
By the analytic continuation of the mapping $\phi $ along the subpath $%
\gamma [0,\tau ]$ for all $\tau <1$, there is a path $\phi \circ \gamma
:[0,1)\rightarrow bD.$ Then we consider the following sequences 
\begin{eqnarray*}
p_{j} &=&\gamma \left( 1-\frac{1}{j}\right) ,\quad \text{for }j\in \Bbb{N}%
^{+}, \\
p_{j}^{\prime } &=&\phi \circ \gamma \left( 1-\frac{1}{j}\right) ,\quad 
\text{for }j\in \Bbb{N}^{+}.
\end{eqnarray*}
Since $bD$ is compact, there is a subsequence $p_{m_{j}}^{\prime }$ and a
point $p^{\prime }\in bD$ such that 
\begin{equation*}
p_{m_{j}}^{\prime }\rightarrow p^{\prime }.
\end{equation*}

Thus we reduce the proof to a local problem near the point $p^{\prime }\in
bD.$ Then, by the definition, there is a connected open neighborhood $W$ of
the point $p^{\prime }$ and a biholomorphic mapping $\varphi $ on $W\cap D$
such that 
\begin{equation*}
\varphi \in H\left( W\cap D\right) \cap C^{1}\left( W\cap \overline{D}%
\right) ,\quad \varphi \left( W\cap bD\right) \subset bB^{n+1}
\end{equation*}
and the induced mapping $\varphi :U\cap bD\rightarrow bB^{n+1}$ is CR
diffeomorphic.

Then we consider the mapping 
\begin{equation*}
\phi ^{-1}\circ \varphi ^{-1}\in H\left( \varphi \left( W\cap D\right)
\right) \cap C^{1}\left( \varphi \left( W\cap \overline{D}\right) \right)
,\quad \phi ^{-1}\circ \varphi ^{-1}\left( \varphi \left( W\cap bD\right)
\right) \subset bB^{n+1}.
\end{equation*}
By Lemma \ref{Lewy-Pinchuk}, the remaining part of the proof repeats the
proof of Lemma \ref{dogginal}. This completes the proof.%
\endproof%

By Lemma \ref{onsphere} and Lemma \ref{dogginal2}, we obtain the following
result by the same argument of the proof of Theorem \ref{creq}.

\begin{theorem}
\label{covering}Let $D$ be a bounded domain in $\Bbb{C}^{n+1}$ with
spherical differentiable boundary $bD.$ Suppose that there is a
biholomorphic mapping $\phi $ on $U\cap D$ for a connected open neighborhood 
$U$ of a point $p\in bD$ satisfying 
\begin{equation*}
\phi \left( U\cap bD\right) \subset bB^{n+1}.
\end{equation*}
Let $E$ be the path space of $bD$ pointed at the point $p\in bD$ mod
homotopy so that $E$ be a universal covering of $bD$ with a natural CR
structure and a natural CR covering map $\psi :E\rightarrow bD$. Then there
is a unique CR equivalence $\varphi :E\rightarrow bB^{n+1}\backslash L$
commuting the diagram 
\begin{equation*}
\begin{array}{lll}
E &  &  \\ 
\downarrow \psi & \overset{\varphi }{\searrow } &  \\ 
bD & \overset{\phi }{\rightarrow } & bB^{n+1}\backslash L
\end{array}
\end{equation*}
where $L\subset bB^{n+1}$ is the Doggaebi variety associated to the mapping $%
\phi .$
\end{theorem}

\subsection{Second Dogginal Lemma}

We have examined the analytic continuation of a biholomorphic mapping $\phi $
from the spherical differentiable boundary $bD$ of a bounded domain $D$ into
the boundary $bB^{n+1}$ of the unit ball $B^{n+1}.$ From now on, we shall
examine the analytic continuation of the mapping $\phi $ into the domain $D.$

\begin{lemma}
\label{branched}Let $D$ be a bounded domain in $\Bbb{C}^{n+1}$ with
spherical differentiable boundary $bD.$ Suppose that there is a
biholomorphic mapping $\phi $ on a connected open neighborhood $U$ of a
point $p\in bD$ satisfying 
\begin{equation*}
\phi \left( U\cap bD\right) \subset bB^{n+1}.
\end{equation*}
Then the inverse mapping $\phi ^{-1}$ is analytically continued to a locally
biholomorphic mapping from $B^{n+1}$ into $D.$
\end{lemma}

\proof%
By Theorem \ref{covering}, there is a Doggaebi variety $L$ on $bB^{n+1}$
such that the mapping $\phi ^{-1}$ is uniquely analytically continued on the
set $bB^{n+1}\backslash L.$ Then the mapping $\phi ^{-1}$ is analytically
continued to a holomorphic mapping on the unit ball $B^{n+1}.$ Hence we have
a holomorphic mapping $\varphi :B^{n+1}\rightarrow D$ which extends to the
mapping $\phi ^{-1}:bB^{n+1}\backslash L\rightarrow bD.$ Note that the zero
set of the determinant of the Jacobian matrix $\varphi ^{\prime }$ of the
mapping $\varphi $ cannot be located on the set $L$ on the boundary $%
bB^{n+1}.$ Thus the mapping $\varphi :B^{n+1}\rightarrow D$ is locally
biholomorphic. This completes the proof.%
\endproof%

\begin{lemma}
\label{transversally}Let $D$ be a bounded domain in $\Bbb{C}^{n+1}$ with
spherical differentiable boundary $bD$ and $\phi $ be a biholomorphic
mapping on an open neighborhood $U$ of a point $r\in bD$ satisfying 
\begin{equation*}
\phi \in H\left( U\cap D\right) \cap C^{1}\left( U\cap \overline{D}\right)
\quad \text{and}\quad \phi \left( U\cap bD\right) \subset bB^{n+1}.
\end{equation*}
Let $L$ be the Doggaebi variety associated to the mapping $\phi $ and $%
\varphi :B^{n+1}\rightarrow D$ be a locally biholomorphic mapping to be an
analytic continuation of the mapping $\phi ^{-1}$. Suppose that there is a
line segment $\gamma :[0,1]\rightarrow D$ with 
\begin{equation*}
p=\gamma \left( 1\right) \quad \text{and}\quad p_{j}=\gamma \left( 1-\frac{1%
}{j}\right)
\end{equation*}
and the germ of a locally biholomorphic mapping $\phi =\varphi ^{-1}$ at the
point $p_{1}\equiv \gamma \left( 0\right) $ is analytically continued along
the subpath $\gamma [0,\tau ]$ with all $\tau <1$ as a locally biholomorphic
mapping. Let $p_{j}^{\prime }\equiv \phi \left( p_{j}\right) \in B^{n+1}$ be
the sequence obtained by the analytic continuation of the mapping $\phi
=\varphi ^{-1}$ along the path $\gamma [0,1).$ Then there is a subsequence $%
p_{m_{j}}^{\prime }$ and a point $p^{\prime }\in L\subset bB^{n+1}$ such
that the subsequence $p_{m_{j}}^{\prime }$ converges to the point $p^{\prime
}$ to a direction transversal to the complex tangent hyperplane at the point 
$p^{\prime }\in bB^{n+1}$.
\end{lemma}

\proof%
Since the closed ball $\overline{B^{n+1}}$ is compact, there is a
subsequence $p_{m_{j}}^{\prime }$ which converges to a point $p^{\prime }\in 
\overline{B^{n+1}}.$ Since the mapping $\varphi $ is locally biholomorphic
on $B^{n+1}$, the point $p^{\prime }$ should be on the boundary $bB^{n+1}.$
Further, since $\varphi =\phi ^{-1}:bB^{n+1}\backslash L\rightarrow bD$ is
locally biholomorphic, we obtain 
\begin{equation*}
p_{m_{j}}^{\prime }\rightarrow p^{\prime }\in L.
\end{equation*}
Let $\pi _{p}$ be the complex line containing the line segment $\gamma
[0,1]. $ Since $p_{m_{j}}^{\prime }\equiv \phi \left( p_{m_{j}}\right) ,$
the analytic continuation of the mapping $\phi =\varphi ^{-1}$ on the
complex line $\pi _{p}$ yields a complex curve $\varphi ^{-1}\left( \pi
_{p}\cap D\right) $ which touches the point $p^{\prime }\in L\subset
bB^{n+1}.$ Since the mapping $\varphi :B^{n+1}\rightarrow D$ is locally
biholomorphic up to the boundary subset $bB^{n+1}\backslash L$ and the
boundary $bB^{n+1}$ is strongly pseudoconvex, the extension 
\begin{equation*}
\varphi ^{-1}\left( \pi _{p}\cap \overline{D}\right)
\end{equation*}
of the complex curve $\varphi ^{-1}\left( \pi _{p}\cap D\right) $ is
transversal to the boundary $bB^{n+1}$ near the point $p^{\prime }\in
L\subset bB^{n+1}.$

We claim that the extension $\varphi ^{-1}\left( \pi _{p}\cap \overline{D}%
\right) $ is transversal to the complex tangent hyperplane at the point $%
p^{\prime }\in L.$ We take a path $\lambda :[0,1]\rightarrow bB^{n+1}$ such
that $p^{\prime }=\lambda (1)$ and 
\begin{equation*}
\lambda [0,1]\subset \overline{\varphi ^{-1}\left( \pi _{p}\cap bD\right) }%
\cap bB^{n+1}.
\end{equation*}
Note that there are finitely many closed paths $\gamma _{j}$ on $bD$ such
that 
\begin{equation*}
\pi _{p}\cap bD=\bigcup_{j}\gamma _{j}.
\end{equation*}
Since $bD$ is strongly pseudoconvex, each path $\gamma _{j}$ is transversal
to the complex tangent hyperplane on $bD$ at each point.

We consider the following sequence 
\begin{equation*}
q_{j}=\varphi \circ \lambda \left( 1-\frac{1}{j}\right) \in \pi _{p}\cap
bD,\quad q_{j}^{\prime }=\lambda \left( 1-\frac{1}{j}\right) \quad \text{for 
}j\in \Bbb{N}^{+}
\end{equation*}
such that 
\begin{equation*}
q_{j}^{\prime }\rightarrow p^{\prime }\quad \text{as }j\rightarrow \infty .
\end{equation*}
Since the set $\pi _{p}\cap bD$ is compact, there is a subsequence $%
q_{m_{j}} $ and a point $q\in \pi _{p}\cap bD$ such that 
\begin{equation*}
q_{m_{j}}\rightarrow q\quad \text{as }j\rightarrow \infty .
\end{equation*}
Let $\phi _{m_{j}}$ be the germ of the analytic continuation of the mapping $%
\phi =\varphi ^{-1}$ at the point $q_{m_{j}}.$ By Theorem \ref{anypath} and
Lemma \ref{localhull}, we may assume that there is an open neighborhood $W$
of the point $q$ such that the mapping $\phi _{m_{j}}$ is holomorphic on $%
W\cap D.$

Since $bD$ is spherical, there is a biholomorphic mapping $\psi $ on $W,$ if
necessary, shrinking $W,$ such that 
\begin{equation*}
\psi \in H\left( W\cap D\right) \cap C^{1}\left( W\cap \overline{D}\right)
\quad \text{and}\quad \psi \left( W\cap bD\right) \subset bB^{n+1}.
\end{equation*}
We may assume that 
\begin{equation*}
q^{\prime \prime }\equiv \psi \left( q\right) \neq p^{\prime }\quad \text{and%
}\quad q_{m_{j}}^{\prime \prime }\equiv \psi \left( q_{m_{j}}\right) .
\end{equation*}
Then, by Lemma \ref{sphere}, the composition 
\begin{equation*}
\eta _{j}\equiv \phi _{m_{j}}\circ \psi ^{-1}\in Aut\left( B^{n+1}\right)
\end{equation*}
are automorphisms of the unit ball $B^{n+1}.$ Note that we have the
condition 
\begin{equation*}
\eta _{j}\left( q_{m_{j}}^{\prime \prime }\right) =q_{m_{j}}^{\prime
}\rightarrow p^{\prime }
\end{equation*}
and the sequence $\eta _{j}$ converges to the point $p^{\prime }$ uniformly
on every compact subset of $B^{n+1},$ if necessary, passing to a
subsequence. Otherwise, $p^{\prime }\notin L.$

Let $\pi _{j}^{\prime }$ be the complex line passing through the two points $%
p^{\prime }$ and $q_{m_{j}}^{\prime },$ and $\delta _{j}^{\prime }$ be the
analytic disk 
\begin{equation*}
\delta _{j}^{\prime }=\pi _{j}^{\prime }\cap B^{n+1}.
\end{equation*}
The area $\left| \eta _{j}^{-1}\left( \delta _{j}^{\prime }\right) \right| $
of the analytic disk $\eta _{j}^{-1}\left( \delta _{j}^{\prime }\right) $ is
bounded from the below. Otherwise, we would have $p^{\prime }\notin L.$ Then
we decompose the mapping $\eta _{j}$ as follows 
\begin{equation*}
\eta _{j}=\mu _{j}\circ \nu _{j}
\end{equation*}
where 
\begin{equation*}
\mu _{j}\in Aut_{q_{m_{j}}^{\prime }}\left( bB^{n+1}\right) ,\quad \nu
_{j}\in Aut_{p^{\prime }}\left( bB^{n+1}\right) .
\end{equation*}
Then, by First and Second Scaling Lemmas, we take a sequence of local
automorphisms $\sigma _{j}\in Aut_{q_{m_{j}}^{\prime }}\left(
bB^{n+1}\right) $ defined by the normalizing parameters 
\begin{equation*}
U_{\sigma _{j}}=id_{n\times n},\quad \rho _{\sigma _{j}}=\rho _{\mu
_{j}},\quad a_{\sigma _{j}}=0,\quad r_{\sigma _{j}}=0.
\end{equation*}
Then the following sequence 
\begin{equation*}
\beta _{j}=\sigma _{j}^{-1}\circ \eta _{j}
\end{equation*}
converges to an automorphism of the unit ball $B^{n+1},$ if necessary,
passing to a subsequence.

Let $\pi _{j}$ be the complex line containing the tangent vector of the path 
\begin{equation*}
\psi \left( \varphi \circ \lambda [0,1]\cap W\right) \subset bB^{n+1}
\end{equation*}
at the point $q_{m_{j}}^{\prime \prime },$ and $\delta _{j}$ be the analytic
disk 
\begin{equation*}
\delta _{j}=\pi _{j}\cap B^{n+1}.
\end{equation*}
Clearly, the area $\left| \delta _{j}\right| $ of the analytic disk $\delta
_{j}$ is bounded from the below. Note that, by the mapping $\sigma _{j}\in
Aut_{q_{m_{j}}^{\prime }}\left( bB^{n+1}\right) ,$ the area $\left| \sigma
_{j}^{-1}\left( \delta _{j}^{\prime \prime }\right) \right| $ of the
analytic disk $\sigma _{j}\left( \delta _{j}^{\prime \prime }\right) $ is
bounded from the below whenever the area $\left| \delta _{j}^{\prime \prime
}\right| $ of the analytic disk $\delta _{j}^{\prime \prime }$ is bounded
from the below. Thus the area $\left| \eta _{j}\left( \delta _{j}\right)
\right| $ of the analytic disk 
\begin{equation*}
\eta _{j}\left( \delta _{j}\right) =\sigma _{j}\circ \beta _{j}\left( \delta
_{j}\right)
\end{equation*}
is bounded from the below. Since the analytic disk $\eta _{j}\left( \delta
_{j}\right) $ is the intersection of $B^{n+1}$ and the complex line
containing the tangent vector of the path $\lambda [0,1]$ at the point 
\begin{equation*}
q_{m_{j}}^{\prime }=\lambda \left( 1-\frac{1}{m_{j}}\right) .
\end{equation*}
Thus the path $\lambda [0,1]$ is transversal to the complex tangent
hyperplane at the point $p^{\prime }=\lambda \left( 1\right) $. Therefore,
the point $p_{j}^{\prime }$ approaches to the point $p^{\prime }\in L$ to a
direction transversal to the complex tangent hyperplane at the point $%
p^{\prime }\in bB^{n+1}.$ This completes the proof.%
\endproof%

\begin{lemma}[Second Dogginal Lemma]
\label{doggi}Let $D$ be a bounded domain in $\Bbb{C}^{n+1}$ with spherical
differentiable boundary $bD$ and $\phi $ be a biholomorphic mapping on an
open neighborhood $U$ of a point $r\in bD$ satisfying 
\begin{equation*}
\phi \in H\left( U\cap D\right) \cap C^{1}\left( U\cap \overline{D}\right)
\quad \text{and}\quad \phi \left( U\cap bD\right) \subset bB^{n+1}.
\end{equation*}
Let $L$ be the Doggaebi variety associated to the mapping $\phi $ and $%
\varphi :B^{n+1}\rightarrow D$ be a locally biholomorphic mapping to be an
analytic continuation of the mapping $\phi ^{-1}$. Suppose that there is a
line segment $\gamma :[0,1]\rightarrow D$ with 
\begin{equation*}
p=\gamma \left( 1\right) \quad \text{and}\quad p_{j}=\gamma \left( 1-\frac{1%
}{j}\right)
\end{equation*}
and the germ of a locally biholomorphic mapping $\phi =\varphi ^{-1}$ at the
point $p_{1}\equiv \gamma \left( 0\right) $ is analytically continued along
the subpath $\gamma [0,\tau ]$ with all $\tau <1$ as a locally biholomorphic
mapping. Let $p_{j}^{\prime }\equiv \phi \left( p_{j}\right) \in B^{n+1}$ be
the sequence obtained by the analytic continuation of the mapping $\phi
=\varphi ^{-1}$ along the path $\gamma [0,1).$ Suppose that there is a point 
$p^{\prime }\in L\subset bB^{n+1}$ such that the sequence $p_{j}^{\prime }$
converges to the point $p^{\prime }$ to a direction transversal to the
complex tangent hyperplane at the point $p^{\prime }\in bB^{n+1}$. Let $\pi
_{p}$ be the complex line containing the line segment $\gamma [0,1]$. Then
there is a distinguished complex hyperplane $H_{p}\subset T_{p}\Bbb{C}^{n+1}$
at the point $p\in D$ satisfying 
\begin{equation*}
T_{p}\Bbb{C}^{n+1}=H_{p}\oplus \pi _{p}
\end{equation*}
and an open neighborhood $V$ of the line segment $\gamma [0,1)$ such that

\begin{enumerate}
\item  $bV\cap \pi _{p}\cap B\left( p;\delta \right) $ is an angle for a
sufficient small $\delta >0,$ which contains the path $\gamma [0,1],$

\item  $bV\cap R_{p}\cap B\left( p;\delta \right) $ is a slanted paraboloid
for a sufficiently small $\delta >0,$ which smoothly touches the complex
hyperplane $H_{p}$ at the point $p,$

\item  the germ of a locally biholomorphic mapping $\phi =\varphi ^{-1}$ at
the point $p_{1}\equiv \gamma \left( 0\right) $ is analytically continued on 
$V$ as a locally biholomorphic mapping,
\end{enumerate}

\noindent where $R_{p}$ is a real hyperplane containing the line segment $%
\gamma [0,1]$ and 
\begin{equation*}
H_{p}\subset R_{p}.
\end{equation*}
\end{lemma}

\proof%
We take a sequence of automorphisms $\phi _{j}\in Aut\left( B^{n+1}\right) $
such that 
\begin{equation*}
\phi _{j}\left( p_{1}^{\prime }\right) =p_{j}^{\prime }\rightarrow p^{\prime
}\in L.
\end{equation*}
Then the composition $\tau _{j}\equiv \varphi \circ \phi
_{j}:B^{n+1}\rightarrow D$ forms a normal family so that there is a
subsequence $\tau _{m_{j}}\equiv \varphi \circ \phi _{m_{j}}$ which
converges to a holomorphic mapping $\tau :B^{n+1}\rightarrow D$ uniformly on
every compact subset of $B^{n+1}.$

We claim that $\tau $ is a locally biholomorphic mapping. We take a point $%
q\in bB^{n+1}\backslash L$ and $q_{m_{j}}^{\prime }=\tau _{m_{j}}\left(
q\right) \in bD$ such that 
\begin{equation*}
q_{m_{j}}^{\prime }=\tau _{m_{j}}\left( q\right) \rightarrow q^{\prime }\in
bD,
\end{equation*}
if necessary, passing to a subsequence. Since $bD$ is spherical, there is an
open neighborhood $W$ of the point $q^{\prime }$ and a biholomorphic mapping 
$\psi $ on $W\cap D$ satisfying 
\begin{equation*}
\psi \in H\left( W\cap D\right) \cap C^{1}\left( W\cap \overline{D}\right)
\quad \text{and}\quad \psi \left( W\cap bD\right) \subset bB^{n+1}.
\end{equation*}
Then we consider the composition $\chi _{m_{j}}\equiv \psi \circ \varphi
\circ \phi _{m_{j}}$ so that, by Lemma \ref{sphere} and by abuse of
notation, 
\begin{equation*}
\chi _{m_{j}}\in Aut\left( B^{n+1}\right) .
\end{equation*}
Suppose that the assertion is not true. Then the sequence $\chi _{m_{j}}$
converges to a boundary point $\psi \left( q^{\prime }\right) \in bB^{n+1}$
uniformly on every compact subset of $B^{n+1},$ if necessary, passing to a
subsequence. Hence the sequence $\tau _{m_{j}}\equiv \varphi \circ \phi
_{m_{j}}$ would converge to a boundary point $q^{\prime }\in bD$ uniformly
on every compact subset of $B^{n+1}$. By the way, note that 
\begin{equation*}
\tau _{m_{j}}\left( p_{1}^{\prime }\right) =p_{m_{j}}\rightarrow p\in D.
\end{equation*}
This is a contradiction. Hence the holomorphic mapping $\tau $ is locally
biholomorphic. Thus there is a real number $\delta >0$ such that the mapping 
$\tau _{m_{j}}$ and the inverse $\tau _{m_{j}}^{-1}$ are analytically
continued to a biholomorphic mapping respectively on 
\begin{equation*}
B\left( p_{m_{j}};\delta \right) ,\quad B\left( p_{1}^{\prime };\delta
\right) .
\end{equation*}

Let $\varepsilon _{j}$ be the euclidean length between the two points $%
p^{\prime }$ and $p_{j}^{\prime }$. By First Scaling Lemma, there is a
sequence of automorphisms 
\begin{equation*}
\sigma _{j}\in Aut_{p^{\prime }}\left( bB^{n+1}\right)
\end{equation*}
such that, for a subsequence $\phi _{m_{j}},$%
\begin{equation*}
U_{\sigma _{j}}=id_{n\times n},\quad \rho _{\sigma _{j}}=\varepsilon
_{m_{j}},\quad a_{\sigma _{j}}=0,\quad r_{\sigma _{j}}=0
\end{equation*}
and the composition 
\begin{equation*}
\eta _{j}\equiv \sigma _{j}^{-1}\circ \phi _{m_{j}}
\end{equation*}
uniformly converges to an automorphism of the unit ball $B^{n+1}.$ Then
there is a real number $\delta >0$ such that the mapping $\eta _{j}$ and the
inverse $\eta _{j}^{-1}$ are analytically continued respectively on 
\begin{equation*}
B\left( p_{1}^{\prime };\delta \right) ,\quad B\left( p_{m_{j}}^{\prime
\prime };\delta \right)
\end{equation*}
where 
\begin{equation*}
p_{m_{j}}^{\prime \prime }=\eta _{j}\left( p_{m_{j}}^{\prime }\right) .
\end{equation*}

Then we obtain 
\begin{equation*}
\varphi =\tau _{m_{j}}\circ \eta _{j}^{-1}\circ \sigma _{j}^{-1}.
\end{equation*}
Thus the mapping $\phi =\varphi ^{-1}$ is analytically continued onto the
open set 
\begin{equation*}
\tau _{m_{j}}\circ \eta _{j}^{-1}\left( \sigma _{j}\left( B\left(
p_{m_{j}}^{\prime \prime };\delta \right) \right) \right)
\end{equation*}
which is centered at the point $p_{m_{j}}.$ For a canonical normalizing
mapping $\mu _{p^{\prime }},$ we obtain 
\begin{equation*}
\sigma _{j}^{\prime }\equiv \mu _{p^{\prime }}\circ \sigma _{j}\circ \mu
_{p^{\prime }}^{-1}:\left\{ 
\begin{array}{l}
z^{*}=\sqrt{\varepsilon _{m_{j}}}z \\ 
w^{*}=\varepsilon _{m_{j}}w
\end{array}
\right. .
\end{equation*}
We set 
\begin{equation*}
\mu _{p^{\prime }}\left( p_{m_{j}}^{\prime \prime }\right) \rightarrow
p^{\prime \prime \prime }\in \mu _{p^{\prime }}\left( B^{n+1}\right)
\end{equation*}
so that, shrinking $\delta >0,$ if necessary, the mapping $\phi =\varphi
^{-1}$ is analytically continued onto the open set 
\begin{equation*}
\tau _{m_{j}}\circ \eta _{j}^{-1}\circ \mu _{p^{\prime }}\left( \sigma
_{j}^{\prime }\left( B\left( p^{\prime \prime \prime };\delta \right)
\right) \right)
\end{equation*}
which is centered at the point $p_{m_{j}}.$

Since $\varepsilon _{m_{j}}$ is the euclidean distance between the two
points $p$ and $p_{m_{j}},$ the analytically continued region of the mapping 
$\phi =\varphi ^{-1}$ along the line segment $\gamma [0,1]\subset D$
contains an open set along the line segment $\gamma [0,1]$ which touches to
the point $p=\gamma \left( 1\right) $ by an edge shape on the complex line $%
\pi _{p}$ and by a slanted paraboloid shape on the real hyperplane $R_{p}.$
This completes the proof.%
\endproof%

\begin{lemma}
Let $D$ be a bounded domain in $\Bbb{C}^{n+1}$ with spherical differentiable
boundary $bD$ and $\phi $ be a biholomorphic mapping on an open neighborhood 
$U$ of a point $r\in bD$ satisfying 
\begin{equation*}
\phi \in H\left( U\cap D\right) \cap C^{1}\left( U\cap \overline{D}\right)
\quad \text{and}\quad \phi \left( U\cap bD\right) \subset bB^{n+1}.
\end{equation*}
Let $L$ be the Doggaebi variety associated to the mapping $\phi $ and $%
\varphi :B^{n+1}\rightarrow D$ be a locally biholomorphic mapping to be an
analytic continuation of the mapping $\phi ^{-1}$. Let $H_{p}$ be the
complex hyperplane at the point $p\in D$ in Lemma \ref{doggi} and $\pi _{p}$
be a complex line satisfying 
\begin{equation*}
T_{p}\Bbb{C}^{n+1}=H_{p}\oplus \pi _{p}.
\end{equation*}
Suppose that, along a line segment $\gamma :[0,1]\rightarrow \pi _{p}\cap D$
with $p=\gamma \left( 1\right) ,$ the germ of a locally biholomorphic
mapping $\phi =\varphi ^{-1}$ at the point $\gamma \left( 0\right) $ is
analytically continued along the subpath $\gamma [0,\tau ]$ with all $\tau
<1 $ as a locally biholomorphic mapping such that the limit point 
\begin{equation*}
\lim_{\tau \rightarrow 1}\phi \circ \gamma \left( \tau \right)
\end{equation*}
is the point $p^{\prime }\in L\subset bB^{n+1}$ in Lemma \ref{doggi}. Then
there is an open neighborhood $V$ of the line segment $\gamma [0,1)$ such
that

\begin{enumerate}
\item  $bV\cap \pi _{p}\cap B\left( p;\delta \right) $ is an angle for a
sufficient small $\delta >0,$ which contains the path $\gamma [0,1],$

\item  $bV\cap R_{p}\cap B\left( p;\delta \right) $ is a slanted paraboloid
for a sufficiently small $\delta >0,$ which smoothly touches the complex
hyperplane $H_{p}$ at the point $p,$

\item  the germ of a locally biholomorphic mapping $\phi =\varphi ^{-1}$ at
the point $\gamma \left( 0\right) $ is analytically continued on $V$ as a
locally biholomorphic mapping,
\end{enumerate}

\noindent where $R_{p}$ is a real hyperplane containing the line segment $%
\gamma [0,1]$ and 
\begin{equation*}
H_{p}\subset R_{p}.
\end{equation*}
\end{lemma}

\proof%
By Lemma \ref{transversally}, $\phi \circ \gamma \left( \tau \right) $
approaches to the point $p^{\prime }\in L$ to a direction transversal to the
complex tangent hyperplane at the point $p^{\prime }$ as $\tau \rightarrow
1. $ Then there is a complex hyperplane $H_{p}^{\prime }$ satisfying the
conditions in Lemma \ref{doggi}.

Note that the complex hyperplane $H_{p}$ is determined by the mapping $%
\varphi :B^{n+1}\rightarrow D$ and the complex tangent hyperplane $%
H_{p^{\prime }}$ at the point $p^{\prime }\in L\subset bB^{n+1},$ but $H_{p}$
does not depend on the approaching direction of $\phi \circ \gamma \left(
\tau \right) \rightarrow p^{\prime }$ as $\tau \rightarrow 1.$ Thus $%
H_{p}^{\prime }=H_{p}$. This completes the proof.%
\endproof%

\subsection{Biholomorphic equivalence}

\begin{lemma}
\label{lifting}Let $D$ be a bounded domain in $\Bbb{C}^{n+1}$ with spherical
differentiable boundary $bD$ and $G$ be the path space of $D$ mod homotopic
relation so that $G$ is a universal covering Riemann domain with a natural
complex structure and a natural holomorphic covering map $\kappa
:G\rightarrow D.$ Suppose that there is a biholomorphic mapping 
\begin{equation*}
\phi \in H\left( U\cap D\right) \cap C^{1}\left( U\cap \overline{D}\right)
\end{equation*}
for a connected open neighborhood $U$ of a point $p\in bD$ satisfying 
\begin{equation*}
\phi \left( U\cap bD\right) \subset bB^{n+1}.
\end{equation*}
Then the analytic continuation of the inverse mapping $\phi ^{-1}$ is lifted
to a holomorphic mapping $\psi $ from $B^{n+1}$ onto $G$ such that $\psi
:B^{n+1}\rightarrow G$ is a proper locally biholomorphic mapping, i.e., a
finite covering of $G,$ satisfying the following relation: 
\begin{equation*}
\begin{array}{lcl}
&  & G \\ 
& \overset{\psi }{\nearrow } & \downarrow \kappa \\ 
B^{n+1} & \overset{\varphi }{\longrightarrow } & D
\end{array}
.
\end{equation*}
\end{lemma}

\proof%
By Theorem \ref{branched}, there is a locally biholomorphic mapping 
\begin{equation*}
\varphi :B^{n+1}\rightarrow D
\end{equation*}
such that $\varphi $ is the analytic continuation of the mapping $\phi ^{-1}$
and 
\begin{equation*}
\varphi =\phi ^{-1}\quad \text{on}\quad bB^{n+1}\backslash L
\end{equation*}
where $L$ is the Doggaebi variety associated to the mapping $\phi .$

A piecewise line segment curve is a path $\gamma :[0,1]\rightarrow D$
consisting of finitely many line segments. Fix a point $b\in U\cap D.$ We
shall show that the mapping $\phi $ is analytically continued along any
piecewise line segment curve. Suppose that $\gamma :[0,1]\rightarrow D$ be a
piecewise line segment curve with $\gamma \left( 0\right) =b$ such that the
germ of the mapping $\phi $ at the point $\gamma \left( 0\right) $ is
analytically continued along all subpath $\gamma [0,\tau ]$ with $\tau <1,$
but not the whole path $\gamma [0,1].$ Let $\pi _{p}$ be a complex line
satisfying, for sufficiently small $\varepsilon >0,$%
\begin{equation*}
\gamma [1-\varepsilon ,1]\subset \pi _{p}.
\end{equation*}
Then, by Second Dogginal Lemma, there is a distinguished complex hyperplane $%
H_{p}$ at the point $p=\gamma \left( 1\right) $ satisfying 
\begin{equation*}
T_{p}\Bbb{C}^{n+1}=H_{p}\oplus \pi _{p}
\end{equation*}
and an open neighborhood $V$ containing the line segment $\gamma [0,\tau ]$
with all $\tau <1$ such that

\begin{enumerate}
\item  $bV\cap \pi _{p}\cap B\left( p;\varepsilon \right) $ is an angle for
a sufficient small $\varepsilon >0,$ which contains the path $\gamma [0,1],$

\item  $bV\cap R_{p}\cap B\left( p;\varepsilon \right) $ is a slanted
paraboloid for a sufficiently small $\varepsilon >0,$ which smoothly touches
the complex hyperplane $H_{p}$ at the point $p,$

\item  the germ of a locally biholomorphic mapping $\phi =\varphi ^{-1}$ at
the point $\gamma \left( 0\right) $ is analytically continued on $V$ as a
locally biholomorphic mapping,
\end{enumerate}

\noindent where $R_{p}$ is a real hyperplane satisfying, for sufficiently
small $\varepsilon >0,$%
\begin{equation*}
\gamma [1-\varepsilon ,1]\subset R_{p}\quad \text{and}\quad H_{p}\subset
R_{p}.
\end{equation*}

Thus we can find a sequence of points $p_{j}\in V\cap B\left( p;\varepsilon
\right) $ and a sequence of germs $\phi _{j}$ of the analytic continuation
of the mapping $\phi $ at the point $p_{j}$ such that the sequence $p_{j}$
converges to the point $p$ to a direction tangential to the complex
hyperplane $H_{p}$ at the point $p.$ We set 
\begin{equation*}
p_{j}^{\prime }=\phi _{j}\left( p_{j}\right) \in B^{n+1}.
\end{equation*}
Then there is a point $p^{\prime }\in L\subset bB^{n+1}$ such that $%
p_{j}^{\prime }\rightarrow p^{\prime }$ to a direction tangential to the
complex tangent hyperplane at the point $p^{\prime }\in bB^{n+1}.$ We take a
sequence of automorphisms $\varphi _{j}\in Aut\left( B^{n+1}\right) $ such
that 
\begin{equation*}
\varphi _{j}\left( p_{1}^{\prime }\right) =p_{j}^{\prime }\rightarrow
p^{\prime }\in L.
\end{equation*}
Then the composition $\tau _{j}\equiv \varphi \circ \varphi
_{j}:B^{n+1}\rightarrow D$ forms a normal family so that there is a
subsequence $\tau _{m_{j}}\equiv \varphi \circ \varphi _{m_{j}}$ which
converges to a holomorphic mapping $\tau :B^{n+1}\rightarrow D$ uniformly on
every compact subset of $B^{n+1}.$

We claim that $\tau $ is a locally biholomorphic mapping. We take a point $%
q\in bB^{n+1}\backslash L$ and $q_{m_{j}}^{\prime }=\tau _{m_{j}}\left(
q\right) \in bD$ such that 
\begin{equation*}
q_{m_{j}}^{\prime }=\tau _{m_{j}}\left( q\right) \rightarrow q^{\prime }\in
bD,
\end{equation*}
if necessary, passing to a subsequence. Since $bD$ is spherical, there is an
open neighborhood $W$ of the point $q^{\prime }$ and a biholomorphic mapping 
$\psi $ on $W\cap D$ satisfying 
\begin{equation*}
\psi \in H\left( W\cap D\right) \cap C^{1}\left( W\cap \overline{D}\right)
\quad \text{and}\quad \psi \left( W\cap bD\right) \subset bB^{n+1}.
\end{equation*}
Then we consider the composition $\chi _{m_{j}}\equiv \psi \circ \varphi
\circ \varphi _{m_{j}}$ so that, by Lemma \ref{sphere} and by abuse of
notation, 
\begin{equation*}
\chi _{m_{j}}\in Aut\left( B^{n+1}\right) .
\end{equation*}
Suppose that the assertion is not true. Then the sequence $\chi _{m_{j}}$
converges to a boundary point $\psi \left( q^{\prime }\right) \in bB^{n+1}$
uniformly on every compact subset of $B^{n+1},$ if necessary, passing to a
subsequence. Hence the sequence $\tau _{m_{j}}\equiv \varphi \circ \varphi
_{m_{j}}$ would converge to a boundary point $q^{\prime }\in bD$ uniformly
on every compact subset of $B^{n+1}$. By the way, note that 
\begin{equation*}
\tau _{m_{j}}\left( p_{1}^{\prime }\right) =p_{m_{j}}\rightarrow p\in D.
\end{equation*}
This is a contradiction. Hence the holomorphic mapping $\tau $ is locally
biholomorphic. Thus there is a real number $\delta >0$ such that the mapping 
$\tau _{m_{j}}$ and the inverse $\tau _{m_{j}}^{-1}$ are analytically
continued to a biholomorphic mapping respectively on 
\begin{equation*}
B\left( p_{m_{j}};\delta \right) ,\quad B\left( p_{1}^{\prime };\delta
\right) .
\end{equation*}

Let $\varepsilon _{j}$ be the euclidean length between the two points $%
p^{\prime }$ and $p_{j}^{\prime }$. By Second Scaling Lemma, there is a
sequence of automorphisms 
\begin{equation*}
\sigma _{j}\in Aut_{p^{\prime }}\left( bB^{n+1}\right)
\end{equation*}
such that, for a subsequence $\varphi _{m_{j}},$%
\begin{equation*}
U_{\sigma _{j}}=id_{n\times n},\quad \rho _{\sigma _{j}}=\varepsilon
_{m_{j}}^{2},\quad a_{\sigma _{j}}=0,\quad r_{\sigma _{j}}=0
\end{equation*}
and the composition 
\begin{equation*}
\eta _{j}\equiv \sigma _{j}^{-1}\circ \varphi _{m_{j}}
\end{equation*}
uniformly converges to an automorphism of the unit ball $B^{n+1}.$ Then
there is a real number $\delta >0$ such that the mapping $\eta _{j}$ and the
inverse $\eta _{j}^{-1}$ are analytically continued respectively on 
\begin{equation*}
B\left( p_{1}^{\prime };\delta \right) ,\quad B\left( p_{m_{j}}^{\prime
\prime };\delta \right)
\end{equation*}
where 
\begin{equation*}
p_{m_{j}}^{\prime \prime }=\eta _{j}\left( p_{m_{j}}^{\prime }\right) .
\end{equation*}

Then we obtain 
\begin{equation*}
\varphi =\tau _{m_{j}}\circ \eta _{j}^{-1}\circ \sigma _{j}^{-1}.
\end{equation*}
Thus the mapping $\phi =\varphi ^{-1}$ is analytically continued onto the
open set 
\begin{equation*}
\tau _{m_{j}}\circ \eta _{j}^{-1}\left( \sigma _{j}\left( B\left(
p_{m_{j}}^{\prime \prime };\delta \right) \right) \right)
\end{equation*}
which is centered at the point $p_{m_{j}}.$ For a canonical normalizing
mapping $\mu _{p^{\prime }},$ we obtain 
\begin{equation*}
\sigma _{j}^{\prime }\equiv \mu _{p^{\prime }}\circ \sigma _{j}\circ \mu
_{p^{\prime }}^{-1}:\left\{ 
\begin{array}{l}
z^{*}=\varepsilon _{m_{j}}z \\ 
w^{*}=\varepsilon _{m_{j}}^{2}w
\end{array}
\right. .
\end{equation*}
We set 
\begin{equation*}
\mu _{p^{\prime }}\left( p_{m_{j}}^{\prime \prime }\right) \rightarrow
p^{\prime \prime \prime }\in \mu _{p^{\prime }}\left( B^{n+1}\right)
\end{equation*}
so that, shrinking $\delta >0,$ if necessary, the mapping $\phi =\varphi
^{-1}$ is analytically continued onto the open set 
\begin{equation*}
\tau _{m_{j}}\circ \eta _{j}^{-1}\circ \mu _{p^{\prime }}\left( \sigma
_{j}^{\prime }\left( B\left( p^{\prime \prime \prime };\delta \right)
\right) \right)
\end{equation*}
which is centered at the point $p_{m_{j}}.$

Since the sequence $p_{m_{j}}$ converges to $p$ to a direction tangential to
the complex hyperplane $H_{p}$ at the point $p,$ and $\varepsilon _{m_{j}}$
is the euclidean distance between the two points $p$ and $p_{m_{j}},$ the
analytically continued region of the mapping $\phi =\varphi ^{-1}$ along the
sequence $p_{m_{j}}\in V\cap D$ contains an open set along the converging
direction of the sequence $p_{m_{j}}$ by an edge shape tangential to the
complex hyperplane $H_{p}$ at the point $p$ and a $\sqrt{\left| x\right| }$
curve shape normal to $H_{p}.$

Thus there is a complex line $\pi _{p}^{\prime }$ satisfying 
\begin{equation*}
\pi _{p}^{\prime }\subset H_{p}\quad \text{and}\quad p\in \pi _{p}^{\prime },
\end{equation*}
and an open neighborhood $W$ of the point $p$ such that the mapping $\phi $
is analytically continued on 
\begin{equation*}
W\cap \pi _{p}^{\prime }.
\end{equation*}
Then we take a line segment $\gamma :[0,1]\rightarrow \pi _{p}^{\prime
}\subset H_{p}$ such that 
\begin{equation*}
\gamma :[0,1)\subset W\cap \left( \pi _{p}^{\prime }\backslash p\right)
\quad \text{and}\quad \gamma \left( 1\right) =p.
\end{equation*}
Note that 
\begin{equation*}
\lim_{\tau \rightarrow 1}\phi \circ \gamma \left( \tau \right) =p^{\prime
}\in L\subset bB^{n+1}.
\end{equation*}
By Second Dogginal Lemma, the distinguished complex tangent hyperplane $%
H_{p} $ satisfies 
\begin{equation*}
T_{p}\Bbb{C}^{n+1}=H_{p}\oplus \pi _{p}^{\prime }
\end{equation*}
so that 
\begin{equation*}
\pi _{p}^{\prime }\cap H_{p}=\left\{ p\right\} .
\end{equation*}
This is a contradiction so that the germ of the mapping $\phi $ at the point 
$b\in U\cap D$ is analytically continued along the whole path $\gamma [0,1].$
Hence the mapping $\phi $ is analytically continued along any piecewise line
segment curve on $D$ as a locally biholomorphic mapping.

We claim that the mapping $\phi :U\cap D\rightarrow B^{n+1}$ is analytically
continued along any path on $D$ as a locally biholomorphic mapping. Let $%
\gamma :[0,1]\rightarrow D$ be a path on $D.$ Then we take a piecewise line
segment curve $\lambda :[0,1]$ and a continuous function $\Gamma
:[0,1]\times [0,1]\rightarrow D$ such that 
\begin{eqnarray*}
\gamma \left( 0\right) &=&\lambda \left( 0\right) ,\quad \gamma \left(
1\right) =\lambda \left( 1\right) \\
\Gamma \left( 0,\tau \right) &=&\gamma \left( \tau \right) ,\quad \Gamma
\left( 1,\tau \right) =\lambda \left( \tau \right) \quad \text{for all }\tau
\in [0,1]
\end{eqnarray*}
and the path $\Gamma \left( \cdot ,\tau \right) :[0,1]\rightarrow D$ for
each $\tau \in [0,1]$ is a piecewise line segment curve on $D.$ Clearly, the
mapping $\phi $ is analytically continued along the whole path $\gamma [0,1]$
so that the mapping $\phi $ is analytically continued along any path on $D$
as a locally biholomorphic mapping.

Let $G$ be the path space of $D$ pointed by a point of $U\cap D$ mod
homotopic relation so that $G$ is a universal covering of $D$ with a natural
complex structure and a natural holomorphic covering map $\kappa
:G\rightarrow D.$ By the analytic continuation of the mapping $\phi :U\cap
D\rightarrow B^{n+1}$ on $D,$ the mapping $\varphi :B^{n+1}\rightarrow D$
has its natural proper locally biholomorphic lift $\psi :B^{n+1}\rightarrow
G.$ Since the mapping $\psi :B^{n+1}\rightarrow G$ is proper and locally
biholomorphic, the mapping $\psi $ is a finite covering of $G$. This
completes the proof.%
\endproof%

\begin{theorem}
Let $D$ be a bounded domain in $\Bbb{C}^{n+1}$ with spherical differentiable
boundary $bD$ such that the fundamental group $\pi _{1}(D)$ is finite.
Suppose that there is a biholomorphic mapping 
\begin{equation*}
\phi \in H\left( U\cap D\right) \cap C^{1}\left( U\cap \overline{D}\right)
\end{equation*}
for a connected open neighborhood $U$ of a point $p\in bD$ satisfying 
\begin{equation*}
\phi \left( U\cap bD\right) \subset bB^{n+1}.
\end{equation*}
Then $D$ is necessarily simply connected and the mapping $\phi $ is
analytically continued to a biholomorphic mapping from $D$ onto $B^{n+1}.$
\end{theorem}

\proof%
Let $G$ be the path space of $D$ pointed by a point of $U\cap D$ mod
homotopic relation with a natural complex structure and a natural locally
biholomorphic covering map $\kappa :G\rightarrow D.$ Then, by Lemma \ref
{lifting}, there is a locally biholomorphic finite covering lift $\psi
:B^{n+1}\rightarrow G$ satisfying 
\begin{equation*}
\varphi =\kappa \circ \psi
\end{equation*}
where $\varphi :B^{n+1}\rightarrow D$ is an analytic continuation of the
inverse mapping $\phi ^{-1}.$

Since the fundamental group $\pi _{1}\left( D\right) $ is finite, the
mapping $\kappa :G\rightarrow D$ is a finite covering. Thus the mapping 
\begin{equation*}
\varphi =\kappa \circ \psi :B^{n+1}\rightarrow D
\end{equation*}
is a locally biholomorphic finite covering. Therefore, the analytic
continuation of the inverse mapping $\phi ^{-1}$ on the boundary $bB^{n+1}$
yields finitely many germs at each point of $bB^{n+1}.$ By Lemma \ref{germs}%
, the Doggaebi variety $L$ associated to the mapping $\phi $ is empty so
that the mapping 
\begin{equation*}
\varphi =\phi ^{-1}:bB^{n+1}\rightarrow bD
\end{equation*}
is also a locally biholomorphic finite covering. Thus the mapping 
\begin{equation*}
\varphi :\overline{B^{n+1}}\rightarrow \overline{D}
\end{equation*}
is a finite covering mapping. Then the fixed point property of the close
ball $\overline{B^{n+1}}$ implies that the mapping $\varphi :\overline{%
B^{n+1}}\rightarrow \overline{D}$ is a simple cover. Otherwise, a nontrivial
deck transformation of the closed ball $\overline{B^{n+1}}$ yields a
continuous function on $\overline{B^{n+1}}$ with no fixed point. Therefore,
the analytic continuation of the mapping $\phi :U\cap D\rightarrow B^{n+1}$
is analytically continued to a biholomorphic mapping 
\begin{equation*}
\varphi ^{-1}:D\rightarrow B^{n+1}.
\end{equation*}
This completes the proof.%
\endproof%

\section{Locally Biholomorphic Mappings}

\subsection{Estimates of normalizing parameters}

\begin{lemma}
\label{est1}Let $M$ be a nonspherical analytic real hypersurface in $\Bbb{C}%
^{n+1}$ and $Aut_{p}(M)$ be the isotropy subgroup at a point $p\in M.$
Suppose that there is a real number $c\geq 1$ satisfying 
\begin{equation*}
\sup_{\varphi \in Aut_{p}(M)}\left| U_{\varphi }\right| \leq c<\infty .
\end{equation*}
Then there is a real number $e\geq 1$ satisfying 
\begin{equation*}
\left| a_{\phi }\right| \leq e,\quad e^{-1}\leq \left| \rho _{\phi }\right|
\leq e,\quad \left| r_{\phi }\right| \leq e
\end{equation*}
for every local automorphism $\phi \in Aut_{p}(M).$
\end{lemma}

\proof%
The real hypersurface $\mu _{p}\left( M\right) $ in normal form is expanded
as follows: 
\begin{equation*}
v=\langle z,z\rangle +\sum_{k=4}^{\infty }F_{k}\left( z,\overline{z},u\right)
\end{equation*}
where 
\begin{equation*}
F_{k}\left( \nu z,\nu \overline{z},\nu ^{2}u\right) =\nu ^{k}F_{k}\left( z,%
\overline{z},u\right) .
\end{equation*}
Since $M$ is nonspherical, not all $F_{k}\left( z,\overline{z},u\right) $
are zero. Then we make an estimate of the normalizing parameters $a_{\phi
},\rho _{\phi },r_{\phi }$ for $\phi \in Aut_{p}(M)$(cf.\cite{Pa3}). This
completes the proof.%
\endproof%

\begin{theorem}
\label{strongly}Let $M$ be a nonspherical strongly pseudoconvex analytic
real hypersurface in a complex manifold. Then the local isotropy subgroup $%
Aut_{p}\left( M\right) $ is compact for every point $p\in M.$
\end{theorem}

\proof%
Because the situation is local, we may assume by taking a coordinate chart,
if necessary, that the complex manifold is $\Bbb{C}^{n+1}.$ Then the
pseudoconvexity of $M$ leads to 
\begin{equation*}
\left| U_{\varphi }\right| =1\quad \text{for}\quad \varphi \in Aut_{p}\left(
M\right) .
\end{equation*}
By Lemma \ref{est1}, we obtain the desired result. This completes the proof.%
\endproof%

\begin{lemma}
\label{germ-biho}Let $M,$ $M^{\prime }$ be nonspherical analytic real
hypersurfaces in $\Bbb{C}^{n+1}$ and $p,p^{\prime }$ be points respectively
of $M,M^{\prime }$ such that the two germs $M$ at $p$ and $M^{\prime }$ at $%
p^{\prime }$ are biholomorphically equivalent. Suppose that the isotropy
subgroup $Aut_{p}\left( M\right) $ is compact. Then there is a real number $%
\delta _{p}>0$ such that each germ of a biholomorphic mapping $\phi $
sending the germ $M$ at $p$ to the germ $M^{\prime }$ at $p^{\prime }$ is
analytically continued to the open ball $B(p;\delta _{p}).$
\end{lemma}

\proof%
We take a biholomorphic mapping $\phi $ on a connected open neighborhood $U$
of the point $p\in M$ such that 
\begin{equation*}
\phi \left( U\cap M\right) \subset M^{\prime }.
\end{equation*}
Then every germ of a biholomorphic mapping sending the germ $M$ at $p$ to
the germ $M^{\prime }$ at $p^{\prime }$ is one of the following 
\begin{equation*}
\phi \circ \varphi \quad \text{for}\quad \varphi \in Aut_{p}\left( M\right) .
\end{equation*}
Then the compactness of the group $Aut_{p}\left( M\right) $ leads to the
desired result. This completes the proof.%
\endproof%

\begin{lemma}
\label{uniform1}Let $M$ be a nonspherical analytic real hypersurface in $%
\Bbb{C}^{n+1}$ such that the isotropy subgroup $Aut_{p}(M)$ is compact at
every point $p\in M.$ Then, for each compact subset $K\subset \subset M,$
there is a real number $e\geq 1$ satisfying 
\begin{equation*}
e^{-1}\leq \left| U_{\phi }\right| \leq e,\quad \left| a_{\phi }\right| \leq
e,\quad e^{-1}\leq \left| \rho _{\phi }\right| \leq e,\quad \left| r_{\phi
}\right| \leq e
\end{equation*}
for every point $p\in K$ and every local automorphism $\phi \in Aut_{p}(M).$
\end{lemma}

\proof%
The real hypersurface $\mu _{p}\left( M\right) $ in normal form is expanded
as follows: 
\begin{equation*}
v=\langle z,z\rangle +\sum_{\left| I\right| ,\left| J\right| \geq 2,k\geq
1}\lambda _{IJk}\left( p\right) z^{I}\overline{z}^{J}u^{k}.
\end{equation*}
Suppose that the normalization $N_{e},$ $e=\left( U,a,\rho ,r\right) \in H,$
transforms the real hypersurface $\mu _{p}\left( M\right) $ to a real
hypersurfaces expanded as follows: 
\begin{equation*}
v=\langle z,z\rangle +\sum_{\left| I\right| ,\left| J\right| \geq 2,k\geq
1}\eta _{IJk}\left( p;U,a,\rho ,r\right) z^{I}\overline{z}^{J}u^{k}.
\end{equation*}
Note that the element $\left( U_{\phi },a_{\phi },\rho _{\phi },r_{\phi
}\right) $ for $\phi \in Aut_{p}\left( M\right) $ is characterized by the
following equalities: 
\begin{equation*}
\lambda _{IJk}\left( p\right) =\eta _{IJk}\left( p;U,a,\rho ,r\right) ,\quad
\left| I\right| ,\left| J\right| \geq 2,k\geq 1.
\end{equation*}
Then, since the algebraic subset of $\left( U,a,\rho ,r\right) $ is
characterized by finitely many equalities, there is an integer $K$ such that
the following equalities 
\begin{equation}
\lambda _{IJk}\left( p\right) =\eta _{IJk}\left( p;U,a,\rho ,r\right) ,\quad
\left| I\right| ,\left| J\right| ,k\leq K  \label{variety}
\end{equation}
characterize the element $\left( U_{\phi },a_{\phi },\rho _{\phi },r_{\phi
}\right) $ for $\phi \in Aut_{p}\left( M\right) .$

Since the isotropy subgroup $Aut_{p}\left( M\right) $ is compact, there is a
real number $e_{p}\geq 1$ such that 
\begin{equation*}
e_{p}^{-1}\leq \left| U_{\phi }\right| \leq e_{p},\quad \left| a_{\phi
}\right| \leq e_{p},\quad e_{p}^{-1}\leq \left| \rho _{\phi }\right| \leq
e_{p},\quad \left| r_{\phi }\right| \leq e_{p}
\end{equation*}
for every element $\phi \in Aut_{p}\left( M\right) .$ Thus the algebraic set
defined by the equalities in \ref{variety} is bounded. Note that the
boundedness of the algebraic set is preserved on an open neighborhood of the
point $p$ so that there are a real number $\delta _{p}>0$ and a real number $%
c_{p}\geq 1$ satisfying 
\begin{equation*}
c_{p}^{-1}\leq \left| U_{\phi }\right| \leq c_{p},\quad \left| a_{\phi
}\right| \leq c_{p},\quad c_{p}^{-1}\leq \left| \rho _{\phi }\right| \leq
c_{p},\quad \left| r_{\phi }\right| \leq c_{p}
\end{equation*}
for every point $p\in B\left( p;\delta _{p}\right) \cap M$ and every local
automorphism $\phi \in Aut_{p}(M).$

Since the subset $K$ is compact, there are finitely many points $%
p_{1},\cdots ,p_{l}$ such that 
\begin{equation*}
K\subset \bigcup_{k=1}^{l}B\left( p_{k};\delta _{p_{k}}\right) \cap M.
\end{equation*}
Then we take 
\begin{equation*}
e=\max \left\{ e_{p_{k}}:1\leq k\leq l\right\} .
\end{equation*}
This completes the proof.%
\endproof%

\begin{lemma}
\label{compact}Let $M$ be nonspherical analytic real hypersurfaces in $\Bbb{C%
}^{n+1}$ such that the isotropy subgroup $Aut_{p}\left( M\right) $ is
compact at every point $p\in M.$ Then, for each compact subset $K\subset
\subset M,$ there is a real number $\delta >0$ such that each germ of a
local automorphism $\phi \in Aut_{p}\left( M\right) ,$ $p\in K,$ is
analytically continued to the open ball $B(p;\delta ).$
\end{lemma}

\proof%
By the construction of the normalizing map $\mu _{p},$ the real number $%
\delta _{p}$ depend only on the point $p\in M.$ Hence there is a real number 
$\delta _{p}>0$ for each $p\in M$ such that, for every $q\in B\left(
p;\delta _{p}\right) \cap M,$ the mapping $\mu _{q}$ is biholomorphic on $%
B\left( q;\delta _{p}\right) $ and the inverse $\mu _{q}^{-1}$ is
biholomorphic on $B\left( 0;\delta _{p}\right) .$ Since $K$ is a compact
subset, there are finitely many points $p_{1},\cdots ,p_{l}$ such that 
\begin{equation*}
K\subset \bigcup_{k=1}^{l}B\left( p_{k};\delta _{p_{k}}\right) \cap M.
\end{equation*}
Then we take 
\begin{equation*}
\delta _{1}=\max \left\{ \delta _{p_{k}}:1\leq k\leq l\right\}
\end{equation*}
so that $\mu _{p}$ is biholomorphic on $B\left( p;\delta _{1}\right) $ and $%
\mu _{p}^{-1}$ is biholomorphic on $B\left( 0;\delta _{1}\right) $ for every 
$p\in K.$

By Lemma \ref{uniform1}, there is a real number $\delta _{2}$ such that
every local automorphism $\varphi \in Aut_{0}\left( \mu _{p}\left( M\right)
\right) $ for every $p\in K$ is biholomorphically continued to the
neighborhood $B\left( 0;\delta _{2}\right) $. Then we take 
\begin{equation*}
\delta =\min \left\{ \delta _{1},\delta _{2}\right\} .
\end{equation*}
This completes the proof.%
\endproof%

\begin{lemma}
\label{continuation}Let $M$ be a nonspherical analytic real hypersurface and 
$\gamma :[0,1]\rightarrow M$ be a chain-segment. Let $M^{\prime }$ be a
nonspherical analytic real hypersurface in Moser-Vitushkin normal form(cf. 
\cite{Pa1}). Suppose that there is a connected open neighborhood $U$ of the
point $\gamma \left( 0\right) $ and a biholomorphic mapping $\phi $ on $U$
such that 
\begin{equation*}
\phi \left( U\cap M\right) \subset M^{\prime }
\end{equation*}
and the image $\phi \left( U\cap \gamma [0,1]\right) $ is on the
straightened chain of $M^{\prime }.$ Then the mapping $\phi $ is
analytically continued along the whole chain-segment $\gamma [0,1]$ as a
local biholomorphic mapping.
\end{lemma}

\proof%
Suppose that the assertion is not true. Then there is a real number $\lambda
,$ $0<\lambda \leq 1,$ such that the mapping $\phi $ is analytically
continued along all the subpath $\gamma [0,\tau ],\tau <\lambda ,$ but not
the whole path $\gamma [0,\lambda ],$ as a local biholomorphic mapping$.$

Since $\gamma [0,1]$ is a chain-segment, there is a biholomorphic mapping $%
\varphi $ on an open neighborhood $V$ of the point $\gamma \left( \lambda
\right) $ such that the image $\varphi \left( V\cap \gamma [0,1]\right) $ is
on the straightened chain of $M^{\prime }.$ Then we take a point $p\in
\varphi \left( V\cap \gamma [0,\lambda )\right) $ and an open neighborhood $%
W $ of the point $p$ such that the mapping $\phi $ is analytically continued
on $W$ along the chain-segment $\gamma [0,1]$ and $W\subset V$ so that 
\begin{equation*}
\phi \circ \varphi ^{-1}\left( \varphi \left( W\right) \cap M^{\prime
}\right) \subset M^{\prime }
\end{equation*}
and the mapping $\phi \circ \varphi ^{-1}$ maps the straightened chain $%
\varphi \left( W\cap \gamma [0,1]\right) $ of $M^{\prime }$ onto the
straightened chain of $M^{\prime }.$

Note that the composition $\psi =\phi \circ \varphi ^{-1}$ is necessarily
analytically continued along the whole straightened chain of $M^{\prime }$%
(cf. \cite{Pa3}). Then, by abuse of notation, the composition $\psi \circ
\varphi $ is the analytic continuation of the mapping $\phi $ on the point $%
\gamma \left( \lambda \right) $ along the chain-segment $\gamma [0,1].$ This
is a contradiction. This completes the proof.%
\endproof%

\subsection{Analytic continuation on a real hypersurface}

We defined a canonical normalization $\nu _{p}$ of a nondegenerate analytic
real hypersurface $M$ at a point $p\in M$ to Moser-Vitushkin normal form(cf. 
\cite{Pa1}) by the same way of the canonical normalization $\mu _{p}$ to
Moser normal form.

\begin{lemma}
\label{chainsegment}Let $M,M^{\prime }$ be nonspherical analytic real
hypersurfaces in complex manifolds such that $M^{\prime }$ is compact and
the isotropy subgroup $Aut_{p}\left( M\right) $ is compact at every point $%
p\in M.$ Let $\gamma :[0,1]\rightarrow M$ be a chain-segment and $\phi $ be
a biholomorphic mapping on a connected open set $U$ of the point $\gamma
\left( 0\right) $ such that 
\begin{equation*}
\phi \left( U\cap M\right) \subset M^{\prime }.
\end{equation*}
Then the mapping $\phi $ is analytically continued along the whole
chain-segment $\gamma [0,1]$ as a local biholomorphic mapping.
\end{lemma}

\proof%
Suppose that the assertion is not true. Then, without loss of generality, we
may assume that the mapping $\phi $ is analytically continued along all the
subpath $\gamma [0,\tau ],\tau <1,$ but not the whole chain-segment $\gamma
[0,1]$ as a local biholomorphic mapping. By the analytic continuation of the
mapping $\phi ,$ there is a path $\phi \circ \gamma :[0,1)\rightarrow
M^{\prime }.$ Then we consider the following sequences 
\begin{eqnarray*}
p_{j} &\equiv &\gamma \left( 1-\frac{1}{j}\right) ,\quad j\in \Bbb{N}^{+}, \\
p_{j}^{\prime } &\equiv &\phi \circ \gamma \left( 1-\frac{1}{j}\right)
,\quad j\in \Bbb{N}^{+}.
\end{eqnarray*}
Since $M^{\prime }$ is compact, passing to a subsequence, say $m_{j},$ there
is a point $p^{\prime }\in M^{\prime }$ such that $p_{m_{j}}^{\prime
}\rightarrow p^{\prime }$ as $j\rightarrow \infty $ so that the closure 
\begin{equation*}
\overline{\left\{ p_{m_{j}}^{\prime }\equiv \phi \circ \gamma \left( 1-\frac{%
1}{m_{j}}\right) ,\quad j\in \Bbb{N}^{+}\right\} }
\end{equation*}
is a compact subset in a coordinate chart of the complex manifold.

Without loss of generality, we may assume that the chain $\gamma
[0,1]\subset M$ in a coordinate chart of the complex manifold. By abuse of
notation, we assume that $M,M^{\prime }$ are in $\Bbb{C}^{n+1}.$ Let $\nu
_{p}$ be the canonical normalization of a nondegenerate real hypersurface $M$
at a point $p\in M$ to Moser-Vitushkin normal form.

Then the mapping $\phi $ yields the following germs of a biholomorphic
mapping:

\begin{equation*}
\varphi _{j}\equiv \nu _{p_{m_{j}}^{\prime }}\circ \phi \circ \nu
_{p_{m_{j}}}^{-1}:\nu _{p_{m_{j}}}\left( M\right) \rightarrow \nu
_{p_{m_{j}}^{\prime }}\left( M^{\prime }\right) .
\end{equation*}
Since $\gamma :[0,1]\rightarrow M$ is a chain-segment, there is a sequence
of local automorphisms 
\begin{equation*}
\sigma _{j}\in Aut_{p_{m_{j}}}\left( M\right)
\end{equation*}
such that the mapping $\nu _{p_{m_{j}}^{\prime }}\circ \phi \circ \sigma
_{j} $ sends the germ of the chain $\gamma $ at the point $p_{m_{j}}$ onto
the straightened chain of the real hypersurface $\mu _{p_{m_{j}}^{\prime
}}\left( M^{\prime }\right) $ in Moser-Vitushkin normal form. By Lemma \ref
{continuation}, the composition $\psi _{j}=\nu _{p_{m_{j}}^{\prime }}\circ
\phi \circ \sigma _{j}$ is analytically continued along the whole
chain-segment $\gamma [0,1].$ Therefore, there is a real number $\delta >0$
such that, by abuse of notation, the mapping $\psi _{j}$ is biholomorphic on
the open neighborhood $B\left( p_{m_{j}};\delta \right) .$

By Lemma \ref{compact}, the set $\left\{ \nu _{p_{m_{j}}}\circ \sigma
_{j}\circ \nu _{p_{m_{j}}}^{-1}:j\in \Bbb{N}^{+}\right\} $ is relatively
compact.. Thus, by passing to a subsequence and shrinking $\delta >0$, if
necessary, we may assume that the composition $\tau _{j}=\nu
_{p_{m_{j}}}\circ \sigma _{j}\circ \nu _{p_{m_{j}}}^{-1}$ and its converse $%
\tau _{j}^{-1}$ are biholomorphically continued on the open neighborhood $%
B\left( 0;\delta \right) .$

Since $p_{m_{j}}\rightarrow \gamma \left( 1\right) $ and $p_{m_{j}}^{\prime
}\rightarrow p^{\prime },$ by passing to a subsequence and shrinking $\delta
>0$, if necessary, we may assume that the canonical normalizations $\nu
_{p_{m_{j}}},\nu _{p_{m_{j}}^{\prime }}$ and their inverses $\nu
_{p_{m_{j}}}^{-1},\nu _{p_{m_{j}}^{\prime }}^{-1}$ are biholomorphically
continued respectively on the open neighborhoods 
\begin{equation*}
B\left( p_{m_{j}};\delta \right) ,\quad B\left( p_{m_{j}}^{\prime };\delta
\right) ,\quad B\left( 0;\delta \right) ,\quad B\left( 0;\delta \right) .
\end{equation*}
Then, by abuse of notation, the composition 
\begin{equation*}
\chi _{j}=\nu _{p_{m_{j}}^{\prime }}^{-1}\circ \psi _{j}\circ \nu
_{p_{m_{j}}}^{-1}\circ \tau _{j}^{-1}\circ \nu _{p_{m_{j}}}
\end{equation*}
is biholomorphic on the open neighborhood $B\left( p_{m_{j}};\delta \right)
, $ if necessary, by shrinking $\delta >0.$

Note that the mapping $\chi _{j}$ is a local biholomorphic continuation of
the germ of the mapping $\phi $ at the point $p_{m_{j}}$ for each $j\in \Bbb{%
N}^{+}.$ Thus we take an integer $K$ such that 
\begin{equation*}
\gamma \left( 1\right) \in B\left( p_{m_{K}};\frac{\delta }{2}\right)
\end{equation*}
so that the mapping $\chi _{K}$ is an analytic continuation of the mapping $%
\phi $ on the point $\gamma \left( 1\right) $ along the chain-segment $%
\gamma [0,1].$ This contradiction completes the proof.%
\endproof%

\begin{lemma}
\label{alongpath}Let $M,M^{\prime }$ be nonspherical analytic real
hypersurfaces in complex manifolds such that $M^{\prime }$ is compact and
the isotropy subgroup $Aut_{p}\left( M\right) $ is compact at every point $%
p\in M.$ Suppose that there is a biholomorphic mapping $\phi $ on a
connected open set $U$ of a point $p\in M$ such that 
\begin{equation*}
\phi \left( U\cap M\right) \subset M^{\prime }.
\end{equation*}
Then the biholomorphic mapping $\phi $ is analytically continued along any
path on $M$ as a local biholomorphic mapping.
\end{lemma}

\proof%
For each point $p\in M,$ we make a biholomorphically equivalent deformation
of the real hypersurfaces $\mu _{p}\left( M\right) $ in normal form
continuously to a real hyperquadric by using the scaling mapping 
\begin{equation*}
\left\{ 
\begin{array}{l}
z^{*}=\lambda z \\ 
w^{*}=\lambda ^{2}w
\end{array}
\right. ,\quad \lambda \in \Bbb{R}^{+}.
\end{equation*}
Because the chain is characterized by an order differential equation, the
continuous family of chains on the real hyperquadric is continuously
deformed by the parameter $\lambda $ on a real hypersurface biholomorphic to 
$M$ near the point $p$(cf. \cite{Pa3}).

Let $\gamma :[0,1]\rightarrow M$ be a path on $M$ such that $\gamma \left(
0\right) \in U\cap M.$ Then, for each $\tau \in [0,1],$ there is a real
number $\varepsilon _{\tau }>0,$ a point $p_{\tau }\in M$ and a continuous
function 
\begin{equation*}
\Gamma _{\tau }:[0,1]\times \left( [0,1]\cap (\tau -\varepsilon _{\tau
},\tau +\varepsilon _{\tau })\right) \rightarrow M
\end{equation*}
such that

\begin{enumerate}
\item  $\Gamma _{\tau }\left( \cdot ,\sigma \right) :[0,1]\rightarrow M$ is
a chain-segment for each $\sigma \in [0,1]\cap (\tau -\varepsilon _{\tau
},\tau +\varepsilon _{\tau }),$

\item  $\Gamma _{\tau }\left( 0,\sigma \right) =\gamma \left( \sigma \right) 
$ for each $\sigma \in [0,1]\cap (\tau -\varepsilon _{\tau },\tau
+\varepsilon _{\tau }),$

\item  $\Gamma _{\tau }\left( 1,\sigma \right) =p_{\tau }$ for all $\sigma
\in [0,1]\cap (\tau -\varepsilon _{\tau },\tau +\varepsilon _{\tau }).$
\end{enumerate}

Note that the family $\left\{ [0,1]\cap (\tau -\varepsilon _{\tau },\tau
+\varepsilon _{\tau }):\tau \in [0,1]\right\} $ is an open covering of the
compact set $[0,1].$ Thus there is a finite subcover 
\begin{equation*}
\left\{ \lbrack 0,1]\cap (\tau _{j}-\varepsilon _{\tau _{j}},\tau
_{j}+\varepsilon _{\tau _{j}}):\tau _{j}\in [0,1],\quad j=1,\cdots
,m\right\} .
\end{equation*}
Then, by Lemma \ref{chainsegment}, the biholomorphic mapping $\phi $ is
analytically continued along the whole path $\gamma [0,1]$ as a local
biholomorphic mapping. This completes the proof.%
\endproof%

\subsection{Holomorphic mapping on the boundary}

\begin{lemma}
\label{polynomial}Let $Q$ be a real hyperquadric defined by 
\begin{equation*}
v=\langle z,z\rangle \equiv z^{1}\overline{z}^{1}+\cdots +z^{n}\overline{z}%
^{n}
\end{equation*}
and $\phi $ be a polynomial mapping as follows: 
\begin{equation*}
\phi :\left\{ 
\begin{array}{l}
z^{*}=f(z,w) \\ 
w^{*}=g(z,w)
\end{array}
\right.
\end{equation*}
where, for $m\geq 2,$%
\begin{eqnarray}
f(\mu z,\mu ^{2}w) &=&\mu ^{m}f(z,w)  \notag \\
g(\mu z,\mu ^{2}w) &=&\mu ^{2m}g(z,w).  \label{condi}
\end{eqnarray}
Suppose that 
\begin{equation*}
\phi \left( Q\right) \subset Q.
\end{equation*}
Then $\phi \equiv 0.$
\end{lemma}

\proof%
The mapping $\phi =(f,g)$ yields the identity 
\begin{equation}
\Im g(z,u+i\langle z,z\rangle )=\langle f(z,u+i\langle z,z\rangle
),f(z,u+i\langle z,z\rangle )\rangle .  \label{identity}
\end{equation}
Suppose that $m$ is even, i.e., $m=2k\geq 2.$ We may consider $z,\overline{z}%
,u$ as independent variables in the identity \ref{identity}. Taking $%
\overline{z}=0$ in the equality \ref{identity} yields 
\begin{equation*}
g(z,u)-\overline{g}(0,u)=2i\langle f(z,u),f(0,u)\rangle .
\end{equation*}
Thus we can put 
\begin{equation}
g(z,u)=2i\langle f(z,u),f(0,u)\rangle +\Re g(0,u)-i\langle
f(0,u),f(0,u)\rangle .  \label{identity2}
\end{equation}
Note that 
\begin{equation*}
f(0,u)=f(0,1)u^{\frac{m}{2}},\quad g(0,u)=g(0,1)u^{m}.
\end{equation*}
Then the identity \ref{identity} yields 
\begin{eqnarray*}
&&2\Re \left\{ \langle f(z,u+i\langle z,z\rangle ),f(0,1)\rangle (u+i\langle
z,z\rangle )^{\frac{m}{2}}\right\} \\
&&+\Re g(0,1)\Im (u+i\langle z,z\rangle )^{m} \\
&&-\langle f(0,1),f(0,1)\rangle \Re (u+i\langle z,z\rangle )^{m} \\
&=&\langle f(z,u+i\langle z,z\rangle ),f(z,u+i\langle z,z\rangle )\rangle .
\end{eqnarray*}

Let $p\left( z,\overline{z}\right) $ be a polynomial of the variables $z,%
\overline{z}$ satisfying 
\begin{equation*}
p\left( \mu z,\nu \overline{z}\right) =\mu ^{l}\nu ^{m}p\left( z,\overline{z}%
\right) .
\end{equation*}
Then the polynomial $p\left( z,\overline{z}\right) $ is said to be of type $%
(l,m).$ Collecting the terms of type $(l,1),$ $l=0,1,\cdots ,m-1,$ yields 
\begin{equation*}
miu^{-1}\langle z,z\rangle \left\{ \langle f(z,u),f(0,u)\rangle -\langle
f(0,u),f(0,u)\rangle +\Re g(0,u)\right\} =0.
\end{equation*}
Hence we obtain 
\begin{equation*}
\langle f(z,u),f(0,u)\rangle -\langle f(0,u),f(0,u)\rangle +\Re g(0,u)=0,
\end{equation*}
so that 
\begin{eqnarray*}
\langle f(z,u),f(0,u)\rangle &=&\langle f(0,u),f(0,u)\rangle \\
\Re g(0,u) &=&0.
\end{eqnarray*}
Thus, from the identity \ref{identity2}, we obtain 
\begin{equation*}
g(z,u)=i\langle f(0,u),f(0,u)\rangle
\end{equation*}
by which the identity \ref{identity} yields 
\begin{equation*}
\langle f(0,1),f(0,1)\rangle \Re (u+i\langle z,z\rangle )^{m}=\langle
f(z,u+i\langle z,z\rangle ),f(z,u+i\langle z,z\rangle )\rangle .
\end{equation*}
Note that 
\begin{equation*}
\langle f(z,u+i\langle z,z\rangle ),f(z,u+i\langle z,z\rangle )\rangle \geq
0,
\end{equation*}
but 
\begin{equation*}
\Re (u+i\langle z,z\rangle )^{m}=u^{m}-\frac{m(m-1)}{2}u^{m-2}\langle
z,z\rangle ^{2}+\cdots .
\end{equation*}
Thus collecting terms of type $(2,2)$ yields 
\begin{eqnarray*}
&&-\frac{m(m-1)}{2}\langle f(0,1),f(0,1)\rangle u^{m-2}\langle z,z\rangle
^{2} \\
&=&\left| \sum_{\alpha ,\beta }\frac{z^{\alpha }z^{\beta }}{2}\left( \frac{%
\partial ^{2}f}{\partial z^{\alpha }\partial z^{\beta }}\right) (0,u)\right|
^{2}+\frac{mu^{-2}}{2}\langle z,z\rangle ^{2}\langle f(0,u),f(0,u)\rangle \\
&\geq &0.
\end{eqnarray*}
Since $m\geq 2,$ we obtain 
\begin{equation*}
\langle f(0,1),f(0,1)\rangle =0
\end{equation*}
which yields 
\begin{equation*}
\langle f(z,u+i\langle z,z\rangle ),f(z,u+i\langle z,z\rangle )\rangle =0,
\end{equation*}
i.e., 
\begin{equation*}
f(z,w)=g(z,w)=0.
\end{equation*}

Suppose that $m$ is odd, i.e., $m=2k+1\geq 3.$ We may consider $z,\overline{z%
},u$ as independent variables in the identity \ref{identity}. Then, by the
condition \ref{condi}, taking $\overline{z}=0$ in the equality \ref{identity}
yields 
\begin{equation*}
g(z,u)=\overline{g}(0,u)
\end{equation*}
since 
\begin{equation*}
f(0,u)=0.
\end{equation*}
Thus we can put 
\begin{equation*}
g(z,u)=g(0,u),\quad g(0,1)\in \Bbb{R}.
\end{equation*}
Then the identity \ref{identity} yields 
\begin{equation}
g(0,1)\Im (u+i\langle z,z\rangle )^{m}=\langle f(z,u+i\langle z,z\rangle
),f(z,u+i\langle z,z\rangle )\rangle .  \label{identity3}
\end{equation}
Note that 
\begin{equation*}
\langle f(z,u+i\langle z,z\rangle ),f(z,u+i\langle z,z\rangle )\rangle \geq
0,
\end{equation*}
but 
\begin{eqnarray*}
\Im (u+i\langle z,z\rangle )^{m} &=&mu^{m-1}\langle z,z\rangle \\
&&-\frac{m(m-1)(m-2)}{2}u^{m-3}\langle z,z\rangle ^{3}+\cdots .
\end{eqnarray*}
Thus collecting terms of type $(1,1)$ yields 
\begin{eqnarray*}
&&mg(0,1)u^{m-1}\langle z,z\rangle \\
&=&\left| \sum_{\alpha }z^{\alpha }\left( \frac{\partial f}{\partial
z^{\alpha }}\right) (0,u)\right| ^{2}\geq 0.
\end{eqnarray*}
Collecting terms of type $(3,3)$ yields 
\begin{eqnarray*}
&&-\frac{m(m-1)(m-2)}{2}g(0,1)u^{m-3}\langle z,z\rangle ^{3} \\
&=&\left| \sum_{\alpha ,\beta ,\gamma }\frac{z^{\alpha }z^{\beta }z^{\gamma }%
}{6}\left( \frac{\partial ^{3}f}{\partial z^{\alpha }\partial z^{\beta
}\partial z^{\gamma }}\right) (0,u)\right| ^{2} \\
&&+\frac{m-1}{2}\left| \sum_{\alpha }z^{\alpha }\left( \frac{\partial f}{%
\partial z^{\alpha }}\right) (0,u)\right| ^{2}u^{-2}\langle z,z\rangle ^{2}
\\
&\geq &0.
\end{eqnarray*}
Since $m\geq 3,$ we obtain 
\begin{equation*}
g(0,1)\langle z,z\rangle =0
\end{equation*}
by which the identity \ref{identity3} yields 
\begin{equation*}
\langle f(z,u+i\langle z,z\rangle ),f(z,u+i\langle z,z\rangle )\rangle =0,
\end{equation*}
i.e., 
\begin{equation*}
f(z,w)=g(z,w)=0.
\end{equation*}
This completes the proof.%
\endproof%

\begin{lemma}
\label{nonconstant}Let $M,M^{\prime }$ be strongly pseudoconvex analytic
real hypersurfaces in $\Bbb{C}^{n+1}.$ Let $\phi $ be a holomorphic mapping
on an open neighborhood $U$ of a point $p\in M$ such that 
\begin{equation*}
\phi (U\cap M)\subset M^{\prime }.
\end{equation*}
Then the mapping $\phi $ is either a constant mapping or a biholomorphic
mapping on $U,$ if necessary, shrinking $U.$
\end{lemma}

\proof%
We take $q=\phi \left( p\right) $ so that 
\begin{equation*}
\varphi \equiv \mu _{q}\circ \phi \circ \mu _{p}^{-1}:\mu _{p}\left(
M\right) \rightarrow \mu _{q}\left( M^{\prime }\right) .
\end{equation*}
Then the mapping $\varphi $ is a holomorphic mapping on an open neighborhood 
$V$ of the origin satisfying 
\begin{equation*}
\varphi \left( \mu _{p}\left( M\right) \cap V\right) \subset \mu _{q}\left(
M^{\prime }\right) .
\end{equation*}

The mapping $\varphi =(f,g)$ in $\Bbb{C}^{n}\times \Bbb{C}$ is decomposed as
follows: 
\begin{equation*}
f(z,w)=\sum_{k=1}^{\infty }f_{k}(z,w),\mathrm{\quad }g(z,w)=\sum_{k=1}^{%
\infty }g_{k}(z,w),
\end{equation*}
where 
\begin{equation*}
f_{m}(\mu z,\mu ^{2}w)=\mu ^{m}f_{m}(z,w),\mathrm{\quad }g_{m}(\mu z,\mu
^{2}w)=\mu ^{m}g_{m}(z,w).
\end{equation*}
We assume that $\mu _{p}\left( M\right) ,\mu _{q}\left( M^{\prime }\right) $
are defined respectively by the equations 
\begin{equation*}
v=F(z,\overline{z},u),\mathrm{\quad }v=\langle z,z\rangle +F^{*}(z,\overline{%
z},u)
\end{equation*}
where 
\begin{equation*}
F(z,\overline{z},u)=\langle z,z\rangle +O\left( \left| z\right| ^{4}\right) ,%
\mathrm{\quad }F^{*}(z,\overline{z},u)=O\left( \left| z\right| ^{4}\right) .
\end{equation*}

Then we obtain the following identity near the origin 
\begin{eqnarray}
&&\Im g(z,u+iF(z,\overline{z},u))  \notag \\
&=&\langle f(z,u+iF(z,\overline{z},u)),f(z,u+iF(z,\overline{z},u))\rangle 
\notag \\
&&+F^{*}(f(z,u+iF(z,\overline{z},u)),\overline{f(z,u+iF(z,\overline{z},u))}%
,\Re g(z,u+iF(z,\overline{z},u))).  \label{iden}
\end{eqnarray}
Then, up to weight $2,$ we obtain 
\begin{align}
\Im g_{1}(z,0)& =0  \notag \\
\Im g_{2}(z,u+i\langle z,z\rangle )& =\langle f_{1}(z,0),f_{1}(z,0)\rangle 
\notag
\end{align}
which yields 
\begin{equation*}
\langle f_{1}(z,0),f_{1}(z,0)\rangle =\langle z,z\rangle \Re g_{2}(0,1)
\end{equation*}
and 
\begin{eqnarray*}
g_{1}(z,0) &=&0 \\
g_{2}(z,w) &=&g_{2}(0,1)w,\mathrm{\quad }g_{2}(0,1)\in \Bbb{R}.
\end{eqnarray*}
Note that $g_{2}(0,1)\neq 0$ if and only if the mapping $\varphi $ is a
biholomorphic mapping at the origin.

Suppose that the assertion is not true. Then we have 
\begin{equation*}
f_{1}(z,0)=g_{1}(z,0)=g_{2}(z,w)=0.
\end{equation*}
As inductive hypothesis, suppose that, for $m\geq 2,$%
\begin{eqnarray}
f_{l}(z,w) &=&0,\quad l=1,\cdots ,m-1,  \notag \\
g_{l}(z,w) &=&0,\quad l=1,\cdots ,2m-2.  \label{hypo}
\end{eqnarray}
By the condition 
\begin{equation*}
F(z,\overline{z},u)=\langle z,z\rangle +O\left( \left| z\right| ^{4}\right)
,\quad F^{*}(z,\overline{z},u)=o\left( \left| z\right| ^{4}\right) ,
\end{equation*}
the identity \ref{iden} yields 
\begin{equation*}
\Im g_{2m-1}(z,u+i\langle z,z\rangle )=0
\end{equation*}
Here we may consider $z,\overline{z},u$ as independent variables so that
taking $\overline{z}=0$ yields 
\begin{equation*}
g_{2m-1}(z,u+i\langle z,z\rangle )=\overline{g_{2m-1}}(0,u-i\langle
z,z\rangle )=0.
\end{equation*}
Hence the hypothesis \ref{hypo} necessarily comes to 
\begin{eqnarray*}
f_{l}(z,w) &=&0,\quad l=1,\cdots ,m-1, \\
g_{l}(z,w) &=&0,\quad l=1,\cdots ,2m-1.
\end{eqnarray*}
Then the identity \ref{iden} yields 
\begin{equation*}
\Im g_{2m}(z,u+i\langle z,z\rangle )=\langle f_{m}(z,u+i\langle z,z\rangle
),f_{m}(z,u+i\langle z,z\rangle )\rangle .
\end{equation*}
Note that the polynomial mapping $\varphi _{m}\equiv (f_{m},g_{2m}),$ $m\geq
2,$ satisfies 
\begin{equation*}
\varphi _{m}\left( Q\right) \subset Q.
\end{equation*}
By Lemma \ref{polynomial}, $\varphi _{m}\equiv 0$ so that 
\begin{eqnarray*}
f_{l}(z,w) &=&0,\quad l=1,\cdots ,m, \\
g_{l}(z,w) &=&0,\quad l=1,\cdots ,2m.
\end{eqnarray*}
This completes the induction so that $\varphi \equiv 0,$ i.e., the mapping $%
\varphi $ is a constant mapping. This completes the proof.%
\endproof%

\subsection{Proper holomorphic mappings}

\begin{lemma}
Let $D,D^{\prime }$ be strongly pseudoconvex bounded domains with
nonspherical real analytic boundaries $bD,bD^{\prime }$ such that the
boundaries $bD,bD^{\prime }$ are both simply connected. Suppose that there
is a biholomorphic mapping $\phi $ on a connected open neighborhood $U$ of a
point $p\in bD$ such that 
\begin{equation*}
\phi \left( U\cap bD\right) \subset bD^{\prime }.
\end{equation*}
Then the mapping $\phi $ is analytically continued to a biholomorphic
mapping from $D$ onto $D^{\prime }.$
\end{lemma}

\proof%
Since $bD,bD^{\prime }$ are simply connected, by Lemma \ref{alongpath}, the
mappings $\phi ,\phi ^{-1}$ are both analytically continued, by abuse of
notation, respectively to a biholomorphic mapping 
\begin{equation*}
\phi :D\rightarrow D
\end{equation*}
and 
\begin{equation*}
\phi ^{-1}:D^{\prime }\rightarrow D.
\end{equation*}
Then, by the identity theorem, the mapping $\phi $ is biholomorphic. This
completes the proof.%
\endproof%

\begin{lemma}
Let $D,D^{\prime }$ be strongly pseudoconvex bounded domains with
nonspherical real analytic boundaries $bD,bD^{\prime }$ such the boundary $%
bD^{\prime }$ is simply connected and the closed set $\overline{D^{\prime }}$
satisfies the fixed point property. Suppose that there is a biholomorphic
mapping $\phi $ on a connected open neighborhood $U$ of a point $p\in bD$
such that 
\begin{equation*}
\phi \left( U\cap bD\right) \subset bD^{\prime }.
\end{equation*}
Then the mapping $\phi $ is analytically continued to a biholomorphic
mapping from $D$ onto $D^{\prime }.$
\end{lemma}

\proof%
Since $bD^{\prime }$ is simply connected, by Lemma \ref{alongpath}, the
inverse mapping $\phi ^{-1}$ is analytically continued to a locally
biholomorphic proper mapping 
\begin{equation*}
\varphi :D^{\prime }\rightarrow D.
\end{equation*}
Note that the mappings $\varphi :D^{\prime }\rightarrow D$ and $\varphi
:bD^{\prime }\rightarrow bD$ are both covering maps. We claim that the
covering is simple. Otherwise, there is a nontrivial deck transformation of $%
\overline{D^{\prime }}$ yields a continuous self mapping of $\overline{%
D^{\prime }}$ with no fixed point. This is a contradiction. Therefore, the
mapping $\varphi :D^{\prime }\rightarrow D$ is biholomorphic. This completes
the proof.%
\endproof%

\begin{lemma}
\label{nonspheric}Let $D,D^{\prime }$ be strongly pseudoconvex bounded
domains with nonspherical real analytic boundaries $bD,bD^{\prime }$ such
that the domain $D$ and the boundary $bD^{\prime }$ are both simply
connected. Suppose that there is a biholomorphic mapping $\phi $ on a
connected open neighborhood $U$ of a point $p\in bD$ such that 
\begin{equation*}
\phi \left( U\cap bD\right) \subset bD^{\prime }.
\end{equation*}
Then the mapping $\phi $ is analytically continued to a biholomorphic
mapping from $D$ onto $D^{\prime }.$
\end{lemma}

\proof%
Since $bD^{\prime }$ is simply connected, by Lemma \ref{alongpath}, the
inverse mapping $\phi ^{-1}$ is analytically continued to a locally
biholomorphic proper mapping 
\begin{equation*}
\varphi :D^{\prime }\rightarrow D.
\end{equation*}
Note that the mapping $\varphi :D^{\prime }\rightarrow D$ is a covering map.
Since $D$ is simply connected, the covering is simple. Therefore, the
mapping $\varphi :D^{\prime }\rightarrow D$ is biholomorphic. This completes
the proof.%
\endproof%

\begin{theorem}
\label{convex}Let $D,D^{\prime }$ be strongly convex bounded domains with
real analytic boundaries $bD,bD^{\prime }$. Suppose that there is a
biholomorphic mapping $\phi $ on a connected open neighborhood $U$ of a
point $p\in bD$ such that 
\begin{equation*}
\phi \left( U\cap bD\right) \subset bD^{\prime }.
\end{equation*}
Then the mapping $\phi $ is analytically continued to a biholomorphic
mapping from $D$ onto $D^{\prime }.$
\end{theorem}

\proof%
Note that a strongly convex bounded domain is homeomorphic to an open ball $%
B^{n+1}.$ Suppose that the boundaries $bD,bD^{\prime }$ are spherical. Then
we take a biholomorphic mapping $\varphi $ on $U,$ if necessary, shrinking $%
U,$ such that 
\begin{equation*}
\varphi \left( U\cap bD\right) \subset bB^{n+1}.
\end{equation*}
Then, by Lemma \ref{boundary}, the mapping $\varphi $ and the composition $%
\psi \equiv \varphi \circ \phi ^{-1}$ are analytically continued, by abuse
of notation, to biholomorphic mappings as follows: 
\begin{equation*}
\varphi :D\rightarrow B^{n+1}\quad \text{and}\quad \psi :D^{\prime
}\rightarrow B^{n+1}.
\end{equation*}
Thus the composition $\psi ^{-1}\circ \varphi :D\rightarrow D^{\prime }$ is
a biholomorphic mapping and the analytic continuation of the mapping $\phi .$

Suppose that the boundaries $bD,bD^{\prime }$ are nonspherical. Then, by
Lemma \ref{nonspheric}, the mapping $\phi $ is analytically continued to a
biholomorphic mapping from $D$ onto $D^{\prime }.$ This completes the proof.%
\endproof%

\begin{lemma}
\label{proper}Let $D,D^{\prime }$ be strongly pseudoconvex bounded domains
with real analytic boundaries $bD,bD^{\prime }.$ Suppose that there is a
proper holomorphic mapping $\phi :D\rightarrow D^{\prime }.$ Then there is
an open neighborhood $U$ of a point $p\in bD$ such that the mapping $\phi $
is analytically continued on $U$ and 
\begin{equation*}
\phi \left( U\cap bD\right) \subset bD^{\prime }.
\end{equation*}
\end{lemma}

\proof%
We may apply the boundary regularity of Lemma \ref{b-regularity} so that the
mapping $\varphi :D\rightarrow D^{\prime }$ is holomorphic on an open
neighborhood of $\overline{D}$. This completes the proof.%
\endproof%

\begin{lemma}
\label{Alexander}Let $\phi :B^{n+1}\rightarrow B^{n+1}$ be a proper
holomorphic mapping. Then $\phi \in Aut\left( B^{n+1}\right) .$
\end{lemma}

\proof%
By Lemma \ref{proper}, there is an open neighborhood $U$ of a point $p\in
bB^{n+1}$ such that $\phi $ is analytically continued on $U$ and 
\begin{equation*}
\phi \left( U\cap bB^{n+1}\right) \subset bB^{n+1}.
\end{equation*}
By Lemma \ref{sphere}, $\phi \in Aut\left( B^{n+1}\right) .$ This completes
the proof.%
\endproof%

\begin{lemma}
\label{equiv}Let $D,D^{\prime }$ be strongly pseudoconvex bounded domains in 
$\Bbb{C}^{n+1}$ with real analytic boundaries $bD,bD^{\prime }.$ Suppose
that there is a proper holomorphic mapping $\phi :D\rightarrow D^{\prime }$.
Then the mapping $\phi :D\rightarrow D^{\prime }$ is a locally biholomorphic
mapping so that $\phi :\overline{D}\rightarrow \overline{D^{\prime }}$ is a
covering map.
\end{lemma}

\proof%
By Lemma \ref{proper}, there is a point $p\in bD$ and an open neighborhood $%
U $ of the point $p$ such that the mapping $\phi $ is analytically continued
to $U$ and 
\begin{equation*}
\phi \left( U\cap bD\right) \subset bD^{\prime }.
\end{equation*}
By Lemma \ref{nonconstant}, the mapping $\phi $ is biholomorphic on $U,$ if
necessary, shrinking $U.$

Suppose that the boundaries $bD,bD^{\prime }$ are spherical. Then, by Lemma 
\ref{sph-proper}, the mapping $\phi :D\rightarrow D^{\prime }$ is a locally
biholomorphic mapping.

Suppose that the boundary $bD,bD^{\prime }$ are nonspherical. Then, by Lemma 
\ref{alongpath}, the mapping $\phi :D\rightarrow D^{\prime }$ is
analytically continued on an open neighborhood of the boundary $bD$ to be
locally biholomorphic. Thus the mapping $\phi $ is a locally biholomorphic
mapping. This completes the proof.%
\endproof%

\begin{theorem}
Let $D$ be a strongly pseudoconvex bounded domain with real analytic
boundary $bD.$ Suppose that there is a proper holomorphic self mapping $\phi
:D\rightarrow D.$ Then $\phi $ is a biholomorphic automorphism of $D$.
\end{theorem}

\proof%
By Lemma \ref{equiv}, the mapping $\phi :\overline{D}\rightarrow \overline{D}
$ is a self covering map. We claim that $\phi :D\rightarrow D$ is a simple
covering. Otherwise, there would be a integer $m>1$ such that, for every $%
p\in \overline{D},$ $m$ is the order of the set 
\begin{equation*}
\left\{ q\in \overline{D}:\phi \left( q\right) =p\right\} .
\end{equation*}
Then we define the $k$ times composition $\phi ^{k}$ of the mapping $\phi $
as follows 
\begin{equation*}
\phi ^{k}\equiv \underbrace{\phi \circ \cdots \circ \phi }_{k}:\overline{D}%
\rightarrow \overline{D}.
\end{equation*}
Note that the inverse image of each point $p\in \overline{D}$ under the
mapping $\phi ^{k}:\overline{D}\rightarrow \overline{D},$ $k\in \Bbb{N}^{+},$
is a set of order $m^{k}.$

We claim that, for a given real number $\delta >0,$ there is a point $q\in
bD $ and a compact subset $K\subset \subset D$ and a subsequence $\phi
^{m_{k}}$ such that 
\begin{equation*}
\phi ^{-m_{k}}\left( B\left( q;\delta \right) \cap D\right) \cap K\neq
\emptyset \quad \text{for all }k\in \Bbb{N}^{+}.
\end{equation*}
Otherwise, for every compact subset $K\subset \subset D,$ there is an
integer $l$ such that 
\begin{equation*}
\phi ^{-k}\left( D\backslash K\right) \subset D\backslash K\quad \text{for
all }k\geq l.
\end{equation*}
Since $\phi ^{k}:D\rightarrow D$ is a finite covering map, we have 
\begin{equation*}
D=\phi ^{-k}\left( K\cup D\backslash K\right) =\phi ^{-k}\left( K\right)
\cup \phi ^{-k}\left( D\backslash K\right)
\end{equation*}
so that 
\begin{equation*}
K\subset \phi ^{-k}\left( K\right) \quad \text{for all }k\geq l.
\end{equation*}
Note that the inverse image of $\phi ^{-k}\left( K\right) $ is a union of $%
m^{k}$ disconnected compact subsets. We may assume that $K$ is connected so
that $K$ is in a compact subset of the union $\phi ^{-k}\left( K\right) .$
This is impossible. Thus, with such a point $q\in bD$, we take an
accumulation point $q^{\prime }\in D^{\prime }$ of the set 
\begin{equation*}
\lim_{k\rightarrow \infty }\phi ^{-m_{k}}\left( \left\{ q\right\} \right)
=\lim_{k\rightarrow \infty }\left\{ p\in bD:\phi ^{m_{k}}\left( p\right)
=q\right\} .
\end{equation*}

Suppose that $bD$ is spherical. Then there are open neighborhoods $%
U,U^{\prime }$ respectively of $q,q^{\prime }$ and biholomorphic mappings $%
\psi ,\psi ^{\prime }$ respectively on $U,U^{\prime }$ such that 
\begin{eqnarray*}
\psi \left( U\cap bD\right) &\subset &bB^{n+1} \\
\psi ^{\prime }\left( U^{\prime }\cap bD\right) &\subset &bB^{n+1}.
\end{eqnarray*}
By Lemma \ref{sphere}, the compositions 
\begin{equation*}
\varphi _{k}\equiv \psi ^{\prime }\circ \phi ^{-k}\circ \psi ^{-1}:\psi
\left( U\cap D\right) \rightarrow \psi ^{\prime }\left( U^{\prime }\cap
D\right)
\end{equation*}
is analytically continued to, by abuse of notation, automorphisms $\varphi
_{k}\in Aut\left( B^{n+1}\right) .$ By the construction, a subsequence $%
\varphi _{m_{k}}$ must converges to the point $\psi ^{\prime }\left(
q^{\prime }\right) $ uniformly on every compact subset of $B^{n+1}$. Thus
there is a sequence of real numbers $\delta _{k}\searrow 0$ such that 
\begin{equation*}
\phi ^{m_{k}}:B\left( q^{\prime };\delta _{k}\right) \cap D\rightarrow U\cap
D
\end{equation*}
and 
\begin{equation*}
\phi ^{m_{k}}=\psi ^{-1}\circ \varphi _{m_{k}}^{-1}\circ \psi ^{\prime
}\quad \text{on }B\left( q^{\prime };\delta _{k}\right) \cap D.
\end{equation*}
Then we take a point $p\in U\cap D$ such that there is a sequence $q_{k}$
satisfying 
\begin{equation*}
q_{k}\in \phi ^{-m_{k}}\left( \left\{ p\right\} \right) \cap B\left(
q^{\prime };\delta _{k}\right)
\end{equation*}
and, for a compact subset $K\subset \subset D,$ there is a sequence $p_{k}$
satisfying 
\begin{equation*}
p_{k}\in \phi ^{-m_{k}}\left( \left\{ p\right\} \right) \cap K.
\end{equation*}

Therefore, there is a sequence of deck transformations $\phi _{k}\in
Aut\left( D\right) $ of the covering map $\phi ^{m_{k}}:D\rightarrow D$ such
that 
\begin{equation*}
\phi _{k}\left( p_{k}\right) =q_{k}\rightarrow q^{\prime }\in bD.
\end{equation*}
Since $p_{k}\in K,$ by Lemma \ref{Wong-Rosay}, there is a biholomorphic
mapping 
\begin{equation*}
\sigma :D\rightarrow B^{n+1}.
\end{equation*}
Then the composition 
\begin{equation*}
\iota \equiv \sigma \circ \phi \circ \sigma ^{-1}:B^{n+1}\rightarrow B^{n+1}
\end{equation*}
would be a proper self mapping, but not an automorphism of the unit ball $%
B^{n+1}.$ This is a contradiction to Lemma \ref{Alexander} that every proper
self mapping of the unit ball $B^{n+1}$ is an automorphism of the unit ball $%
B^{n+1}.$

Suppose that $bD$ is nonspherical. We take a sequence $p_{k}\in bD$ such
that 
\begin{equation*}
p_{k}\in \phi ^{-k}\left( \left\{ q\right\} \right) \subset bD
\end{equation*}
and 
\begin{equation*}
\left| p_{k}-q^{\prime }\right| =\min \left\{ \left| p-q^{\prime }\right|
:p\in \phi ^{-k}\left( \left\{ q\right\} \right) \right\} .
\end{equation*}
Let $\mu _{p_{k}}$ be the canonical normalizing mapping at the point $%
p_{k}\in bD.$ Since $p_{k}\rightarrow q^{\prime },$ there is an open
neighborhood $W$ of the point $q^{\prime }$ such that the mapping $\mu
_{p_{k}}$ is analytically continued on $W.$ Then we define a function $%
\varepsilon _{k}\left( p\right) $ for a sufficiently large $k$ and a point $%
p\in W\cap bD$ such that 
\begin{equation*}
\varepsilon _{k}\left( p\right) \equiv \sum_{j=1}^{n}\left| z_{j}\circ \mu
_{p_{k}}\left( p\right) \right| ^{2}+\left| w\circ \mu _{p_{k}}\left(
p\right) \right|
\end{equation*}
where $z_{j},$ $j=1,\cdots ,n,$ $w$ are the coordinate functions of $\Bbb{C}%
^{n+1}$ such that $z_{j},$ $j=1,\cdots ,n,$ are for the complex tangent
hyperplane and $w$ are for the complex line normal to the complex tangent
hyperplane of $\mu _{p_{k}}\left( W\cap bD\right) $ at the origin. Then we
take a sequence $q_{k}$ such that 
\begin{equation*}
q_{k}\in \phi ^{-k}\left( \left\{ q\right\} \right) \subset bD
\end{equation*}
and 
\begin{equation*}
\varepsilon _{k}\left( q_{k}\right) =\min \left\{ \varepsilon _{k}\left(
p\right) :p\in \phi ^{-k}\left( \left\{ q\right\} \right) \right\} .
\end{equation*}
Let $\pi _{k}$ be a complex line containing $p_{k}$ and $q_{k}.$ Then we
take a subsequence $\pi _{m_{k}}$ so that $\pi _{m_{k}}$ converges to a
complex line $\pi _{q^{\prime }}\subset T_{q^{\prime }}\Bbb{C}^{n+1}.$ In
other words, the a subsequence $q_{m_{k}}$ converges to the point $q^{\prime
}$ to a direction. Then we obtain 
\begin{equation*}
\varphi _{k}\equiv \mu _{p_{k}}\circ \phi ^{-k}\circ \mu _{q}^{-1}:\mu
_{q}\left( bD\right) \rightarrow \mu _{p_{k}}\left( bD\right) .
\end{equation*}
Then we define a biholomorphic mapping 
\begin{equation*}
\sigma _{k}:\left\{ 
\begin{array}{l}
z^{*}=\sqrt{\varepsilon _{m_{k}}}z \\ 
w^{*}=\varepsilon _{m_{k}}w
\end{array}
\right.
\end{equation*}
where 
\begin{equation*}
\varepsilon _{m_{k}}=\varepsilon _{m_{k}}\left( q_{m_{k}}\right) .
\end{equation*}
Then we take the composition 
\begin{equation*}
\kappa _{k}\equiv \sigma _{k}^{-1}\circ \varphi _{m_{k}}:\mu _{q}\left(
bD\right) \rightarrow \sigma _{k}^{-1}\circ \mu _{p_{m_{k}}}\left( bD\right)
.
\end{equation*}
Note that there is a real number $\delta >0$ such that the germs of real
hypersurfaces $\mu _{q}\left( bD\right) $ and $\sigma _{k}^{-1}\circ \mu
_{p_{m_{k}}}\left( bD\right) $ are analytically continued to the open
neighborhood $B\left( 0;\delta \right) .$ Since the boundary $bD$ is
nonspherical, by Lemma \ref{compact}, the mapping $\kappa _{k}$ is
biholomorphic on $B\left( 0;\delta \right) ,$ if necessary, shrinking $%
\delta >0.$ Further, by the construction, the limit of the sequence $\kappa
_{k}$ cannot be a constant mapping on $B\left( 0;\delta \right) $ so that,
by Theorem \ref{nonconstant}, the limit of the sequence $\kappa _{k}$ would
be biholomorphic on $B\left( 0;\delta \right) .$

By the way, the sequence of real hypersurfaces $\sigma _{k}^{-1}\circ \mu
_{p_{m_{k}}}\left( bD\right) $ converges to a real hyperquadric uniformly on
an open neighborhood of the origin. Therefore, the germ of the boundary $bD$
at the point $q^{\prime }$ is spherical so that the boundary $bD$ is
spherical(cf. \cite{Pa3}). This is a contradiction to the fact that the
boundary $bD$ is nonspherical. This completes the proof.%
\endproof%

\subsection{Locally realizable CR manifolds}

Let $M$ be a CR manifold of CR dimension $n$ and CR codimension $1$ with a
CR structure $(D,I)$ where $D$ is $2n$ dimensional smooth subbundle of the
tangent bundle $TM$ and $I$ is an automorphism on $D$ such that 
\begin{equation*}
I^{2}V=-V\quad \text{for}\quad V\in \Gamma D.
\end{equation*}
For each point $p\in M,$ there is a local coordinate chart $\left( U,\varphi
\right) $ such that 
\begin{equation*}
\varphi \left( U\right) \subset \Bbb{R}^{2n+1}.
\end{equation*}
Then $M$ shall be called locally realizable CR manifold if there is an open
neighborhood $U$ of each point $p\in M$ and CR functions $f_{1},\cdots
,f_{n+1}$ on $\varphi \left( U\right) $ satisfying 
\begin{equation*}
df_{1}\wedge \cdots \wedge df_{n+1}\neq 0.
\end{equation*}
Let $z_{j},$ $j=1,\cdots ,n+1$, be the coordinate functions of $\Bbb{C}%
^{n+1}.$ Then we obtain a local embedding $\sigma $ of $U\subset M$ into $%
\Bbb{C}^{n+1}$ such that 
\begin{equation*}
f_{j}\equiv z_{j}\circ \sigma \circ \varphi ^{-1}.
\end{equation*}
The CR manifold $M$ shall be called a locally analytically realizable CR
manifold when $\sigma \left( U\right) $ is real analytic by a local
embedding $\sigma $ on an open neighborhood $U$ of each point $p\in M.$ Note
that $M$ is either spherical or nonspherical(cf. \cite{Pa3}).

\begin{lemma}
Let $M$ be a connected locally analytically realizable CR manifold. Suppose
that there is a nontrivial CR mapping $\varphi $ on an open neighborhood $U$
of a point $p\in M$ such that 
\begin{equation*}
\varphi \left( U\right) \subset bB^{n+1}.
\end{equation*}
Then the mapping $\varphi $ is CR continued along any path on $M$ as a
locally CR diffeomorphic mapping.
\end{lemma}

\proof%
Since $M$ is connected, $M$ is necessarily spherical(cf. \cite{Pa3}).
Suppose that the assertion is not true. Then there is a path $\gamma
:[0,1]\rightarrow M$ with $\gamma \left( 0\right) =p$ such that the mapping $%
\varphi $ is CR continued along all subpath $\gamma [0,\tau ]$ with $\tau
<1, $ but not the whole path $\gamma [0,1].$ Then we take a CR embedding $%
\sigma \left( V\right) \subset \Bbb{C}^{n+1}$ of an open neighborhood $V$ of
the point $\gamma \left( 1\right) $ such that $\sigma \left( V\right) $ is a
spherical analytic real hypersurface and 
\begin{equation*}
\varphi \circ \sigma ^{-1}:\sigma \left( V\right) \rightarrow bB^{n+1}
\end{equation*}
is a locally CR diffeomorphism. Without loss of generality, we may assume
that $\gamma [0,1]\subset V.$ Hence there is an open neighborhood $W$ of the
point $\sigma \left( \gamma \left( 0\right) \right) $ and a biholomorphic
mapping $\phi $ on $W$ such that 
\begin{equation*}
\phi =\varphi \circ \sigma ^{-1}\quad \text{on }\sigma \left( V\right) \cap
W.
\end{equation*}
Then the mapping $\phi $ is analytically continued along the whole path $%
\sigma \left( \gamma [0,1]\right) .$ From, by abuse of notation, the mapping 
$\phi $ at an open neighborhood $W^{\prime }$ of the point $\sigma \left(
\gamma \left( 1\right) \right) $ and a CR embedding $\sigma $ of an open
neighborhood $V^{\prime }$ of the point $\gamma \left( 1\right) $, we obtain
the CR mapping 
\begin{equation*}
\phi \circ \sigma :\sigma ^{-1}\left( W^{\prime }\cap \sigma \left(
V^{\prime }\right) \right) \rightarrow bB^{n+1}
\end{equation*}
which is a CR continuation of the CR mapping $\varphi .$ This completes the
proof.%
\endproof%

\begin{theorem}
Let $M$ be a connected locally analytically realizable CR manifold. Suppose
that $M$ is compact and there is a nontrivial CR mapping $\varphi $ on an
open neighborhood $U$ of a point $p\in M$ such that 
\begin{equation*}
\varphi \left( U\right) \subset bB^{n+1}.
\end{equation*}
Then there is a finite subset $L\subset bB^{n+1}$ such that the inverse CR
mapping $\varphi ^{-1}$ is CR continued along any path on $%
bB^{n+1}\backslash L$ as a locally CR diffeomorphic mapping.
\end{theorem}

\proof%
Since $M$ is connected, $M$ is necessarily spherical(cf. \cite{Pa3}). Let $%
q_{j}\in bB^{n+1}$ be a sequence converging to a point $q\in bB^{n+1}$ and $%
\phi _{j}$ be a sequence of the CR continuation of the inverse mapping $%
\varphi ^{-1}$ at the point $q_{j}$. We set 
\begin{equation*}
q_{j}^{\prime }=\phi _{j}\left( q_{j}\right) \in M.
\end{equation*}
Since $M$ is compact, there is a subsequence $q_{m_{j}}^{\prime }$ and a
point $q^{\prime }\in M$ such that 
\begin{equation*}
q_{m_{j}}^{\prime }\rightarrow q^{\prime }.
\end{equation*}
Then we take a CR embedding $\sigma $ of an open neighborhood $V$ of the
point $q^{\prime }\in M$ such that $\sigma \left( V\right) \subset \Bbb{C}%
^{n+1}$ is a spherical analytic real hypersurface. Then we apply the same
argument as in the previous sections so that the singular locus $L$ is a
finite subset of $bB^{n+1}$ and the inverse CR mapping $\varphi ^{-1}$ is CR
continued along any path on $bB^{n+1}\backslash L.$ This completes the proof.%
\endproof%

\begin{lemma}
\label{cranypath}Let $M,M^{\prime }$ be connected locally analytically
realizable CR manifolds with positive definite Levi form. Suppose that $%
M,M^{\prime }$ are nonspherical and $M^{\prime }$ is compact, and there is a
nontrivial CR mapping $\varphi $ on an open neighborhood $U$ of a point $%
p\in M$ such that 
\begin{equation*}
\varphi \left( U\right) \subset M^{\prime }.
\end{equation*}
Then the mapping $\varphi $ is CR continued along any path on $M$ as a
locally CR diffeomorphic mapping.
\end{lemma}

\proof%
Since $M,M^{\prime }$ are locally analytically realizable and the Levi forms
of $M,M^{\prime }$ are positive definite, the isotropy subgroups $%
Aut_{p}\left( M\right) ,Aut_{p^{\prime }}\left( M^{\prime }\right) $ are
compact for every point $p\in M,p^{\prime }\in M^{\prime }.$

Let $q_{j}\in M$ be a sequence converging to a point $q\in M$ and $\varphi
_{j}$ be a sequence of the CR continuation of the mapping $\varphi $ at the
point $q_{j}$. We set 
\begin{equation*}
q_{j}^{\prime }=\varphi _{j}\left( q_{j}\right) \in M^{\prime }.
\end{equation*}
Since $M^{\prime }$ is compact, there is a subsequence $q_{m_{j}}^{\prime }$
and a point $q^{\prime }\in M^{\prime }$ such that 
\begin{equation*}
q_{m_{j}}^{\prime }\rightarrow q^{\prime }.
\end{equation*}
Then we take CR embeddings $\sigma ,\sigma ^{\prime }$ respectively of open
neighborhoods $V,V^{\prime }$ respectively of the points $q\in M,q^{\prime
}\in M^{\prime }$ such that $\sigma \left( V\right) ,\sigma ^{\prime }\left(
V^{\prime }\right) \subset \Bbb{C}^{n+1}$ are nonspherical analytic real
hypersurface. Then we apply the same argument as in the previous subsection.
This completes the proof.%
\endproof%

\begin{theorem}
Let $M,M^{\prime }$ be connected locally analytically realizable CR
manifolds with positive definite Levi form. Suppose that $M,M^{\prime }$ are
compact and nonspherical, and there is a nontrivial CR mapping $\varphi $ on
an open neighborhood $U$ of a point $p\in M$ such that 
\begin{equation*}
\varphi \left( U\right) \subset M^{\prime }.
\end{equation*}
Then the maximal CR extension of the mapping $\varphi $ is a CR equivalence
between the natural universal covering spaces of $M,M^{\prime }$ to be the
pointed path spaces respectively of $M,M^{\prime }$ mod homotopic relation.
\end{theorem}

\proof%
By Lemma \ref{cranypath}, the mapping $\varphi $ is CR continued along any
path on $M$ as a locally CR diffeomorphic mapping. Note that $M$ is compact
and 
\begin{equation*}
\varphi ^{-1}\left( \varphi \left( U\right) \right) \subset M
\end{equation*}
so that we apply Lemma \ref{cranypath} to the inverse mapping $\varphi
^{-1}. $ Thus the inverse mapping $\varphi ^{-1}$ is CR continued along any
path on $M^{\prime }$ as a locally CR diffeomorphic mapping. Thus the CR
continuation of the mapping $\phi $ induces a CR equivalence between the
natural universal coverings of $M,M^{\prime },$ which are the path spaces
mod homotopic relation respectively over $M,M^{\prime }$ with the natural CR
structure. This completes the proof.%
\endproof%

\end{document}